\documentclass[11pt,a4paper,parskip=half]{scrartcl}
\usepackage[utf8]{inputenc}
\usepackage[english]{babel}
\usepackage[T1]{fontenc}
\usepackage{lmodern}
\usepackage[left=3cm,right=3cm,top=3cm,bottom=3cm]{geometry}
\usepackage[runin]{abstract}
    \abslabeldelim{.}

\usepackage{todonotes}
\usepackage[toc,page]{appendix} 

\usepackage{amsmath, amsfonts, amssymb, amsthm, amstext}
\usepackage{mathtools}
\usepackage{mathrsfs}
\usepackage{latexsym}
\usepackage{mathdots}
\usepackage{cases}

\usepackage{graphicx}
\usepackage{subcaption}
\usepackage{algorithm}

 at11pt
 at9pt
\usepackage{enumitem}

\usepackage{tikz}

\usetikzlibrary{
    positioning,
    arrows.meta,
    shapes.geometric,
    fit,
    calc
}

\usepackage{xcolor}

\usepackage{color}
\usepackage{hyperref} %, nameref}
\usepackage[capitalise]{cleveref}
	
\newcounter{thm}
\numberwithin{thm}{section}
\numberwithin{equation}{section}

	\newtheoremstyle{myplain}		% name of the style to be used
			{}			% measure of space to leave above the theorem. E.g.: 3pt
			{}			% measure of space to leave below the theorem. E.g.: 3pt
			{\itshape}				% name of font to use in the body of the theorem
			{}				% measure of space to indent
			{\sffamily\bfseries}				% name of head font
			{.}		% punctuation between head and body
			{ }				% space after theorem head; " " = normal interword space
			{\thmname{#1}\thmnumber{ #2}\textnormal{\textsf{\thmnote{ (#3)}}}}			% Manually specify head
    \newtheoremstyle{mybreak}
            {}{}{}{}{\sffamily\bfseries}{.}{\newline}
            {\thmname{#1}\thmnumber{ #2}\textnormal{\textsf{\thmnote{ (#3)}}}}
	\newtheoremstyle{mydef}
			{}{}{}{}{\sffamily\bfseries}{.}{ }
			{\thmname{#1}\thmnumber{ #2}}
	\newtheoremstyle{myrem}
			{}{}{}{}{\sffamily\itshape}{.}{ }
			{\thmname{#1}\thmnumber{ #2}}

\theoremstyle{myplain}

\theoremstyle{mybreak}
\theoremstyle{mydef}

\theoremstyle{mydef}
	\newtheorem{example}[thm]{Example}

	\newcommand{\cc}{\mathbb{C}}
		
		\newcommand{\nn}{\mathbb{N}}
	\newcommand{\rr}{\mathbb{R}}

\allowdisplaybreaks

\def\sumprime_#1^#2{
    \setbox0=\hbox{$\scriptstyle{#1}$}
    \setbox1=\hbox{$\scriptstyle{#2}$}
    \setbox2=\hbox{$\displaystyle{\sum}$}
    \setbox4=\hbox{${}^\prime\mathsurround=0pt$}
    \dimen0=.5\wd0 \advance\dimen0 by-.5\wd2
    \ifdim\dimen0>0pt
        \ifdim\dimen0>\wd4 \kern\wd4
        \else\kern\dimen0
        \ifdim\dimen1>\wd4 \kern\wd4
        \else\kern\dimen1
    \fi\fi\fi
\mathop{{\sum}^\prime}_{\kern-\wd4 #1}^{\kern-\wd4 #2}
}
\title{\Large  Rational neural networks for tracking    complex singularities of nonlinear PDEs }

\author{
Nadiia Derevianko%
\thanks{Corresponding author. Email: nadiia.derevianko@tum.de} 
 \thanks{Department of Computer Science, School of Computation, Information and Technology, Technical University of Munich, Boltzmannstrasse 3, 85748 Garching b. München, Germany}
\and
Hans-Joachim Bungartz\footnotemark[2] 
 \thanks{Email: bungartz@cit.tum.de}
\and
Felix Dietrich\footnotemark[2] 
 \thanks{MDSI \& MCML, Technical University of Munich. Email: felix.dietrich@tum.de}
}

\date{\today}

\begin{document}

\maketitle

\begin{abstract}
We present a neural network-based method for  numerical analytic continuation of solutions of nonlinear partial differential equations (PDEs) and detection of their complex singularities. The proposed framework employs ``unsafe'' Padé activation units (PAUs) as activation functions, together with a novel backpropagation-free method for computing the weights and biases of the hidden layers. The training procedure is designed specifically for learning meromorphic  functions with pole-type singularities.  Unlike existing methods with fixed hidden-layer parameters, the proposed approach treats the weights and biases as time-dependent functions, allowing the network to adapt to the evolution of the complex singularities.
Using this method, we can efficiently locate the complex singularities of the extended solution, track their trajectories over time, and, based on this dynamics, infer the formation of characteristic phenomena in the solutions of nonlinear PDEs. 
To demonstrate the performance of our method, we analyze three well-known cases: (1) the nonlinear heat equation, which exhibits finite-time blow up; (2) the nonlinear Burgers equation, which develops a shock; and (3) the nonlinear Schrödinger equation, in which rogue waves form.

%\\[1ex]

\textbf{Mathematics Subject Classification (MSC2020):} 	30D30, 	30E10, 34M99, 41A20,  65D15, 	65M70
\end{abstract}

\section{Introduction}

Blow up, shock formation, and rogue waves are among the most fundamental nonlinear phenomena encountered in mathematical physics, engineering, and the applied sciences. They arise in a wide range of applications, including fluid dynamics \cite{MB2002}, optics and oceanography \cite{Dud2019}, biological systems \cite{Hor2003}, and bacterial communication \cite{JK2013}. Although these phenomena occur in different classes of nonlinear PDEs, they are all closely related to the evolution of complex singularities of the analytic continuations of the corresponding solutions \cite{FH19,K17,W03,W22}. In particular, the dynamics of complex singularities, such as poles and branch points, provide valuable information about the onset, location, and nature of the underlying nonlinear phenomena.
%Consequently, the analytic continuation of solutions into the complex plane has become an important tool for the analysis of nonlinear PDEs.
%In this paper, we present a novel neural network-based approach for constructing a \emph{numerical} analytic continuation of solutions to nonlinear PDEs. The proposed method enables the accurate detection and tracking of complex singularities directly from initial-condition data on the real axis, providing an effective framework for investigating the mechanisms underlying blow up, shock formation, rogue waves, and other singular phenomena.
%Consequently, they provide an ideal set of benchmark problems for assessing the capability of the proposed neural network-based numerical analytic continuation method.
%Many nonlinear   PDEs exhibit complex dynamical phenomena, including shock formation, finite-time blow up, and the emergence of rogue waves. 
%It has been recognized that the analytic structure of the solutions of the corresponding PDEs extended to the complex plane provides valuable information about the onset and evolution of these phenomena \cite{FH19, K17, W03, W22}. 
%In particular, the dynamics of complex singularities, such as poles and branch points, can reveal the nature and location of developing singular features.
%whether a singular feature is developing, where it is located, and on its nature. 
Consequently, the analytic continuation of the solutions into the complex plane has become an important tool for studying nonlinear PDEs.
In this paper, we present a novel neural network-based approach for constructing  \textit{numerical} analytic continuation of the solutions of nonlinear PDEs. The proposed method allows for the accurate detection and tracking of complex singularities directly from  the numerical data on the real axis, providing an effective framework for investigating and understanding aforementioned features of the solutions.

To present an analysis of the corresponding phenomena, we
 consider several well-known  model problems in one spatial dimension, namely the nonlinear heat equation (NLH), the nonlinear Burgers equation, and the nonlinear Schrödinger equation (NLS). We investigate two classes of initial conditions: (i) initial data that    possess singularities in the finite complex plane, and (ii) initial data that are entire functions throughout the finite complex plane upon complexification. For each model problem, we study the evolution of the associated phenomena in its relation to the behavior of complex singularities, with particular emphasis on their trajectories as functions of time.

Our method for analytic continuation to the complex plane and for the singularity tracking is based on the recently proposed neural network-based approach for learning  meromorphic functions  with pole-type singularities \cite{DKD2025}. We employ \textit{``unsafe'' Padé Activation Units} (PAUs) as activation functions. The concept of ``unsafe'' PAUs was introduced in \cite{DKD2025} and consists in adaptively constructing rational activation functions as meromorphic functions, each having a single pole located within the domain under investigation. Using the weights and biases of the hidden layer, we then scale and shift, respectively, the pole of the activation function to obtain estimated locations of the singularities; consequently, the number of neurons in the hidden layer is determined by the number of singularities of the function being approximated. While the hidden-layer weights and biases are tuned to capture the singularities, the weights and biases of the output layer are computed via least-squares fitting, ensuring an accurate approximation of the function over the rest of the domain. As shown in Theorem 3.1 of \cite{DKD2025}, this method  converges for meromorphic functions. For further information on rational neural networks based on ``safe'' PAUs (continuous rational activation functions), we refer the interested reader to \cite{BNT20, MSK, P22, Tel17}, and for other methods for training complex neural networks, to \cite{CL22, L22}.

There is a large body of work devoted to computing and tracking singularities of the extended solutions of nonlinear PDEs.We would like to mention the singularity-tracking method \cite{SS83}, which, although limited to characterizing the singularity closest to the real axis, was followed by more advanced methods such as the Borel–Polya–Van der Hoeven method \cite{PF2007} and the Kida technique \cite{Kida}, both of which can also be used to analyze hidden singularities. 
The most relevant to the present work are methods employing Padé–Fourier methods for numerical analytic continuation \cite{FW24, CG15, W03, W22}. Although our approach inherits some properties of rational Padé  approximation by construction, it differs significantly from the methods employed in \cite{FW24, CG15, W03, W22}. First, our method incorporates the SVD-based procedure introduced in \cite{GGT13} (see also Algorithm 2 in \cite{DKD2025}) to determine the optimal number of neurons, which in our setting coincides with the number of pole-type singularities. Consequently, when the target function is meromorphic, our approach effectively eliminates  spurious poles. For functions with branch cuts, the spurious poles and zeros are located along the branch cuts, as is standard for rational approximation. 
Second, we establish an explicit relationship between the singularities of the target function and the weights and biases of the hidden layer. 
Unlike in backpropagation-free approaches with fixed hidden-layer parameters \cite{DY22, DS20}, 
recent work \cite{DD24} has extended this to dynamically changing weights by resampling.  
In our framework, the weights and biases are considered as continuous functions on time,  allowing the network to adapt to the temporal evolution of the singularity structure. 
%In our framework, the values of weights and biases evolve dynamically in time,  allowing the network to adapt to the temporal evolution of the singularity structure. 
Finally, we note that other rational methods can also be used for numerical analytic continuation. We refer to \cite{NT23} for analytic continuation via the AAA method, and to \cite{Y22, Y222, D2025, YD26} for the approach based on the connection between rational functions and exponential sums.

\begin{figure}[htbp]
    \centering
    
    % Top row: two side-by-side plots
    \begin{subfigure}[b]{0.48\textwidth}
        \centering
        \includegraphics[width=\textwidth]{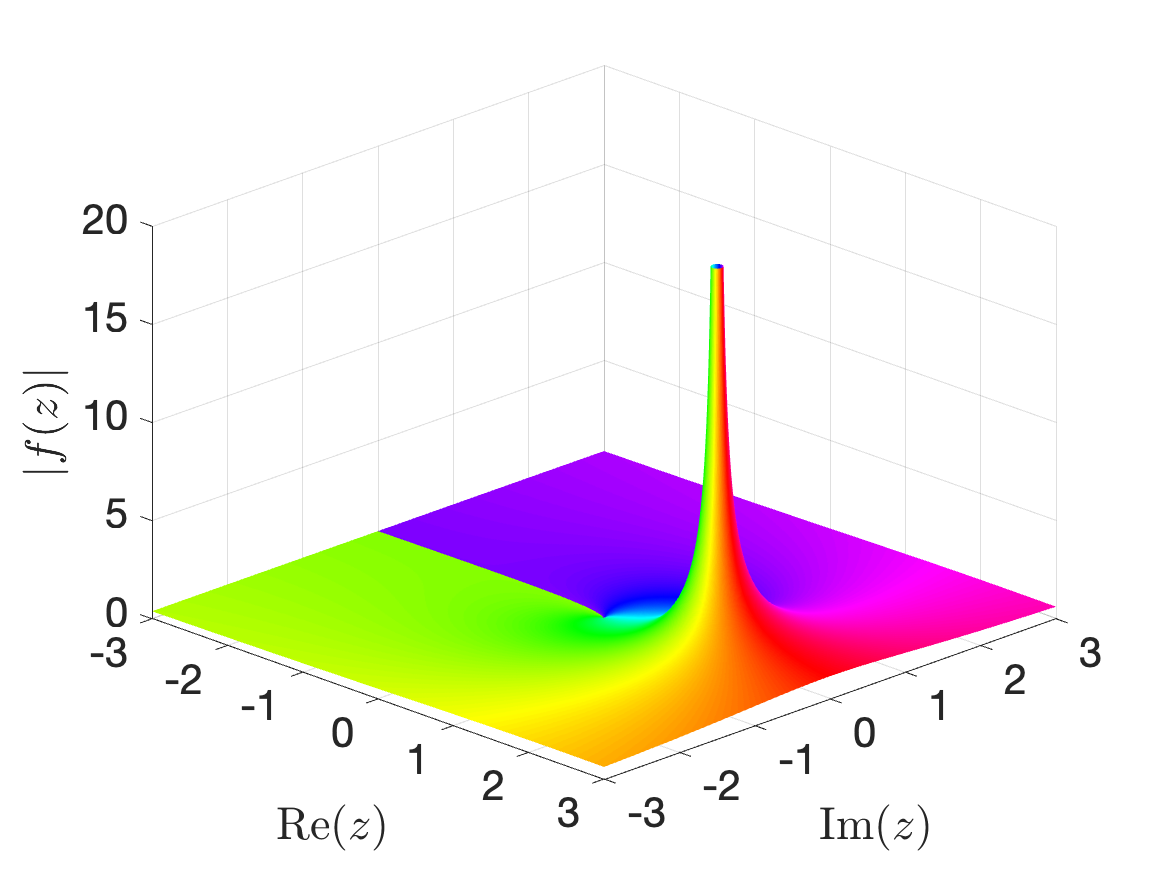}
       % \caption{The analytic landscape }
      %  \label{fig:1a}
    \end{subfigure}
    \hfill
    \begin{subfigure}[b]{0.48\textwidth}
        \centering
        \includegraphics[width=0.8\textwidth]{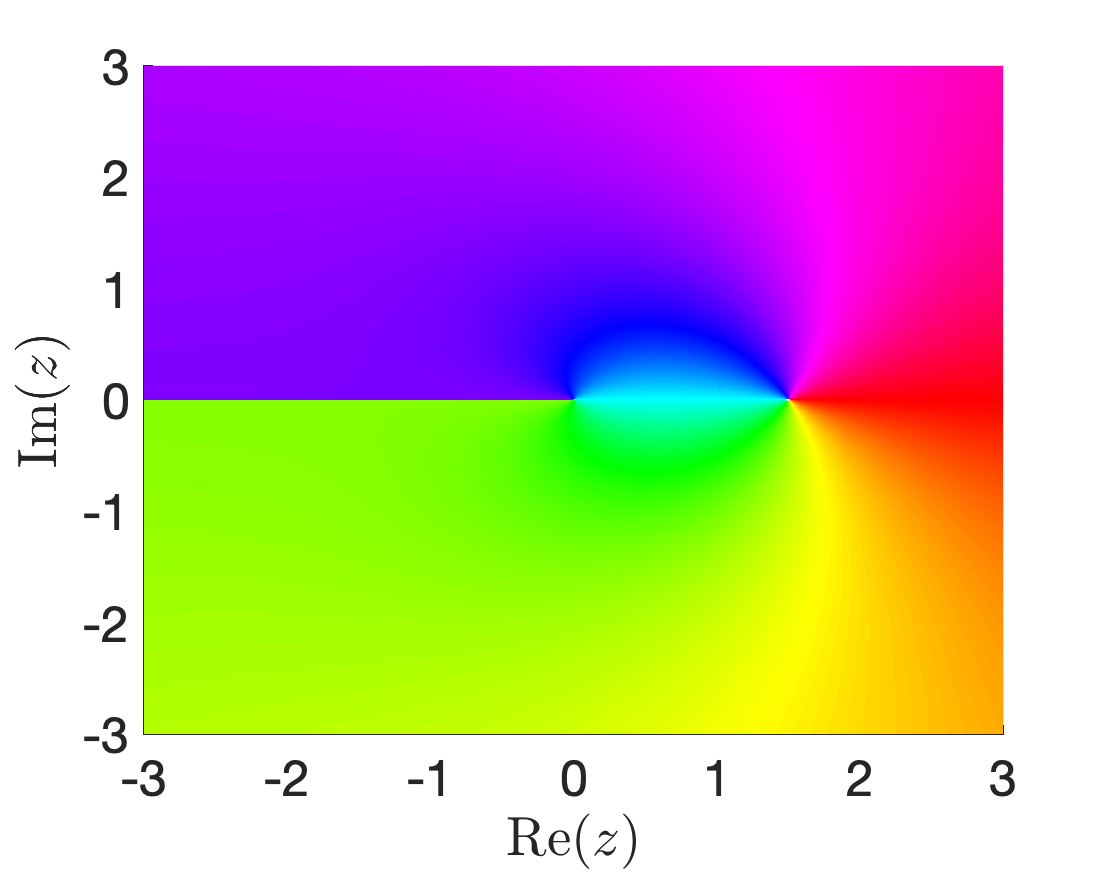}
      %  \caption{}
      %  \label{fig:1b}
    \end{subfigure}
    
  %  \vspace{0.5cm} % spacing between rows
    
    % Bottom: centered color wheel
    \begin{subfigure}[b]{0.35\textwidth}
        \centering
        \includegraphics[width=1.2\textwidth]{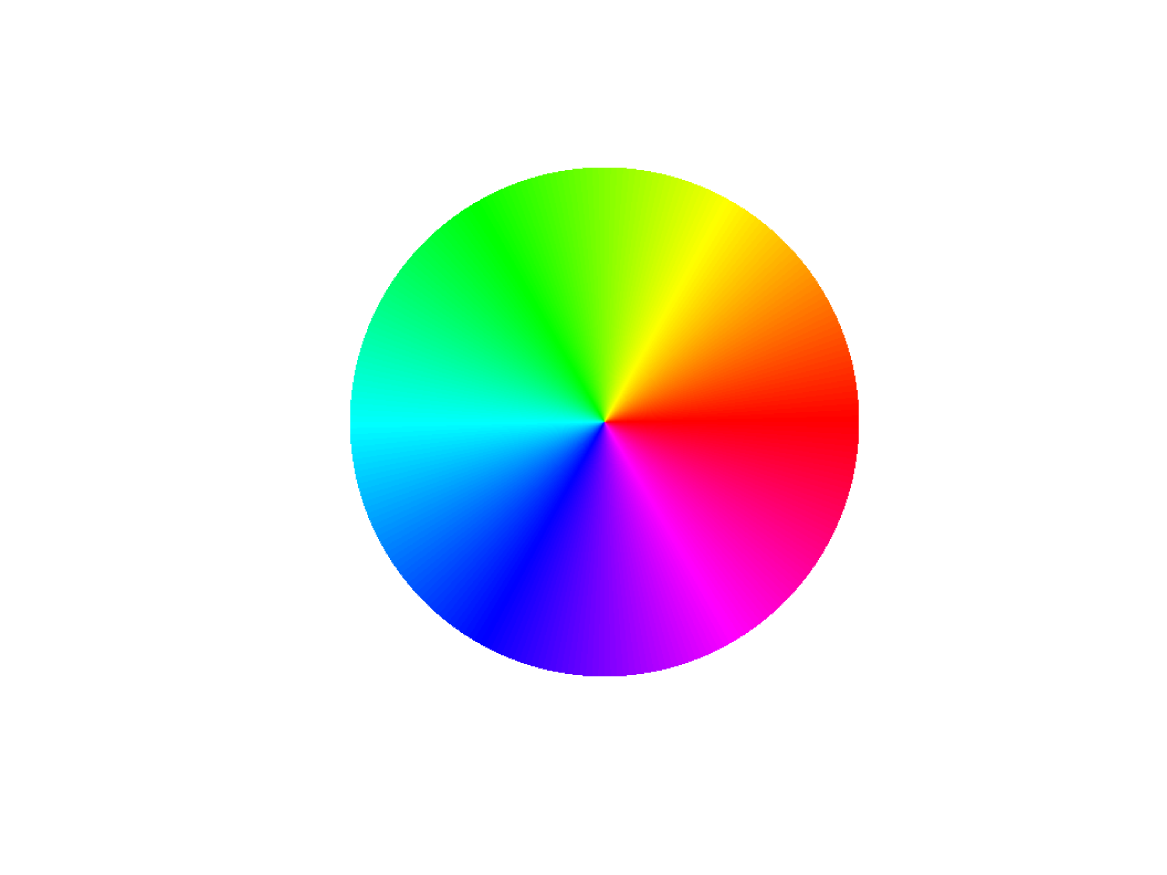}
       % \caption{}
      %  \label{fig:1c}
    \end{subfigure}
    
    \caption{  The analytic landscape (first row, left) and the phase portrait (first row, right) of the function $f(z) = \sqrt{z}/(z-1.5)$, $z \in \mathbb{C}$, colored according to the standard HSV color wheel (second row). The plots reveal the key features of $f(z)$: a simple pole at $z = 1.5$, an algebraic branch point of second order at $z = 0$, and a branch cut along the negative real axis.
    }
    \label{fig_color}
\end{figure}

To identify and visualize poles, branch points  and branch cuts in the complex plane, we employ two complementary representations: the analytic landscape and the phase portrait of the complex-valued function  \cite{Weg12}. 
Writing
\[
f(z)=r e^{\mathrm{i}\phi},
\]
where \(r=|f(z)|\) and \(\phi=\arg f(z)\), the function $f$ is visualized by a three-dimensional analytic landscape (Figure \ref{fig_color}, first row (left)), with the horizontal axes representing the real and imaginary parts of \(z=x+\mathrm{i}y\), the height representing the modulus \(r\), and the color encoding the phase \(\mathrm{e}^{\mathrm{i}\phi}\). The corresponding phase portrait  is encoded by coloring the surface according to the phase \(\mathrm{e}^{\mathrm{i}\phi}\) (Figure \ref{fig_color}, first row (right)).  In both visualizations, we employ the standard HSV color wheel (Figure \ref{fig_color}, second row).
%Writing $f(z)=r e^{\mathrm{i} \phi},$ where $r=|f(z)|$ and $\phi=\arg f(z)$, the software of \cite{Weg12} generates a three-dimensional analytic landscape in which the horizontal axes correspond to the real and imaginary parts of $z=x+\mathrm{i}  y$, while the height represents the modulus $r$, and color represents the phase $\phi$.  Similarly, the phase portrait \cite{Weg122} is encoded by coloring the surface according to the phase $\phi$, using the conventional CMYK color wheel. 
The analytic landscape highlights the magnitude of the function, making poles and zeros readily visible through peaks and valleys, whereas the phase portrait reveals the local phase variation, allowing branch points, poles, and other singularities to be identified through their characteristic color patterns. 

\textbf{Outline.}  In Section \ref{le}, we present the main ideas of our new neral network-based method for construction  of numerical analytic continuation of the solutions of nonlinear PDEs  and detection of their complex singularities. In Section \ref{secapp},  we present applications of our method to three model problems. In Subsections \ref{secnls1} and \ref{secburgers}, we consider the nonlinear heat equation and the nonlinear Burgers equations, respectively. Subsection \ref{secnls} is dedicated to the nonlinear Schrödinger equation. We finish the paper with the conclusions section.

\textbf{Notation.} %As usual, $\nn$, $\rr$, and $\cc$ denote the sets of natural, real, and complex numbers, respectively.  %Throughout the paper,  we  use the matrix notation $\mathbf{X}=(X_{ij})_{i,j=1}^{m,n} \in \cc^{m \times n}$ for matrices of size $m\times n$, and the vector notation $\boldsymbol{x}=(x_j)_{j=1}^n \in \cc^{ n}$ or $\boldsymbol{x}=(x_1,\ldots,x_n)^T \in \cc^{ n}$ for column vectors of size $n$.
Throughout the paper, the notation $X^{(\pm)}$ indicates that both cases, $X^{(+)}$ and $X^{(-)}$, are considered. The statement ``$X^{(\pm)}$ are valid for $Y^{(\pm)}$'' is understood componentwise: that is, $X^{(+)}$ is valid for $Y^{(+)}$, and $X^{(-)}$ is valid for $Y^{(-)}$.

\section{Neural network-based method for numerical analytic continuation and singularity detection}
\label{le}

We now present the main ideas underlying our new neural network-based method for constructing the  extensions of the  solutions of nonlinear PDEs into the complex plane and for detecting their complex singularities.

\subsection{ Description of the method }

\begin{figure}[h!]
  \centering
  %\begin{subfigure}[b]{0.45\linewidth}
    \includegraphics[width=1\linewidth]{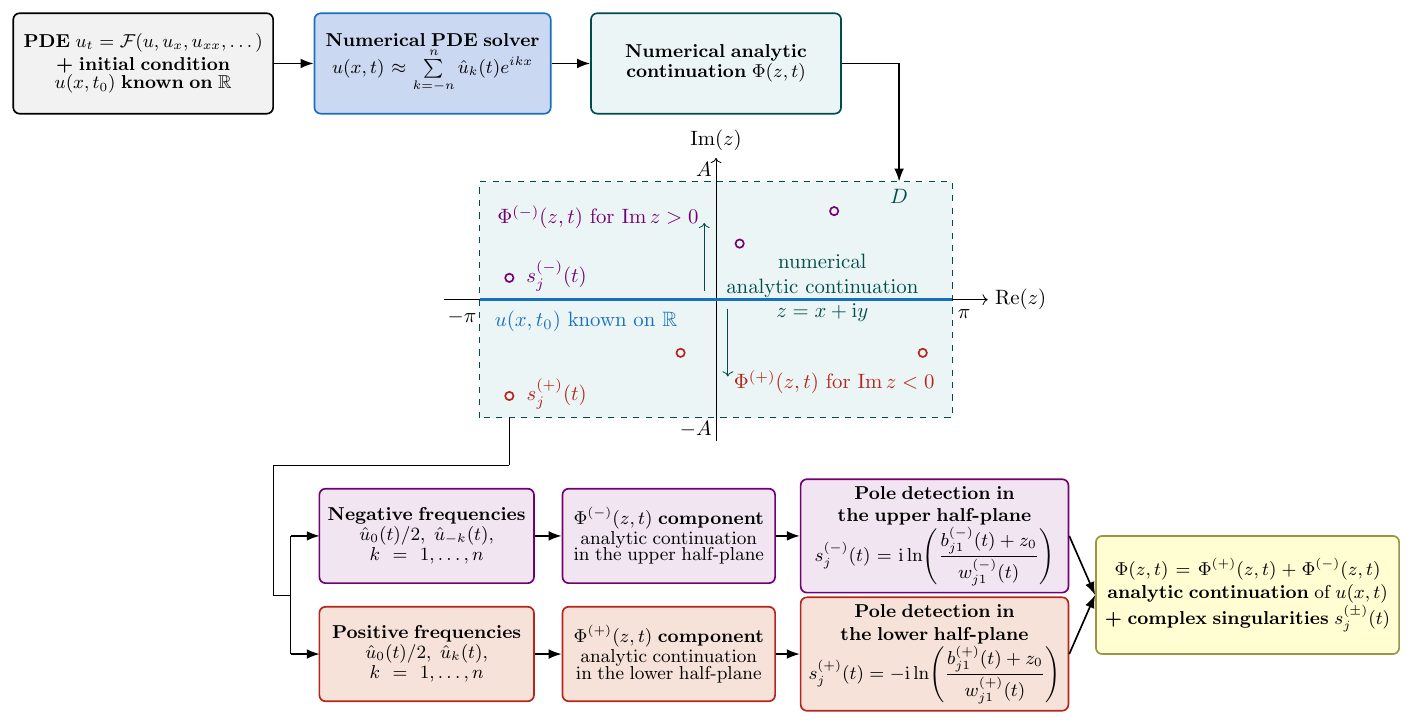}
   % \caption{$t=0.15$}
  %\end{subfigure}
 % \begin{subfigure}[b]{0.45\linewidth}
  %  \includegraphics[width=1\linewidth]{bh1.eps}
   % \caption{$t=0.503$}
  %\end{subfigure}
  \caption{
  A schematic illustration of the procedure for constructing the numerical  analytic continuation $\Phi(z,t)$ of the solution $u(x,t)$ of the PDE (\ref{pde}) and detecting its complex singularities.
  }
  \label{figschema}
\end{figure}

 We consider a one-dimensional evolutionary PDE of the form
\begin{equation}\label{pde}
    u_t=\mathcal{F}(u,u_x,u_{xx},\dots),
\end{equation}
where $u=u(x,t)$, $x \in [-\pi, \pi)$ is the spatial variable, $t\in \rr$ is the time variable, and by $u_t$ and $u_x$, $u_{xx}$ we denote the various partial derivatives of $u$ with respect to $t$ and $x$ at different orders. We also assume the initial condition at $t_0 \in \rr$
\begin{equation}\label{incod}
    u(x,t_0)=u_0(x).
\end{equation}
%As in \cite{DKD2025}, we apply the following approach for construction of analyticcontinuation of the solution $u(x,t)$.
Let $D:=  [-\pi, \pi)\times \mathrm{i} (-A,A)$ for some $A>0$. For each $t\in \rr$, by $u(z,t): D  \rightarrow \cc$ 
we  denote the analytic continuation of the solution of the PDE (\ref{pde}) (i.e., formally instead of $x \in \rr$ we use a complex variable  $z=x+\mathrm{i} y \in \cc$ with $y \in (-A,A)$). 
The proposed method consists of three main steps. First, we compute an approximate solution  of the PDE \eqref{pde} in the interval \([-\pi,\pi)\subset\mathbb{R}\) using a numerical solver. Then we construct its analytic continuation to a domain \(D\subset\mathbb{C}\) by means of the neural network-based method proposed in \cite{D2025}. Finally, we exploit the hidden-layers parameters of the corresponding neural network to identify the complex singularities of the  solution extended into the complex plane. 
 A schematic overview of this pipeline is presented in Figure~\ref{figschema}.

We now present a brief description of these three main steps.
In the first step, we use
 standard numerical methods (such as the Fourier spectral method and the split-step Fourier method) to approximate the corresponding solution of the PDE $u(x,t)$ in $[-\pi,\pi)$ by the truncated Fourier series
\begin{equation}\label{appsr}
    u(x,t)\approx \sum\limits_{k=-n}^n \hat{u}_k(t) \, \mathrm{e}^{\mathrm{i} k x},
\end{equation}
for some $n \in \nn$.  As input, we employ a finite set of samples $u\left( \frac{ \pi j}{ n} ,t_0 \right)$, $j=-n,\dots, n-1$, of the initial condition (\ref{incod}). 
In case of the Fourier spectral method, the Fourier coefficients $\hat{u}_k(t_0)$ are computed by FFT using the initial condition values  $u\left( \frac{ \pi j}{ n} ,t_0 \right)$, $j=-n,\dots, n-1$.  Applying the Fourier spectral method,  the coefficients $\hat{u}_k(t)$ for another $t\neq t_0$ are obtained as the solution to the system of ODEs
%\begin{equation}\label{ode}
\begin{equation}\label{odes}
 \frac{\mathrm{d} \hat{u}_k(t)}{\mathrm{d} t } =  \mathcal{G}_k(\hat{\boldsymbol{u}}(t)), \quad k=-n,\dots,n,
 \end{equation}
 where $ \hat{\boldsymbol{u}}(t):= (\hat{u}_\ell(t))_{\ell=-n}^n$. The system (\ref{odes}) is derived  from  (\ref{pde}) using (\ref{appsr}) as ansatz for the solution. 
 %It can be done by a standard method oflines employing build-in functions \textmd{ode45} or \textmd{ode15s} in Matlab. 
For the split-step Fourier method, the initial condition values $u\left( \frac{ \pi j}{ n} ,t_0 \right)$, $j=-n,\dots, n-1$ are first used to compute the  samples $u\left( \frac{ \pi j}{ n} ,t \right)$, $j=-n,\dots, n-1$ at each time $t\in \rr$. The corresponding Fourier coefficients $\hat{u}_k(t)$ are then obtained by applying FFT to the corresponding function values. Note that in both cases we compute  $\hat{u}_k(t)$ for $k=-n,\ldots,n-1$ and then set $\hat{u}_n(t)=\hat{u}_{-n}(t)$ to obtain approximation (\ref{appsr}) with symmetric frequencies.

%first, employing again the initial condition values  $u\left( \frac{ \pi j}{ n} ,t_0 \right)$, $j=-n,\dots, n-1$, we compute  samples $u\left( \frac{ \pi j}{ n} ,t \right)$, $j=-n,\dots, n-1$ and then  coefficients $\hat{u}_k(t)$ by FFT for $t\in \rr$. 
%\end{equation}

In the second step, we extend the numerical solution $u(x,t)$ into the complex plane. 
Formally, replacing the real variable $x\in\mathbb{R}$ with a complex variable $z\in\mathbb{C}$ in (\ref{appsr}), we obtain
%\begin{equation}\label{appsc}
$$
    u(z,t)\approx \sum\limits_{k=-n}^n \hat{u}_k(t) \, \mathrm{e}^{\mathrm{i} k z}=:  u^{(-)}(z,t)+u^{(+)}(z,t),
$$
%\end{equation}
where  
$$
u^{(-)}(z,t) :=\frac{1}{2} \hat{u}_0(t) +\sum_{k=1}^{n} \hat{u}_{-k}(t) \mathrm{e}^{-\mathrm{i} k z} \quad \text{and} \quad u^{(+)}(z,t):=\frac{1}{2} \hat{u}_0(t)+\sum_{k=1}^n \hat{u}_k(t) \mathrm{e}^{\mathrm{i} k z}.
$$
%$ u^{(-)}(z,t) :=\frac{1}{2} \hat{u}_0(t) +\sum_{k=1}^{n} \hat{u}_{-k}(t) \mathrm{e}^{-\mathrm{i} k z}$ and $ u^{(+)}(z,t):=\frac{1}{2} \hat{u}_0(t)+\sum_{k=1}^n \hat{u}_k(t) \mathrm{e}^{\mathrm{i} k z}$.
 Next, we construct a numerical analytic continuation of $u(x,t)$ in the form of a  neural network  $\Phi: \cc \rightarrow \cc$ defined by
\begin{equation}\label{nnapr}
    \Phi(z,t):=  \Phi^{(+)}(z,t) +  \Phi^{(-)}(z,t),
\end{equation}
such that 
%and we get the final formula for the extended solution in $\cc$
\begin{equation}\label{appsnn}
    u(z,t) \approx \Phi(z,t),
\end{equation}
and the neural network components $\Phi^{(\pm)}$ are given by formulas 
%functions $u^{(\pm)}(z,t)$ are approximated by the corresponding neural network components
\begin{align}
   \Phi^{(+)}(z,t) & :=  \mathbf{W}_2^{(+)}(t) \, r^{(+)} \left( \mathbf{W}_1^{(+)}(t) \mathrm{e}^{\mathrm{i}  z} - \boldsymbol{b}_1^{(+)}(t)  \right) - b_2^{(+)}(t), \label{f1t} \\
     \Phi^{(-)}(z,t) & :=  \mathbf{W}_2^{(-)}(t) \, r^{(-)} \left( \mathbf{W}_1^{(-)} (t)\mathrm{e}^{-\mathrm{i} z} - \boldsymbol{b}_1^{(-)} (t) \right) - b_2^{(-)}(t). \label{f2t}
\end{align}
%such that 
%and we get the final formula for the extended solution in $\cc$
%\begin{equation}\label{appsnn}
 %   u(z,t) \approx \Phi(z,t).
%\end{equation}
In (\ref{f1t})-(\ref{f2t}), the numbers $M^{(\pm)}(t) \in \nn$ of neurons in the hidden layers, weights  $\mathbf{W}_1^{(\pm)}(t)=(w_{11}^{(\pm)}(t),\dots,w_{M^{(\pm)}(t)1}^{(\pm)}(t))^T \in \cc^{M^{(\pm)}(t) \times 1}$ and biases $\boldsymbol{b}_1^{(\pm)}(t)=(b_{1 1}^{(\pm)}(t),\dots,b_{M^{(\pm)}(t) 1}^{(\pm)}(t))^T \in \cc^{M^{(\pm)}(t) }$ of the hidden layers, as well as
weights  $\mathbf{W}_2^{(\pm)}(t) = ( w_{1 2}^{(\pm)}(t),\dots,w_{ M^{(\pm)}(t) 2}^{(\pm)}(t)) \in \cc^{1 \times M^{(\pm)}(t) }$ and biases $b_2^{(\pm)}(t) \in \cc$ of the output layers of the components $\Phi^{(\pm)}$ \textit{depend on time}. By $r^{(\pm)}$ in (\ref{f1t})-(\ref{f2t}), we denote ``unsafe'' PAUs, i.e.  rational  activation functions with a single pole each situated in the domain under investigation \cite{DKD2025}. We construct $r^{(\pm)}$ by Algorithm 3 from \cite{DKD2025} as types $(N^{(\pm)}(t)-M^{(\pm)}(t)+1,1)$ rational Laurent-Padé approximants  to the function 
\begin{equation}\label{actfun}
    \omega(z)= \frac{\varphi(z)}{z-z_0}, \quad z_0 \in \cc,  \, |z_0|>1, \, \varphi(z) \text{ is analytic in } \cc.
\end{equation} 
The numerical analytic continuation $\Phi$ is constructed so that the component $\Phi^{(-)}$ approximates the extended solution $u(z,t)$ and detects its complex singularities in the upper half-plane. Correspondingly, the component $\Phi^{(+)}$ approximates the  $u(z,t)$ and detects its complex singularities in the lower half-plane (see Figure \ref{figschema}).
 We present a brief overview of the procedure used to compute parameters of the neural network $\Phi$ in \eqref{nnapr} in Subsection~\ref{parcom}.

In the third step of our pipeline, we compute estimated locations of the complex singularities of the extended solution $u(z,t)$.
According to \cite{DKD2025}, 
the component $\Phi^{(+)}$  detects the locations of  singularities by the formula
\begin{equation}\label{spl}
\mathrm{e}^{\mathrm{i} s_{j}^{(+)}(t)}=  \frac{b^{(+)}_{j 1}(t)+z_0 }{w^{(+)}_{j 1}(t)}  \quad \text{or }   \quad 
 s_{j}^{(+)}(t)= -\mathrm{i}  \ln \left(  \frac{b^{(+)}_{j 1}(t)+z_0 }{w^{(+)}_{j 1}(t)} \right),
\end{equation}
such that $\mathrm{Im} \, s_{j}^{(+)}(t)<0$ for $j=1,\dots,M^{(+)}(t)$ at given time point $t$. Correspondingly,  the component $\Phi^{(-)}$ detects the locations of singularities by
\begin{equation}\label{spl1}
 \mathrm{e}^{-\mathrm{i} s_{j}^{(-)}(t)}=\frac{b^{(-)}_{j 1}(t)+z_0}{ w^{(-)}_{j 1}(t) }  \quad \text{or }   \quad      s_{j}^{(-)}(t)= \mathrm{i}  \ln \left(  \frac{b^{(-)}_{j 1}(t)+z_0}{ w^{(-)}_{j 1}(t) } \right), 
\end{equation}
such that $\mathrm{Im} \, s_{j}^{(-)}(t)>0$ for $j=1,\dots,M^{(-)}(t)$ at given time point $t$. As can be seen from equations (\ref{spl})-(\ref{spl1}),  the singularities $s_{j}^{(\pm)}(t)$ are detected by scaling and shifting a pole $z_0$ of the activation functions  by weights and biases, respectively, of the hidden layers of the components $\Phi^{(\pm)}$. 
%Note that the weights $w^{(\pm)}_{\ell 1}(t)$  and biases    $b^{(\pm)}_{\ell 1}(t)$ are carrying information about $\mathrm{e}^{\mathrm{i}  s_{\ell}^{(\pm)}(t)}$ (and thus not directly about $s_{\ell}^{(\pm)}(t)$).  
The parameters $M^{(\pm)}(t) \in \nn$  are determined by the numbers of singularities of the extended solution $u(z,t)$ located in the upper-half ($M^{(-)}(t)$) and lower-half ($M^{(+)}(t)$) planes, which can be different for different time points. We compute   $M^{(\pm)}(t)$ by  the SVD-based procedure from \cite{GGT13} (see also Algorithm 2 in \cite{DKD2025}). 

\subsection{Computation of the parameters of the neural network}
\label{parcom}

Finally, we describe the main ideas of the construction of the neural network $\Phi$ in \eqref{nnapr}. To achieve the approximation (\ref{appsnn}),
the parameters of the neural network (\ref{nnapr})  are computed by Algorithm 5 from \cite{DKD2025} with the substitution $\mathrm{e}^{\mathrm{i}  z} = v \in \cc$. 
The main idea of this procedure is that the components  $\Phi^{(\pm)}$ are constructed as \textit{Fourier-Padé approximants} to $u^{(\pm)}$, i.e. $\Phi^{(\pm)}=p^{(\pm)}_{N^{(\pm)}(t)}/ q^{(\pm)}_{M^{(\pm)}(t)}$  and the following properties hold
\begin{equation} \label{pplus} 
    q^{(+)}_{M^{(+)}(t)}(z,t) \,  u^{(+)}(z,t)- p^{(+)}_{N^{(+)}(t)}(z,t) = \mathcal{O}(\mathrm{e}^{\mathrm{i} \, (N^{(+)}(t)+M^{(+)}(t)+1) z}),  
\end{equation}  
and 
\begin{equation} \label{pminus}
    q^{(-)}_{M^{(-)}(t)}(z^{-1},t) \, u^{(-)}(z,t)- p^{(-)}_{N^{(-)}(t)}(z^{-1},t) = \mathcal{O}(\mathrm{e}^{-\mathrm{i} \, (N^{(-)}(t)+M^{(-)}(t)+1) z}), 
\end{equation}   
where by $q^{(\pm)}_{M^{(\pm)}(t)}$ and $p^{(\pm)}_{N^{(\pm)}(t)}$  we denote polynomials of degrees  $M^{(\pm)}(t) \in \nn$ and $N^{(\pm)}(t) \in \nn$, respectively.  

As we already mentioned, ``unsafe'' PAUs  $r^{(\pm)}$ in (\ref{f1t})-(\ref{f2t}) are computed by Algorithm 3 from \cite{DKD2025}  as  types $(N^{(\pm)}(t)-M^{(\pm)}(t)+1,1)$ rational Laurent-Padé approximants  to the function (\ref{actfun}). Thus, the activation functions $r^{(\pm)}$  have the form
\begin{align}
r^{(\pm)}(z,t) & =\frac{n_{N^{(\pm)}(t)+1-M^{(\pm)}(t)}^{(\pm)}(z,t)}{d_1^{(\pm)}(z,t)}:= \frac{\sum_{j=0}^{N^{(\pm)}(t)+1-M^{(\pm)}(t)} \alpha_j^{(\pm)}(t) \,  z^j }{\gamma_1^{(\pm)}(t) \, z +\gamma_0^{(\pm)}(t)}. \label{r1}
%& :=\frac{\alpha^{(\pm)}_{N^{(\pm)}(t)+1-M^{(\pm)}(t)} (t) z^{N^{(\pm)}(t)+1-M^{(\pm)}(t)} +\dots+\alpha_1^{(\pm)}(t) z+\alpha_0^{(\pm)}(t)}{\gamma_1^{(\pm)}(t) z +\gamma_0^{(\pm)}(t)}. \label{r1}
\end{align}
The coefficients $\alpha_j^{(\pm)}(t) \in \cc$, $j=0,\ldots, N^{(\pm)}(t)+1-M^{(\pm)}(t)$ and $\gamma_1^{(\pm)}(t), \, \gamma_0^{(\pm)}(t) \in \cc$ are computed from the conditions $$n_{N^{(\pm)}(t)+1-M^{(\pm)}(t)}^{(\pm)}(z,t) \, \omega(z) - d_1^{(\pm)}(z,t) = \mathcal{O}(z^{N^{(\pm)}(t)+2-M^{(\pm)}(t)}).$$  From the Montessus de Ballore theorem regarding convergence of a Laurent-Pad\'{e} approximation (see, for example, \cite{B87}), we have that if $N^{(\pm)}(t)+1-M^{(\pm)}(t) \rightarrow \infty$, then $r^{(\pm)}(z,t) \rightarrow \omega(z)$  uniformly on any compact subset of $\{z \in \mathbb{C}: |z|< \varrho^{(\pm)} \} \setminus \{z_0\}$ for some $\varrho^{(\pm)}>1$, and $-\gamma_0^{(\pm)}/\gamma_1^{(\pm)} \rightarrow z_0$. The types of rational activation functions $r^{(\pm)}$ sufficiently depend on the numbers $M^{(\pm)}(t)$ of complex singularities of the corresponding solution $u(x,t)$.
 The choice of $z_0$ and $\varphi(z)$ does not have any influence on the approximation accuracy of the neural network $\Phi$. Note that  coefficients of the activation functions $r^{(\pm)}$ are also time-dependent.

%Parameters $M^{(\pm)}(t) \in \nn$ correspond to the number of singularities 
%The components $\Phi^{(+)}$ and $ \Phi^{(-)}$ are constructed using ideas similar to those employed in Algorithm \ref{alg3}, but replacing Laurent coefficients by Fourier coefficients $a_k(t)$ and taking into account their dependency on time. 

To compute  the weights $ w_{j1}^{(\pm)}(t)$ and biases $b_{j 1}^{(\pm)}(t) $, $j=1,\dots,M^{(\pm)}(t)$ of the hidden layers of the neural network components $\Phi^{(\pm)}$, we apply a new  backpropagation-free method developed in  \cite{DKD2025}.
%In \cite{DKD2025}, we presented a new  backpropagation-free method for computing these weights $ w_{j1}^{(\pm)}(t)$ and biases $b_{j 1}^{(\pm)}(t) $ of the hidden layers of neural network components $\Phi^{(\pm)}$ for learning functions with pole-type singularities. 
Further, we present the main ideas of  this  method. Let  $\boldsymbol{q}^{(\pm)}_t:=(q_0^{(\pm)}(t),\dots,q_{M^{(\pm)}(t)}^{(\pm)}(t))^T$ be the coefficient vectors  of the polynomials
$$
q^{(+)}_{M^{(+)}(t)}(z,t): = \sum_{\ell=0}^{M^{(+)}(t)} q_\ell^{(+)}(t) \mathrm{e}^{\mathrm{i} \ell z} \quad \text{ and } \quad q^{(-)}_{M^{(-)}(t)}(z,t): = \sum_{\ell=0}^{M^{(-)}(t)} q_\ell^{(-)}(t) \mathrm{e}^{-\mathrm{i} \ell z}
$$
 as in (\ref{pplus}) and (\ref{pminus}), respectively.   First, employing coefficients $q_\ell^{(\pm)}(t)$, $\ell=0,\dots,M^{(\pm)}(t)$, we determine parameters 
 $C_{j0}^{(\pm)}(t)$ and $C_{j1}^{(\pm)}(t)$, $j=1,\dots,M^{(\pm)}(t)$ from the conditions
\begin{align}
    \prod\limits_{j=1}^{M^{(+)}(t)} (C_{j0}^{(+)}(t)+C_{j1}^{(+)}(t) \, \mathrm{e}^{\mathrm{i}  z} ) & = \sum\limits_{\ell=0}^{M^{(+)}(t)} q_\ell^{(+)}(t) \, \mathrm{e}^{\mathrm{i} \ell z},  \label{ccoef1} \\
    \prod\limits_{j=1}^{M^{(-)}(t)} (C_{j0}^{(-)}(t)+C_{j1}^{(-)}(t) \, \mathrm{e}^{-\mathrm{i}  z} ) &= \sum\limits_{\ell=0}^{M^{(-)}(t)} q_\ell^{(-)}(t) \, \mathrm{e}^{-\mathrm{i} \ell z}. \label{ccoef2}
\end{align}
Then using  coefficients $\gamma_0^{(\pm)}(t)$ and $\gamma_1^{(\pm)}(t)$ of the denominators   of the activation functions $r^{(\pm)}$ in (\ref{r1}), we compute weights and biases of the hidden layers of the components $\Phi^{(\pm)}$ by
\begin{align}
    w_{j1}^{(\pm)}(t) & =\frac{1}{\gamma_1^{(\pm)}(t)} C^{(\pm)}_{j1}(t), \quad j=1,\dots,M^{(\pm)}(t), \label{w1} \\ 
    b_{j 1}^{(\pm)}(t)  & =\frac{1}{\gamma_1^{(\pm)}(t)} (\gamma_0^{(\pm)}(t) - C^{(\pm)}_{j0}(t)),  \quad j=1,\dots,M^{(\pm)}(t),  \label{b1}
\end{align}
respectively. 
Since formulas (\ref{ccoef1})-(\ref{ccoef2}) give us $M^{(\pm)}(t)-1$ degrees of freedom for the choice of the  parameters  $C_{j0}^{(\pm)}(t)$ and $C_{j1}^{(\pm)}(t)$, $j=1,\dots,M^{(\pm)}(t)$, we can compute $C_{j0}^{(\pm)}(t)$, $j=1,\dots,M^{(\pm)}(t)-1$ randomly. 
Note, that from equations (\ref{w1})-(\ref{b1}) we can conclude that the random choice of $C_{j0}^{(\pm)}(t)$, $j=1,\dots,M^{(\pm)}(t)-1$  has an impact on the values of weights and biases of the hidden layers of the components $\Phi^{(\pm)}$, but does not affect the accuracy of computation of the estimated locations of singularities of the extended solution $u(z,t)$, since those are controlled by the parameters $q_\ell^{(\pm)}(t)$, $\ell=0,\ldots,M^{(\pm)}(t)$. 
Note also  that a similar approach to determine the weights and biases of a neural network from the parameters of rational approximation was considered in \cite{P22}, which employs ``safe'' PAUs together with the bisection method and the differential-correction algorithm for rational approximation. Finally, weights and biases of the output layers of the components $\Phi^{(\pm)}$ are computed via
least-squares fitting.

\section{Applications}
\label{secapp}

%Blow-up, shock formation, and rogue waves are among the most important nonlinear phenomena encountered in mathematical physics and engineering. Although they arise in different classes of nonlinear PDEs, they all originate from the evolution of complex singularities in the analytically continued solution. Consequently, they provide an ideal set of benchmark problems for assessing the capability of the proposed neural network-based numerical analytic continuation method.

In this section, we investigate the formation of blow up, shock, and rogue waves phenomena using the method proposed in Section~\ref{le}. To this end, we consider several well-known nonlinear PDEs. 
Note that  the   activation functions  $r^{(\pm)}$ in (\ref{r1}) for each example   are constructed as the Laurent-Pad\'{e} approximation of type $(1,1)$ to $\omega(z)=\frac{\cos z}{z-z_0}$. Unless otherwise indicated, we employ $z_0=-2$. This choice of the pole of the activation function is arbitrary. It affects the numerical values of the weights and biases of the hidden layers but does not affect the neural network's ability to detect singularities or approximate the target function.
On the one hand, we are interested in the accuracy of the singularity detection by equations (\ref{spl})-(\ref{spl1}). On the other hand,
to demonstrate the approximation performance of the neural network   $\Phi(z,t)$ in the domain $D$ at fixed time $t$, we employ the     pole-scaled approximation error defined by
\begin{equation}\label{errorrel}
    %E(z,t)=\frac{|u(z,t)-\Phi(z,t)|}{|u(z,t)|+\varepsilon},
    E(z,t):=|u(z,t)-\Phi(z,t)|\cdot\prod\limits_{j=1}^M |z-z_j(t)|^{m_j}, \quad z \in D,
\end{equation}
where $u(z,t)$ is the explicit extended solution and $z_j(t)$ are its poles in $D$  with multiplicities $m_j$, $j=1,\dots,M$. The error $E(z,t)$ is designed to remove the artificial blow up caused by singularities and observe the approximation quality itself.

\subsection{Nonlinear heat equation (NLH)}
\label{secnls1}

Starting with the category of  simple second-order parabolic PDEs, we study finite-time blow up phenomena in solutions of  NLHs.

\begin{example}
\label{ex31}
To begin, we consider 
 a NLH of the form \cite{W03}
\begin{equation}\label{heateq}
    u_t-\frac{1}{6} u_{xx}+u^2+u=0,
\end{equation}
with the explicit solution  given by 
\begin{equation}\label{heatsol}
    u(x,t)=\frac{\mathrm{e}^{5/3(\sigma-t)+2 \mathrm{i} x}}{(\mathrm{i} \mathrm{e}^{5/6(\sigma-t)+ \mathrm{i} x} +\varrho)^2},
\end{equation}
where  $\sigma$ and $\varrho$ are some constants and $x \in [-\pi,\pi)$. The singularities of the analytic continuation $u(z,t)$ are poles of second order given by
\begin{equation}\label{heatsing}
%$$
z(t)=2 \pi m+ \mathrm{Im} \, (\ln(\mathrm{i} \, \varrho)) + \mathrm{i} \left( \frac{5}{6}(\sigma-t) - \mathrm{Re} \, (\ln(\mathrm{i} \, \varrho))\right), \quad m=0, \pm 1, \pm 2,\dots
%$$
\end{equation}
The poles (\ref{heatsing}) move  toward the real axis and reach it at time $t_b=\sigma-\frac{6}{5} \mathrm{Re} \, (\ln(\mathrm{i} \, \varrho))$ at which occurs   blow up of the solution $u(x,t)$.

\begin{figure}[h!]
  \centering
  \begin{subfigure}[b]{0.24\linewidth}
    \includegraphics[width=1.1\linewidth]{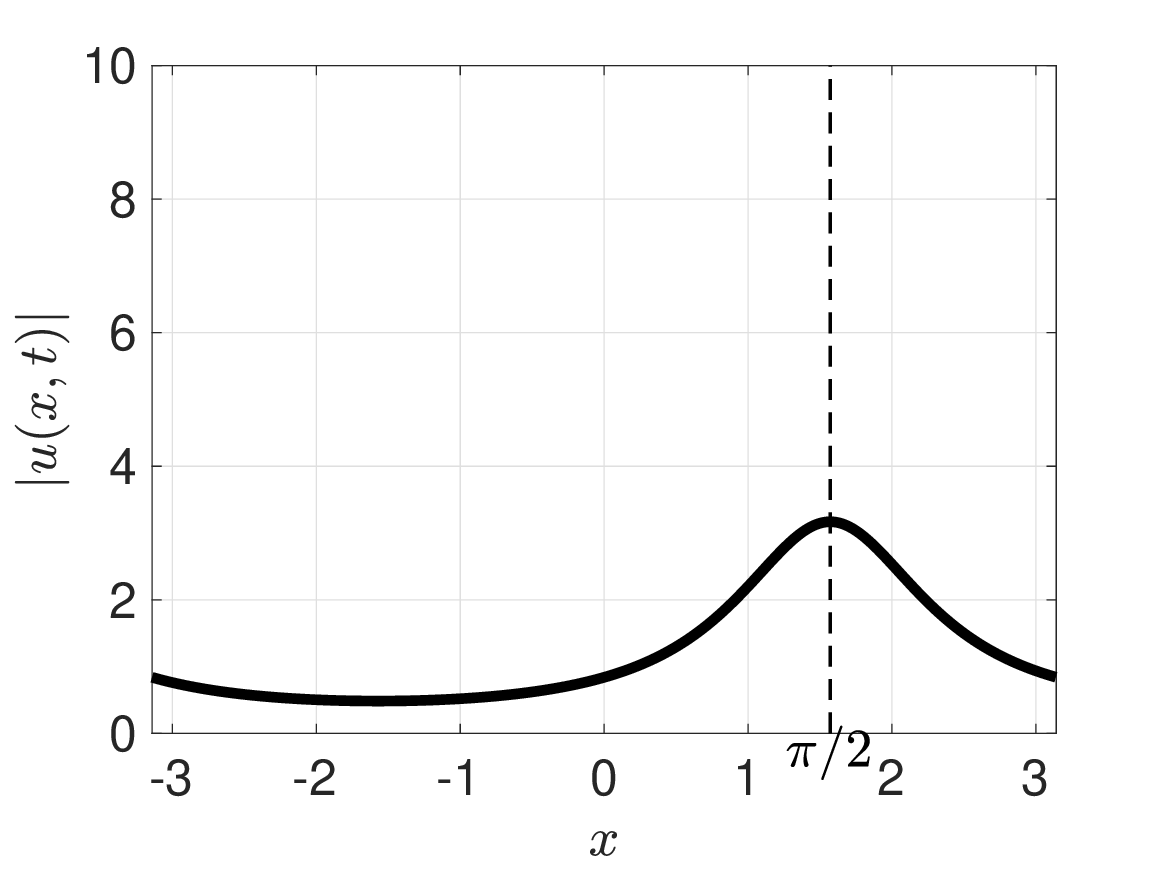}
   % \caption{$t=0$}
  \end{subfigure}
    \begin{subfigure}[b]{0.24\linewidth}
    \includegraphics[width=1.1\linewidth]{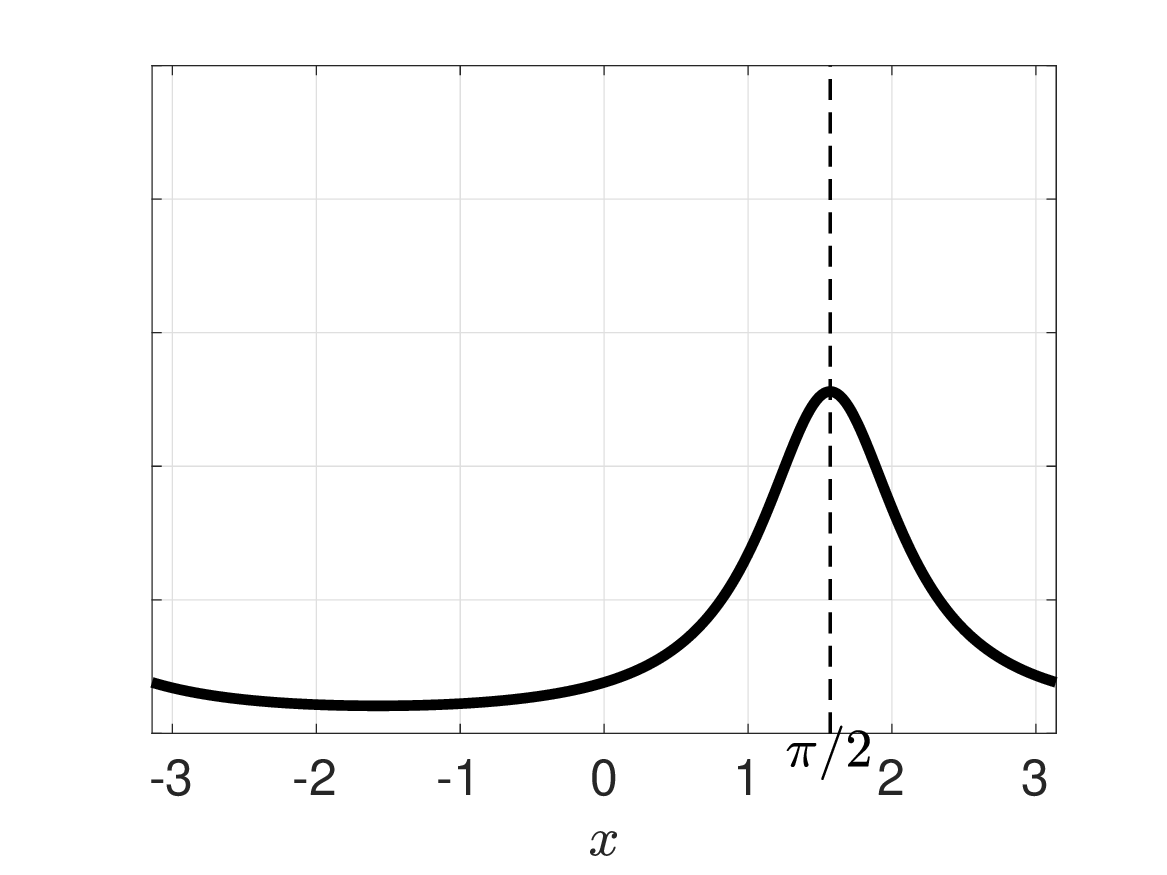}
   % \caption{$t=0.3$}
  \end{subfigure}
  \begin{subfigure}[b]{0.24\linewidth}
    \includegraphics[width=1.1\linewidth]{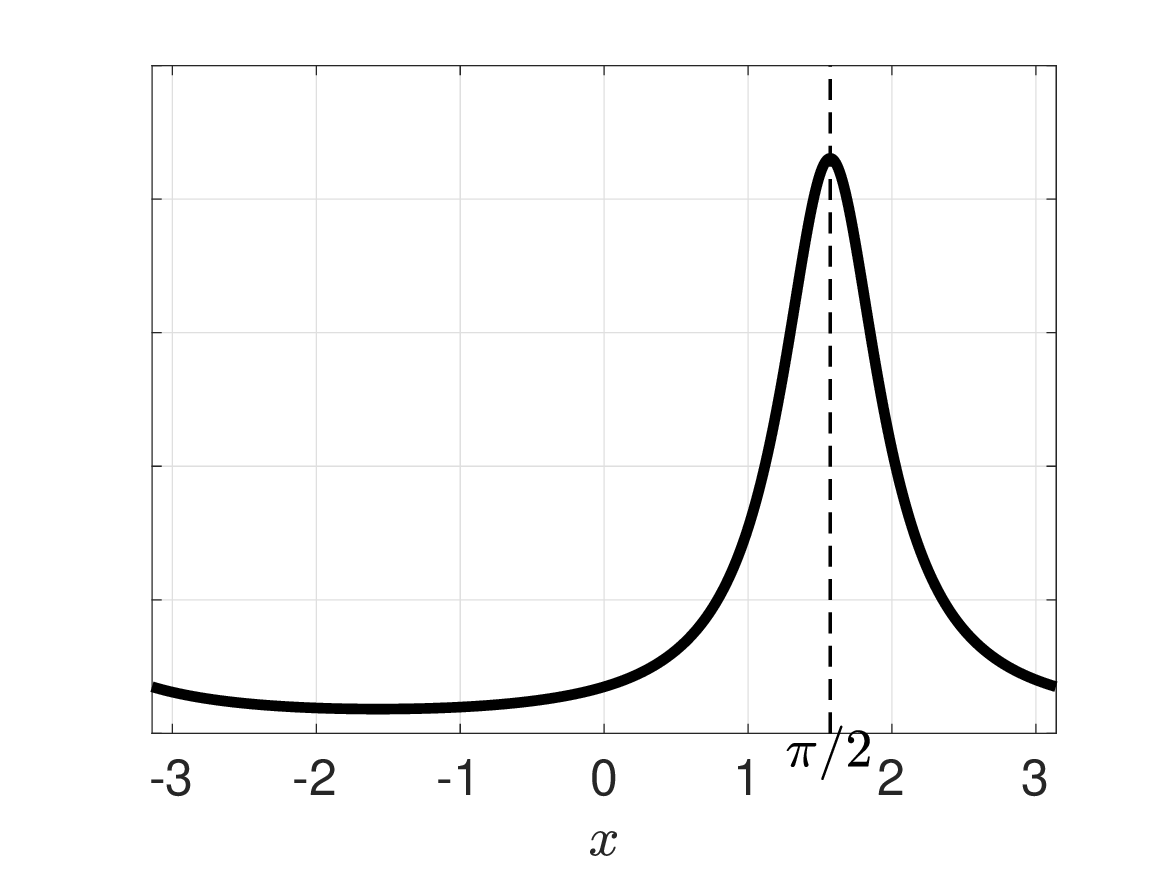}
   % \caption{$t=0.3$}
  \end{subfigure}
  \begin{subfigure}[b]{0.24\linewidth}
    \includegraphics[width=1.1\linewidth]{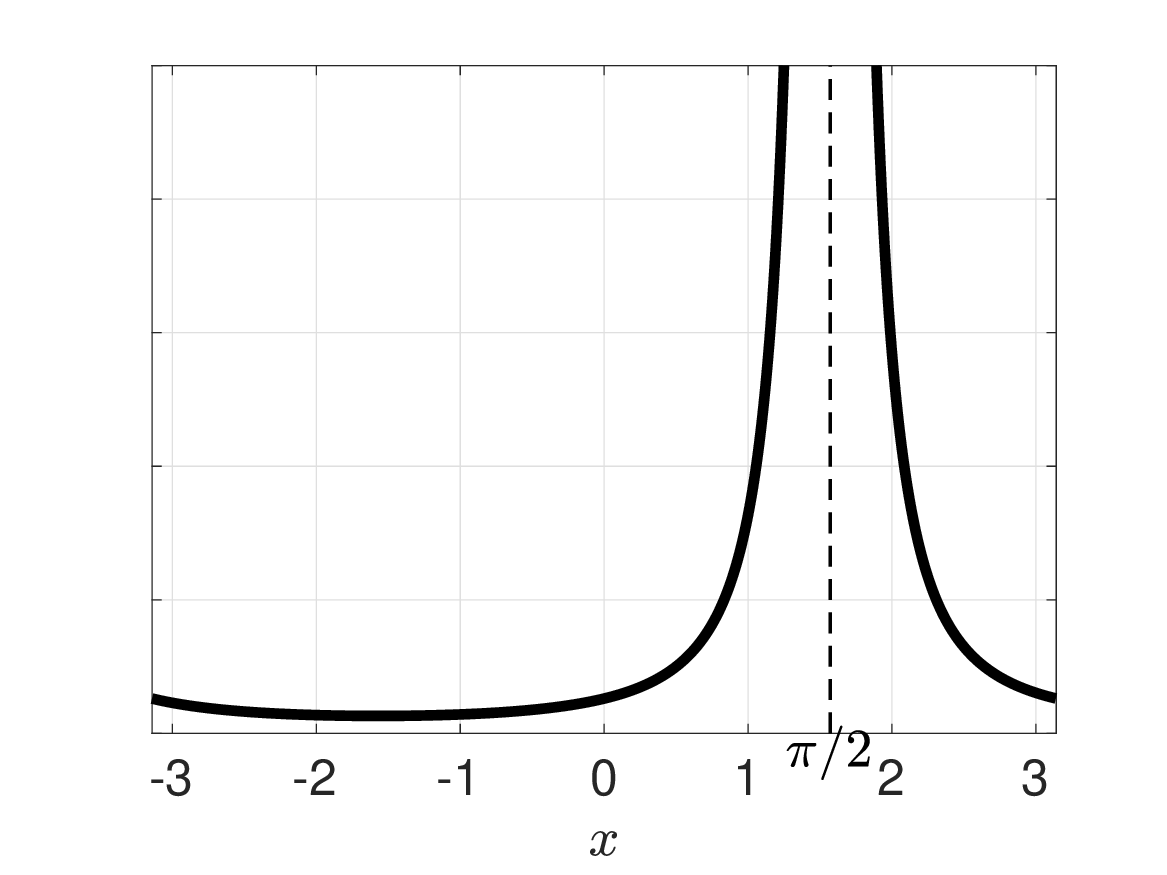}
    %\caption{$t=0.8$}
  \end{subfigure}
  \begin{subfigure}[b]{0.24\linewidth}
    \includegraphics[width=1.1\linewidth]{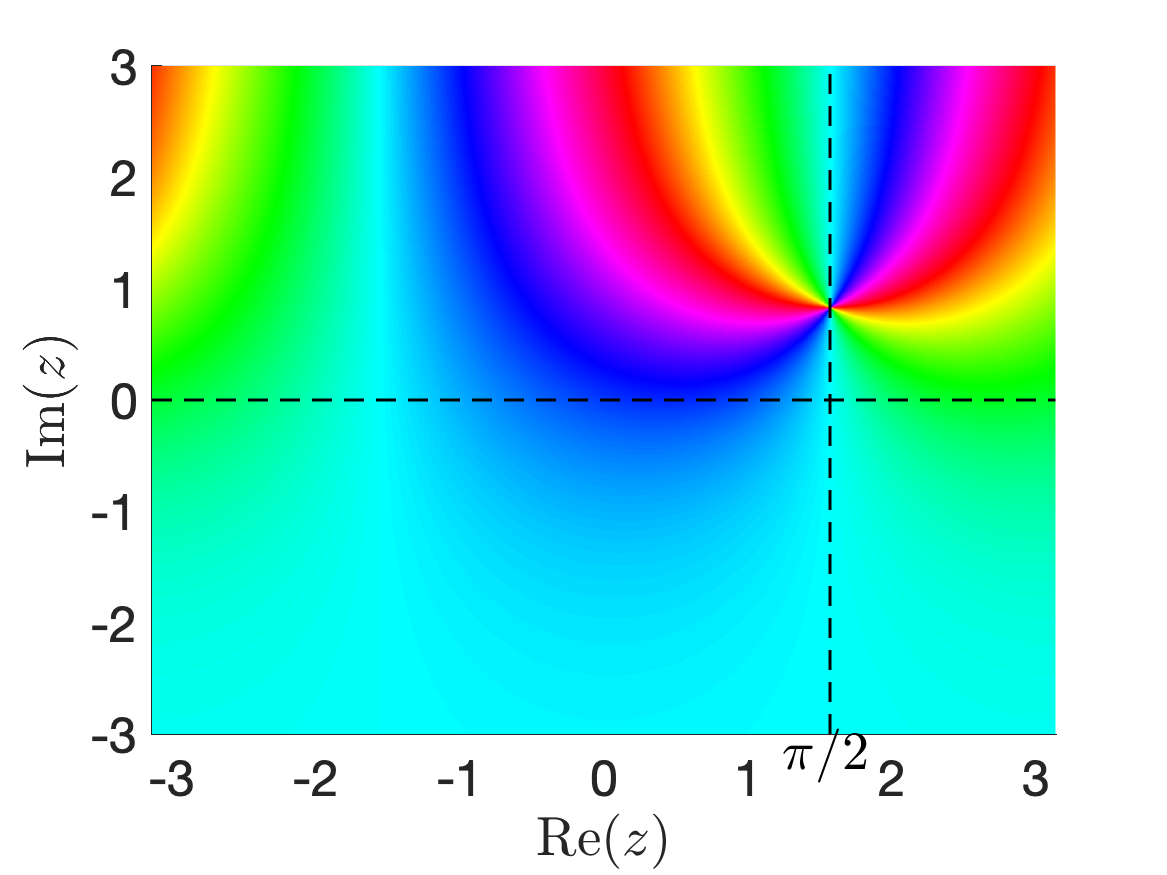}
   % \caption{$t=0.01$}
  \end{subfigure}
   \begin{subfigure}[b]{0.24\linewidth}
    \includegraphics[width=1.1\linewidth]{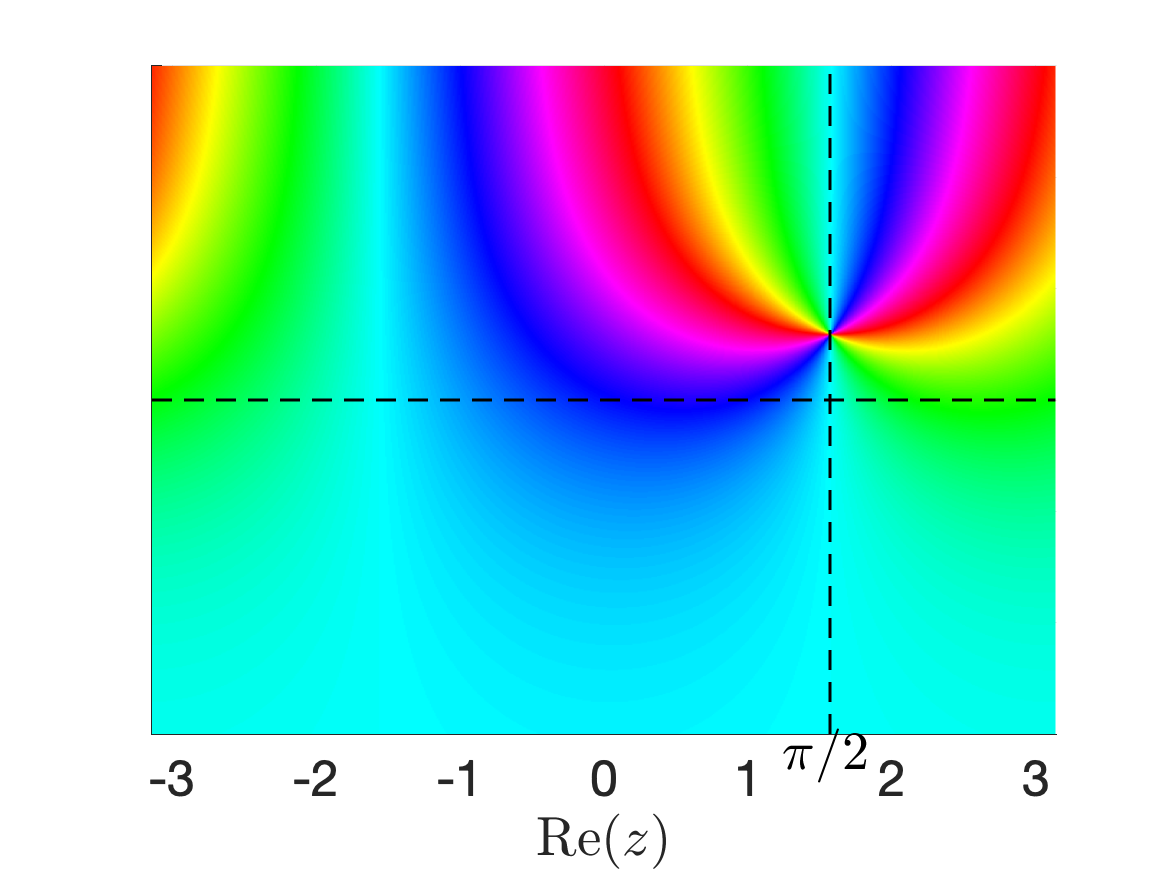}
   % \caption{$t=0.3$}
  \end{subfigure}
  \begin{subfigure}[b]{0.24\linewidth}
    \includegraphics[width=1.1\linewidth]{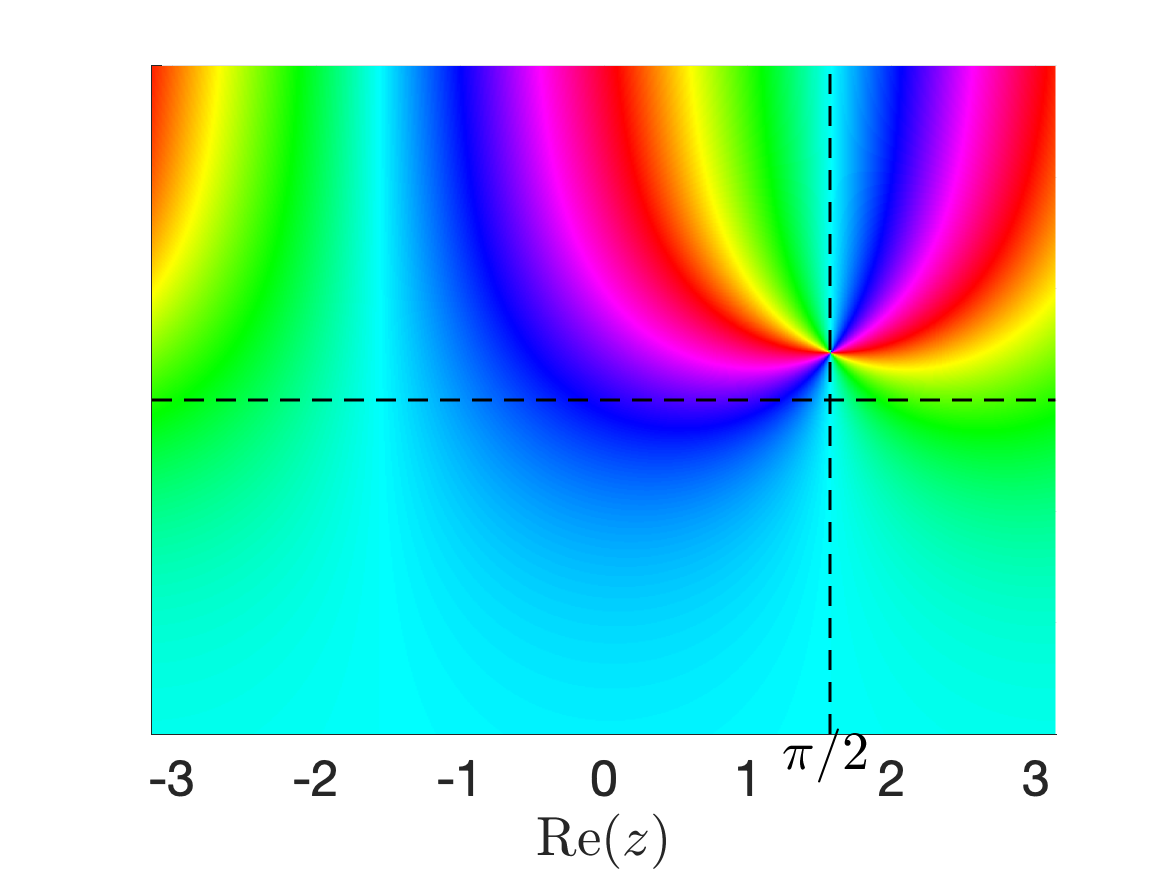}
    %\caption{$t=0.5$}
  \end{subfigure}
  \begin{subfigure}[b]{0.24\linewidth}
    \includegraphics[width=1.1\linewidth]{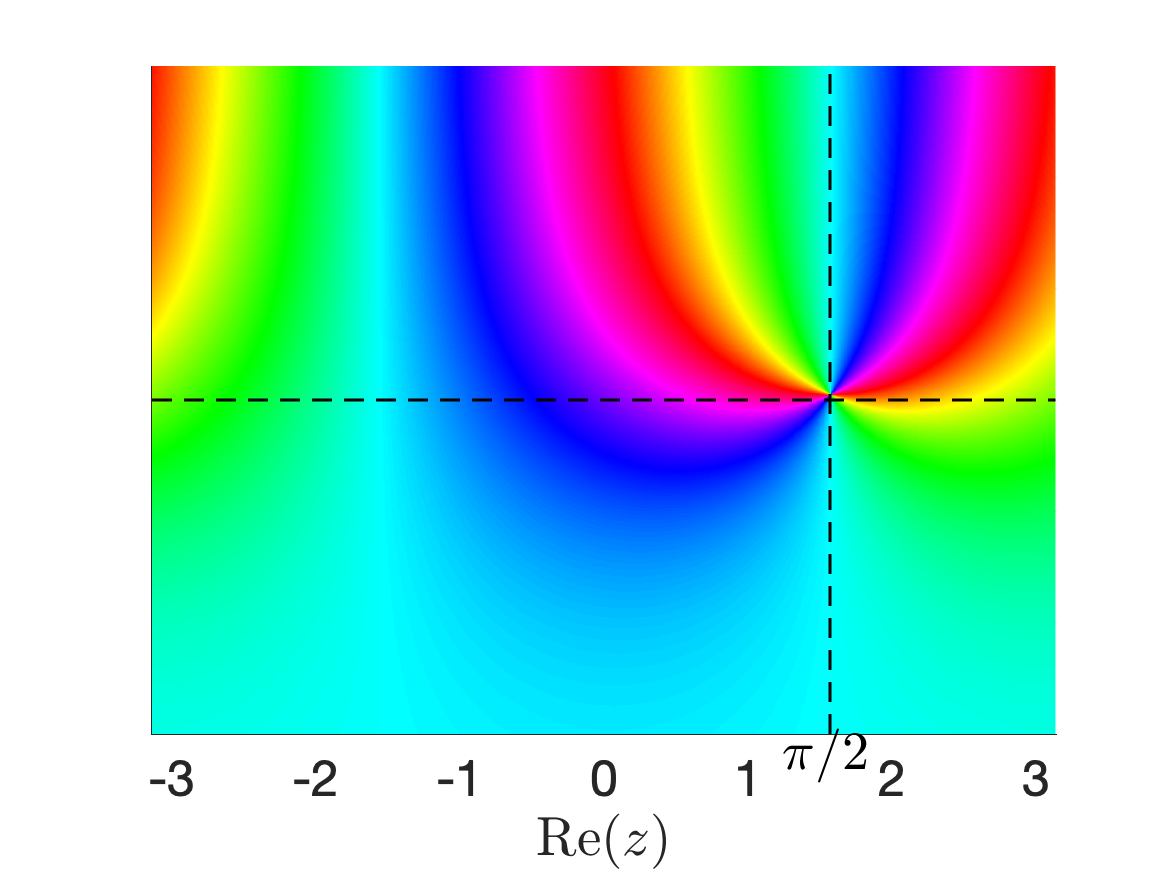}
   % \caption{$t=0.95$}
  \end{subfigure}
  \begin{subfigure}[b]{0.24\linewidth}
    \includegraphics[width=1.1\linewidth]{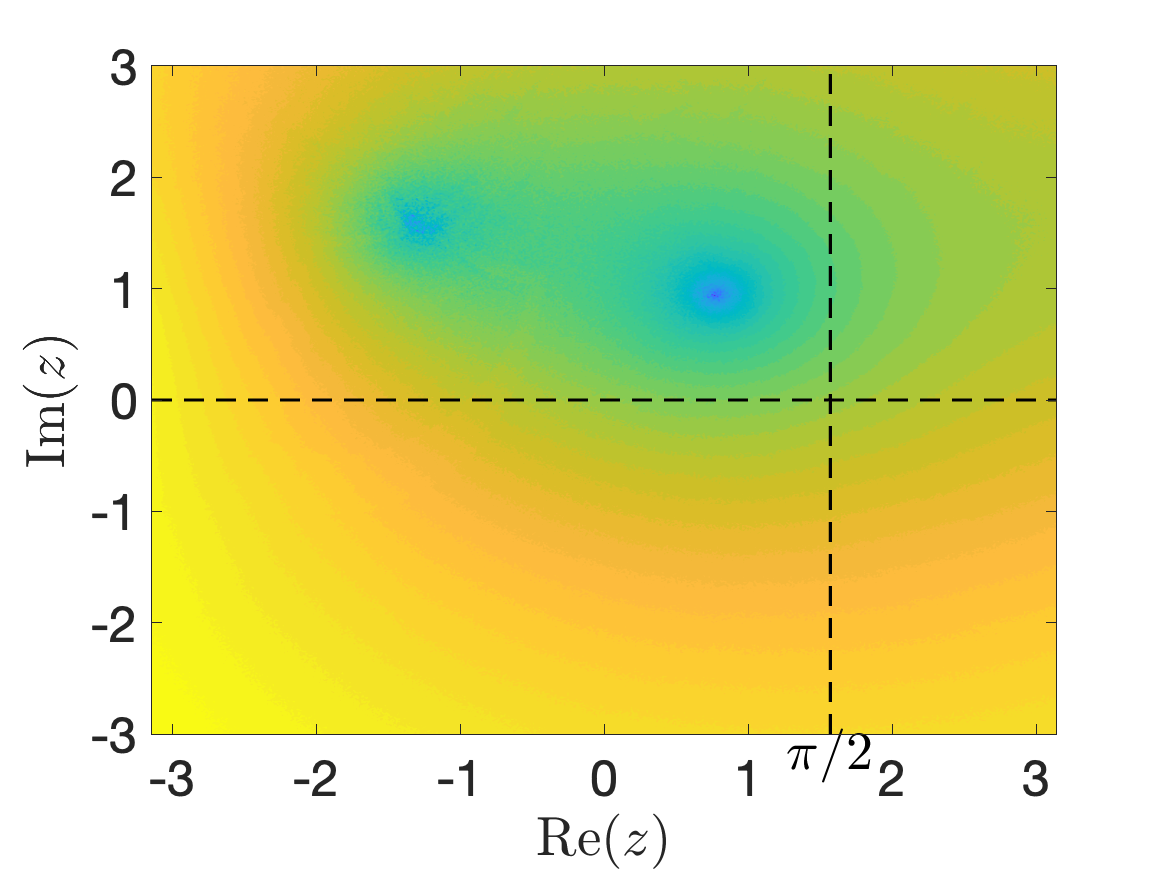}
    \caption{$t=0.01$}
  \end{subfigure}
  \begin{subfigure}[b]{0.24\linewidth}
    \includegraphics[width=1.1\linewidth]{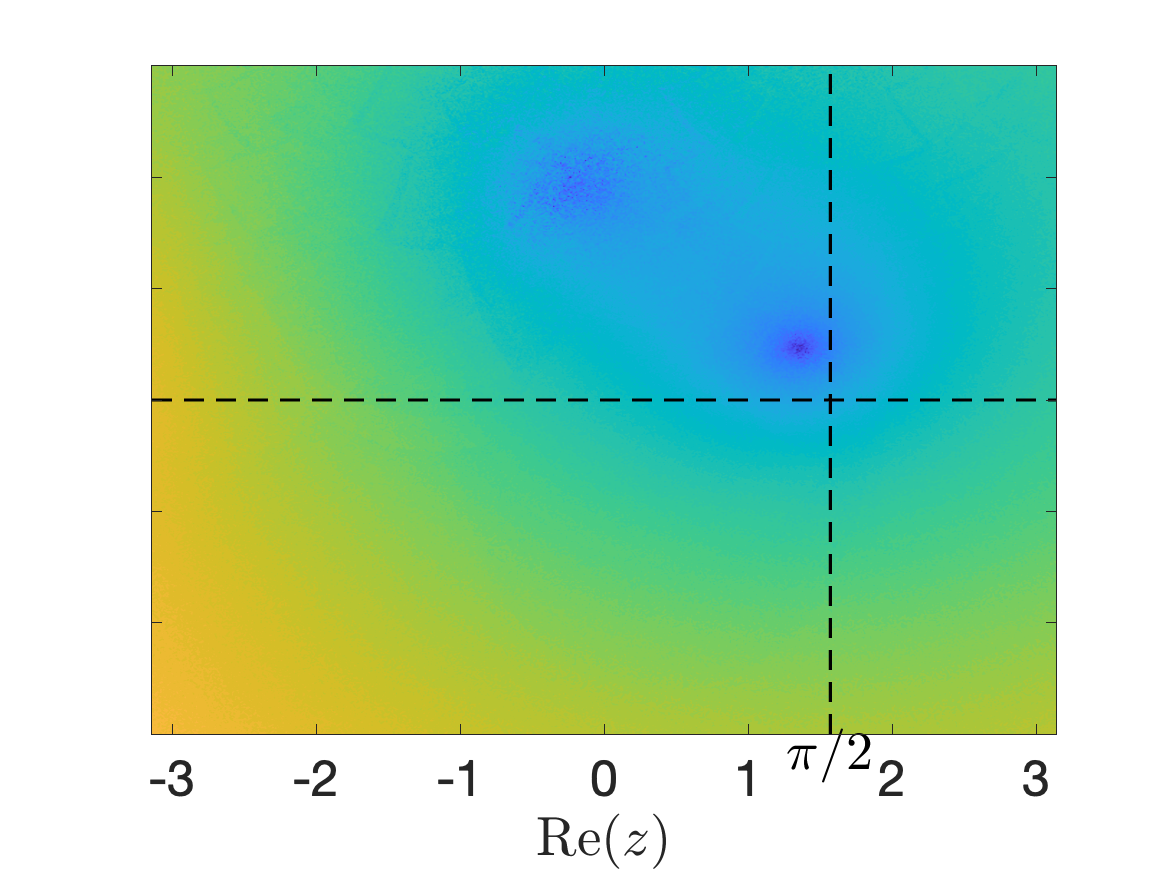}
    \caption{$t=0.3$}
  \end{subfigure}
  \begin{subfigure}[b]{0.24\linewidth}
    \includegraphics[width=1.1\linewidth]{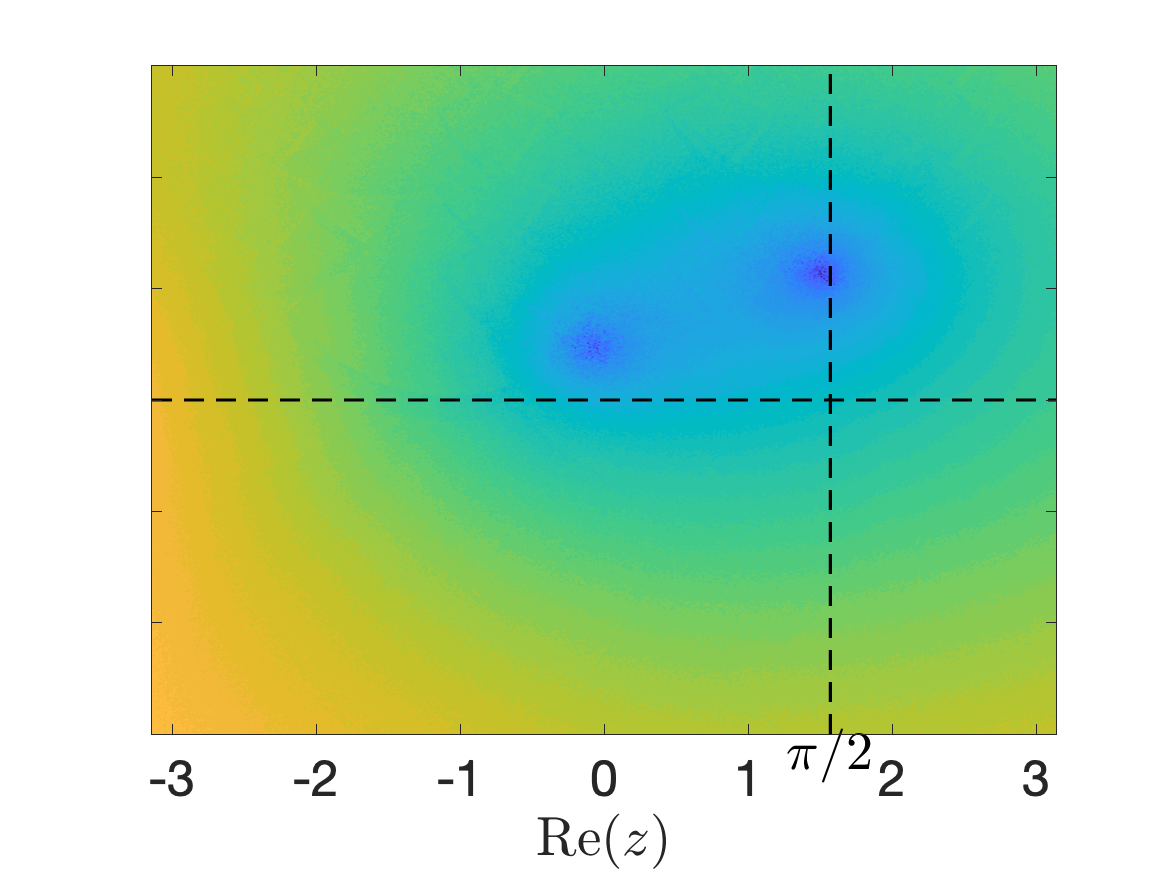}
    \caption{$t=0.5$}
  \end{subfigure}
  \begin{subfigure}[b]{0.24\linewidth}
    \includegraphics[width=1.1\linewidth]{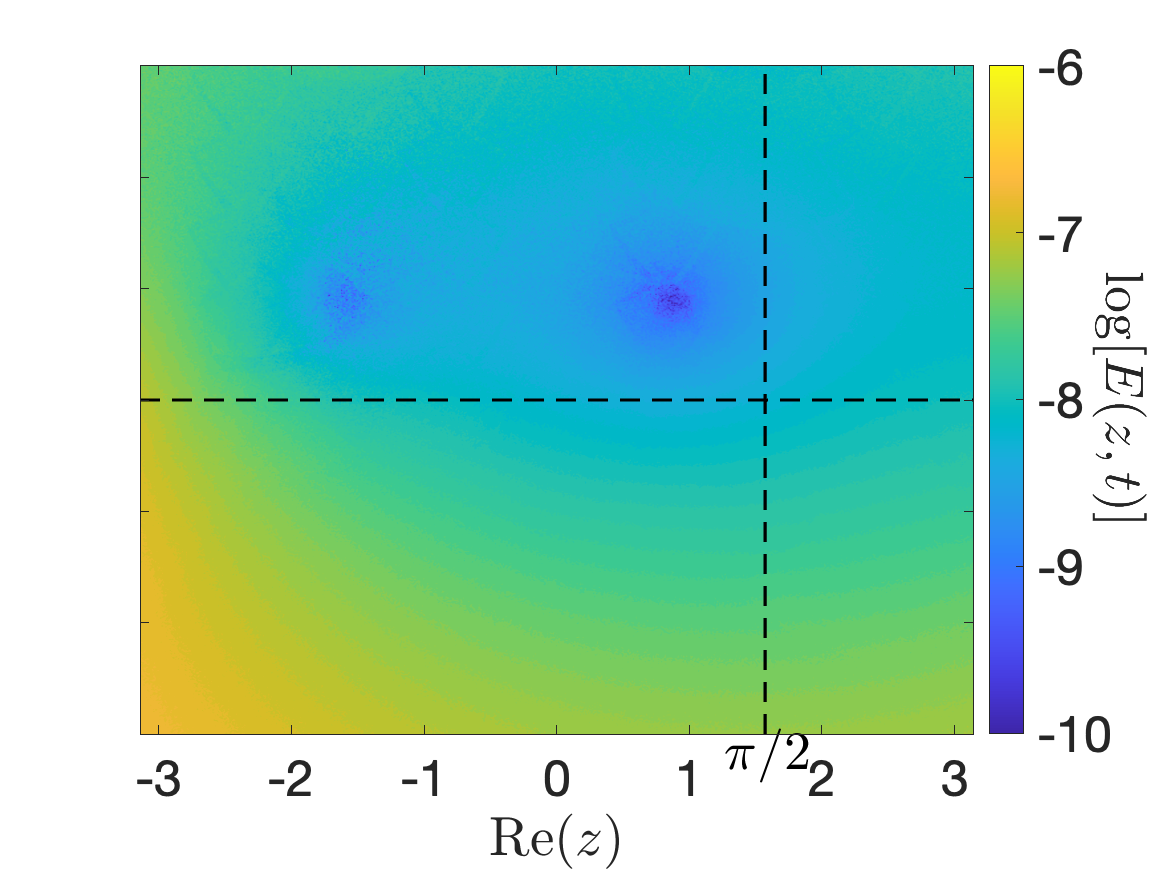}
    \caption{$t=0.95$}
  \end{subfigure}
  \caption{The solution $u(x,t)$ (first row),  the phase portraits of its  numerical analytic continuation $\Phi(z,t)$ as in (\ref{nnapr}) (second row), and the corresponding pole-scaled approximation error  $E(z,t)$ as in (\ref{errorrel}) in the logarithmic scale  (third row) of the NLH (\ref{heateq}) subject to the initial condition $u(x,0)$ as in (\ref{heatsol}) with parameters  $\sigma=\varrho=1$ for time $t=0.01, \, 0.3, \, 0.5, \, 0.95$. The blow up of the solution occurs at time $t_b=1$. }
  \label{figheat1}
\end{figure}

 We choose parameters $\sigma=\varrho=1$  and investigate the extended solution $u(z,t)$ for  $z \in D= [-\pi, \pi)\times  \mathrm{i} \, (-3, 3)$ and $t \in [0,1)$. In the domain $D$, $u(z,t)$ has one pole  with double multiplicity  of the form  
 \begin{equation}\label{sinc}
     z(t)= \frac{\pi}{2}+5\mathrm{i} \, (1- t)/6
 \end{equation}
and the blow up time is $t_b=1$ (see Figure \ref{figheat1}). Next, we construct the numerical analytic continuation $\Phi(z,t)$ of the solution $u(x,t)$  by  (\ref{nnapr}) and compute the estimated location of its singularities by  (\ref{spl})--(\ref{spl1}). As input, we use  $2n$  samples $u\left( \frac{ \pi j}{ n} , 0 \right)$, $j=-n,\dots, n-1$, of the initial condition (\ref{heatsol}) and compute the solution $u(x,t) $ in the form of the truncated Fourier series (\ref{appsr}) by the Fourier spectral method. For this example, the system of ODEs  (\ref{odes}) has the form  
%\begin{equation}\label{odesheat}
$$
 \frac{\mathrm{d} \hat{u}_k(t)}{\mathrm{d} t } + \hat{u}_k(t)(k^2/6+1) + \sum\limits_{\ell=-n}^n \hat{u}_\ell(t) \, \hat{u}_{k-\ell}(t)=0  , \quad k=-n,\dots,n.
 $$
 %\end{equation}

 \begin{figure}[h!]
    \centering
      \begin{subfigure}[b]{0.45\linewidth}
       \includegraphics[width=0.92\linewidth]{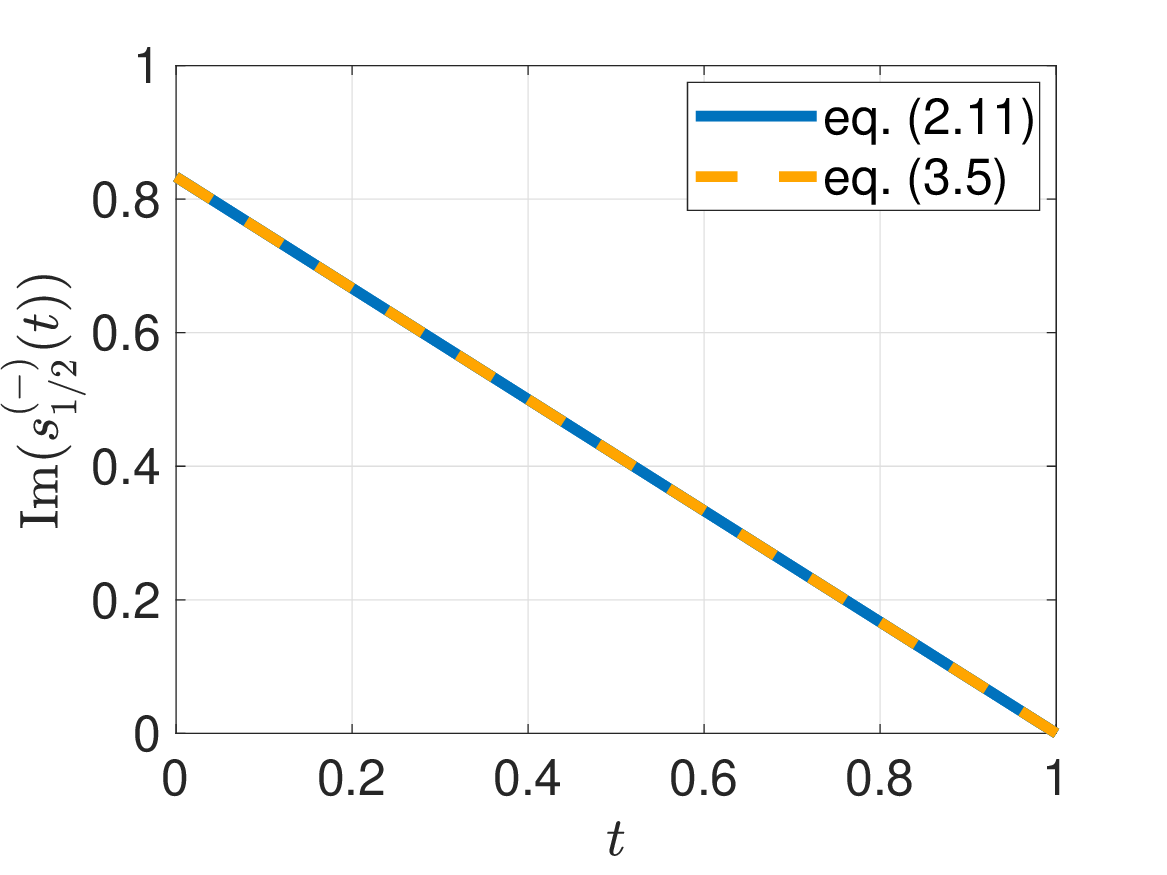}
     % \includegraphics[width=0.92\linewidth]{heatdunc.eps}
%    \includegraphics[width=0.92\linewidth]{sundynh112.eps}
   % \caption{$t=0.15$}
  \end{subfigure}
  \begin{subfigure}[b]{0.45\linewidth}
    \includegraphics[width=0.92\linewidth]{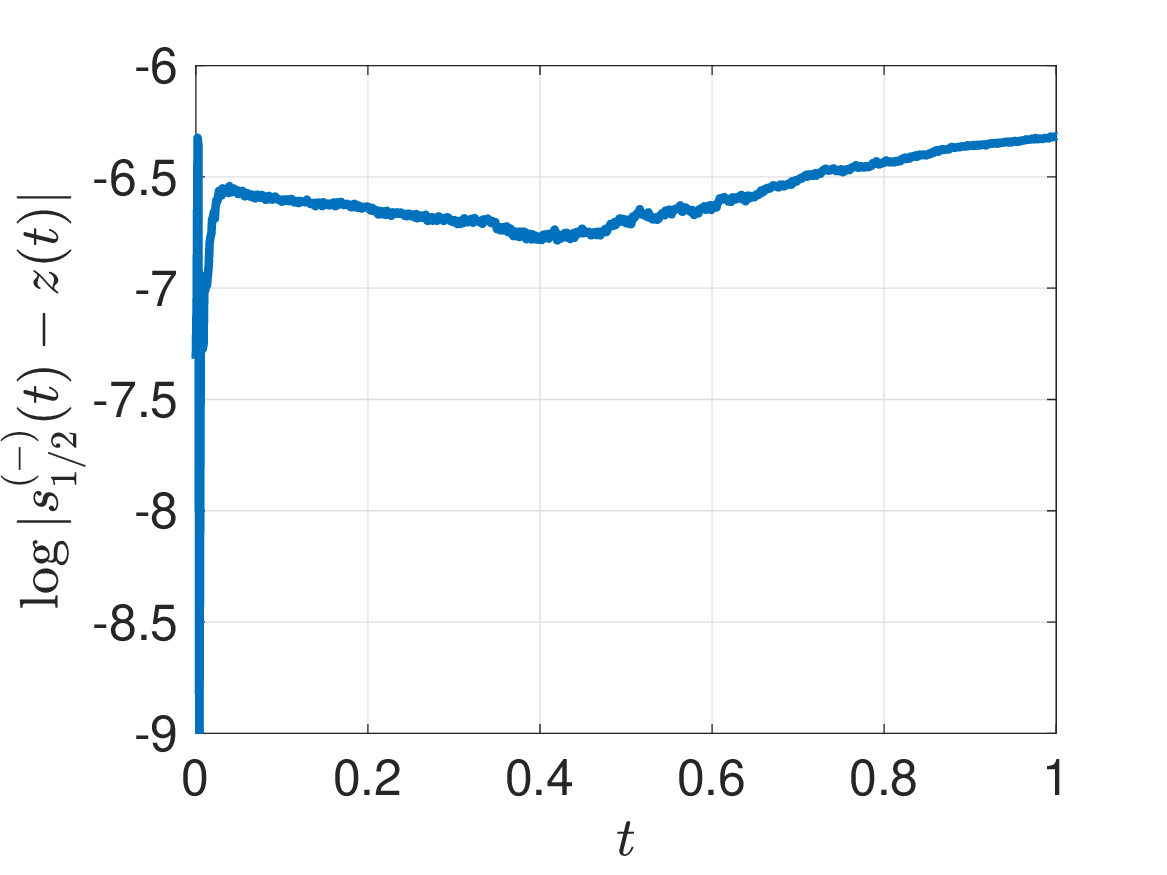}
   % \caption{$t=0.503$}
  \end{subfigure}
    %\caption{$t=0.8$}
      \caption{ Imaginary part of the pole trajectory in the complex plane of the analytic continuation of the solution of the  NLH (\ref{heateq}) subject to the initial condition $u(x,0)$ as in (\ref{heatsol}) with parameters  $\sigma=\varrho=1$  computed by the formula (\ref{spl1})  and  (\ref{sinc})   (left) and  the corresponding absolute approximation error $|s_{1/2}^{(-)}(t)-z(t)|$ in the logarithmic scale (right). }
  \label{sindynfig}
\end{figure}

 The SVD-based procedure from \cite{GGT13}  (Algorithm 2 in \cite{DKD2025}) allows us to compute  the optimal values for the parameters $N^{(-)}(t)=M^{(-)}(t)=2$ and $N^{(+)}(t)=M^{(+)}(t)=0$, which means that $\Phi(z,t)=\Phi^{(-)}(z,t)$. The component $\Phi^{(-)}$ detects two poles, $s_{1}^{(-)}(t)$ and $s_{2}^{(-)}(t)$, according to \eqref{spl1}. For each $0\le t<1$, these poles have identical numerical values, which coincide with those computed using formula \eqref{sinc} (Figure~\ref{sindynfig}, left) with the absolute approximation error $\approx 10^{-6}$ (Figure~\ref{sindynfig}, right).
 %The component $\Phi^{(-)}(z,t)$ detects two poles $s_{1/2}^{(-)}(t)$ according to (\ref{spl1}) with equal numerical values  for each $0 \leq t< 1$, which coincide with  the values computed by the   formula  (\ref{sinc})  (Figure   \ref{sindynfig}, left) and give the absolute approximation error $\approx 10^{-7}$  (Figure   \ref{sindynfig}, right).
 We present the  neural network approximation of the extended solution $\Phi(z,t)$ in Figure \ref{figheat1} (second row) and the corresponding pole-scaled approximation error $E(z,t)$ (\ref{errorrel})  in Figure \ref{figheat1} (third row) for  time  $t=0.01, \, 0.3, \, 0.5, \, 0.95$.
Furthermore, Figure \ref{figheat1} (second row)  enables us to observe the complex singularity of the extended solution as it approaches the real axis. According to (\ref{sinc}), the pole $z(t)$ moves linearly in time in the complex plane toward the real axis along the vertical line $x=\pi/2$ (Figure~\ref{sindynfig}, left).
%According to  (\ref{sinc}), the pole $z(t)$ moves in the complex plane toward the real axis along the line $x=\pi/2$  linearly with  respect to time (Figure \ref{sindynfig}, left).
At time $t_b=1$ it hits the real axis (Figure \ref{figheat1}, second row (d))  and  the solution $u(x,t)$ blows up (Figure \ref{figheat1}, first row (d)).  
Figure \ref{figheat1} (second row) demonstrates that the pole has multiplicity two, as the phase colors wrap twice around the color wheel during a single turn around the singularity.
%The double multiplicity of the pole is demonstrated in Figure \ref{figheat1} (second row)  by the fact that the colors go twice around the color wheel when encircling the singularity. 

Taking into account formula (\ref{spl1}), the weights $w_{j1}^{(-)}(t)$ and biases $b_{j1}^{(-)}(t)$ for $j=1,2$ are expected to vary linearly with time. This prediction is confirmed by the numerical experiments shown in Figure \ref{figwbh}.
%Taking into account  the formula (\ref{spl1}), weights $w_{j1}^{(-)}(t)$ and biases $b_{j1}^{(-)}(t)$ for $j=1,2$  must change linearly with respect to time, which is confirmed by  numerical experiments (see Figure   \ref{figwbh}). 
Figure \ref{figwbh} also illustrates the performance of our method for computing the weights and biases of the hidden layer under different choices of the parameter  $C_{10}^{(-)}(t)$ defined in (\ref{ccoef2}) for $t \in [0,1)$.
%In Figure   \ref{figwbh}, we also  demonstrate the performance of our method for the computation of weights and biases of the hidden layer in dependency on the choice    of parameters  $C_{10}^{(-)}(t)$   as in (\ref{ccoef2}) for $t \in [0,1)$.
Specifically, we consider two cases: the fixed given value  
 and a randomly chosen value of $C_{10}^{(-)}(t)$ for each $t \in [0,1)$.
 Note that  for a pole $z(t)$, $t \in [0,1)$,  with double multiplicity, we obtain two different neurons of the hidden layer of the neural network component $\Phi^{(-)}$, which can be  observed in  Figure  \ref{figwbh}. 
 Note again that  the choice of parameters  $C_{j0}^{(\pm)}(t)$ does not influence  the accuracy of the computation of  singularities. 

\begin{figure}[h!]
  \centering
  \begin{subfigure}[b]{0.45\linewidth}
    \includegraphics[width=1\linewidth]{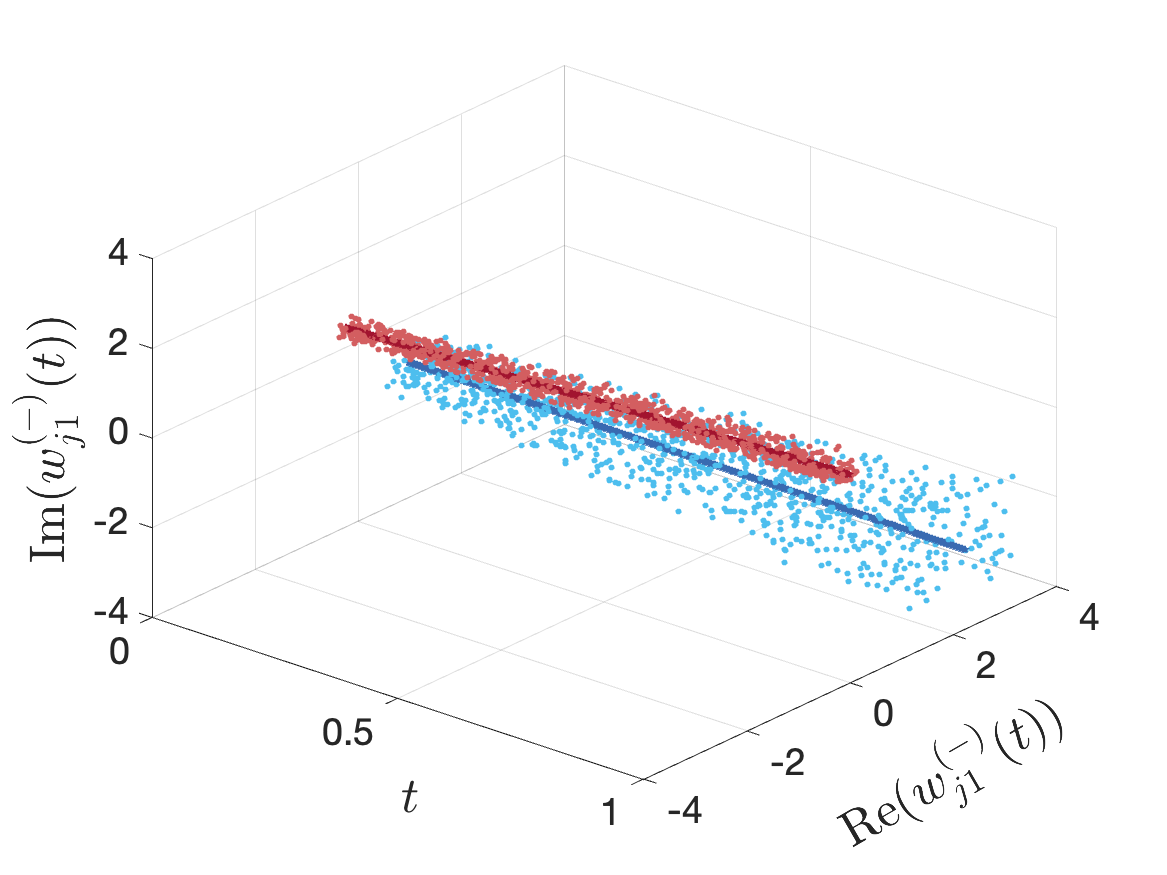}
   % \caption{$t=0.15$}
  \end{subfigure}
  \begin{subfigure}[b]{0.45\linewidth}
    \includegraphics[width=1\linewidth]{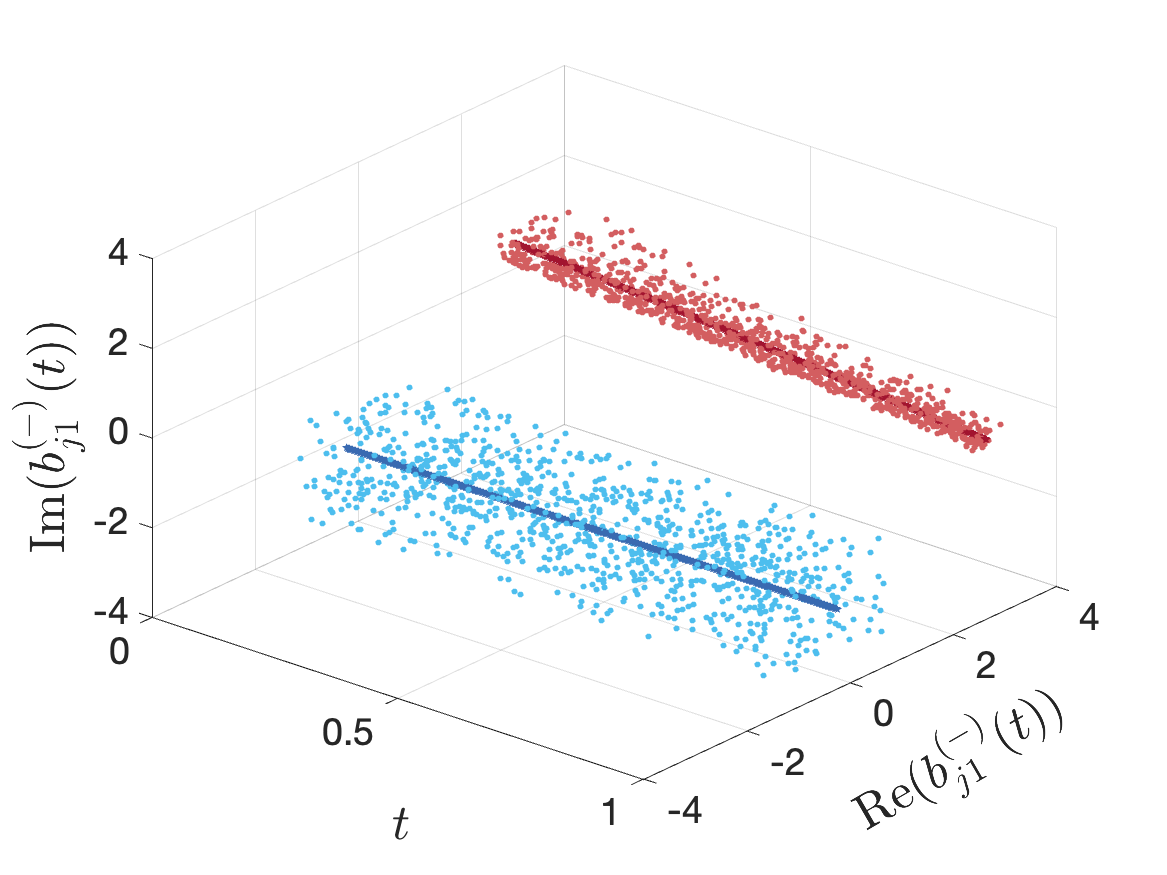}
   % \caption{$t=0.503$}
  \end{subfigure}
  \caption{Weights $w_{j1}^{(-)}(t)$ (left) and biases  $b_{j1}^{(-)}(t)$ (right) used in (\ref{spl1}) for computation of poles  $s_{j}^{(-)}(t)$, $j=1,2$ ($j=1$ -- blue, $j=2$ -- red) of the extended solution of the   NLH (\ref{heateq}) subject to the initial condition $u(x,0)$ as in (\ref{heatsol}) with parameters  $\sigma=\varrho=1$ over the time interval $[0,1)$. Light blue and light red -- weights and biases computed with random choice of parameters $C_{10}^{(-)}(t) \in [0.5,1.5]+ \mathrm{i} \, [0.5,1.5] $, dark blue and dark red --  weights and biases computed with the fixed parameters $C_{10}^{(-)}(t) = 1+\mathrm{i}$ for $t \in [0,1)$. }
  \label{figwbh}
\end{figure}

\end{example}

\begin{example}
\label{ex32}
Next, we consider one more example of the NLH equation
\begin{equation}\label{heat2}
    u_t-u_{xx}-u^{2}=0,
\end{equation}
with the initial condition 
\begin{equation}\label{inheat}
    u(x,0)=\alpha \, \cos(x),    
\end{equation}
where $\alpha>0$ and  $ x \in [-\pi,\pi)$.  Note that in contrast to Example \ref{ex31}, the initial condition  (\ref{inheat}) is an entire function after the extension to the complex plane. Thus, singularities of the analytic continuation $u(z,t)$ are ``born'' at infinity for small $t>0$. As discussed in  \cite{FW24, W03}, the solution of (\ref{heat2}) with the initial condition (\ref{inheat}) develops a finite time blow up.

We investigate the analytic continuation $u(z,t)$ for  $z \in D= [-\pi, \pi)\times  \mathrm{i} \, (-4, 4)$ and $t \in [0,t_b^{(\alpha)})$, where by $t_b^{(\alpha)}$ we denote the blow up time for $\alpha>0$. In the first step, the solution $u(x,t)$  as in (\ref{appsr}) is computed by the Fourier spectral method employing $2n$ function values  $u\left( \frac{ \pi j}{ n} , 0 \right)$, $j=-n,\dots, n-1$, of the initial condition (\ref{inheat}) with $n=200$.  The system of ODEs  (\ref{odes}) has the form 
%\begin{equation}\label{odesheat2}
 $$
 \frac{\mathrm{d} \hat{u}_k(t)}{\mathrm{d} t } + k^2 \, \hat{u}_k(t) - \sum\limits_{\ell=-n}^n \hat{u}_\ell(t) \, \hat{u}_{k-\ell}(t)=0  , \quad k=-n,\dots,n.
 $$
 %\end{equation}

\begin{figure}[h!]
  \centering
    \begin{subfigure}[b]{0.24\linewidth}
    \includegraphics[width=1.1\linewidth]{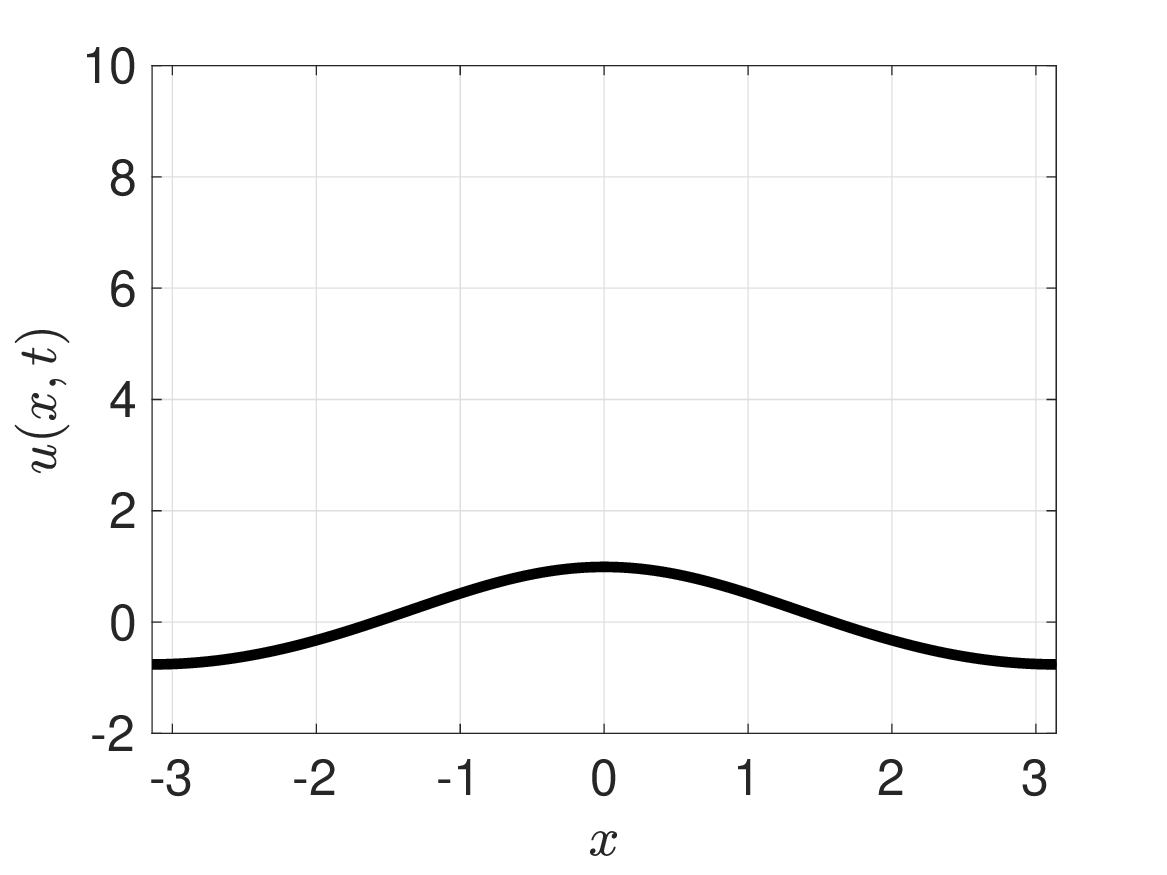}
    %\caption{$t=0.15$}
  \end{subfigure}
  \begin{subfigure}[b]{0.24\linewidth}
    \includegraphics[width=1.1\linewidth]{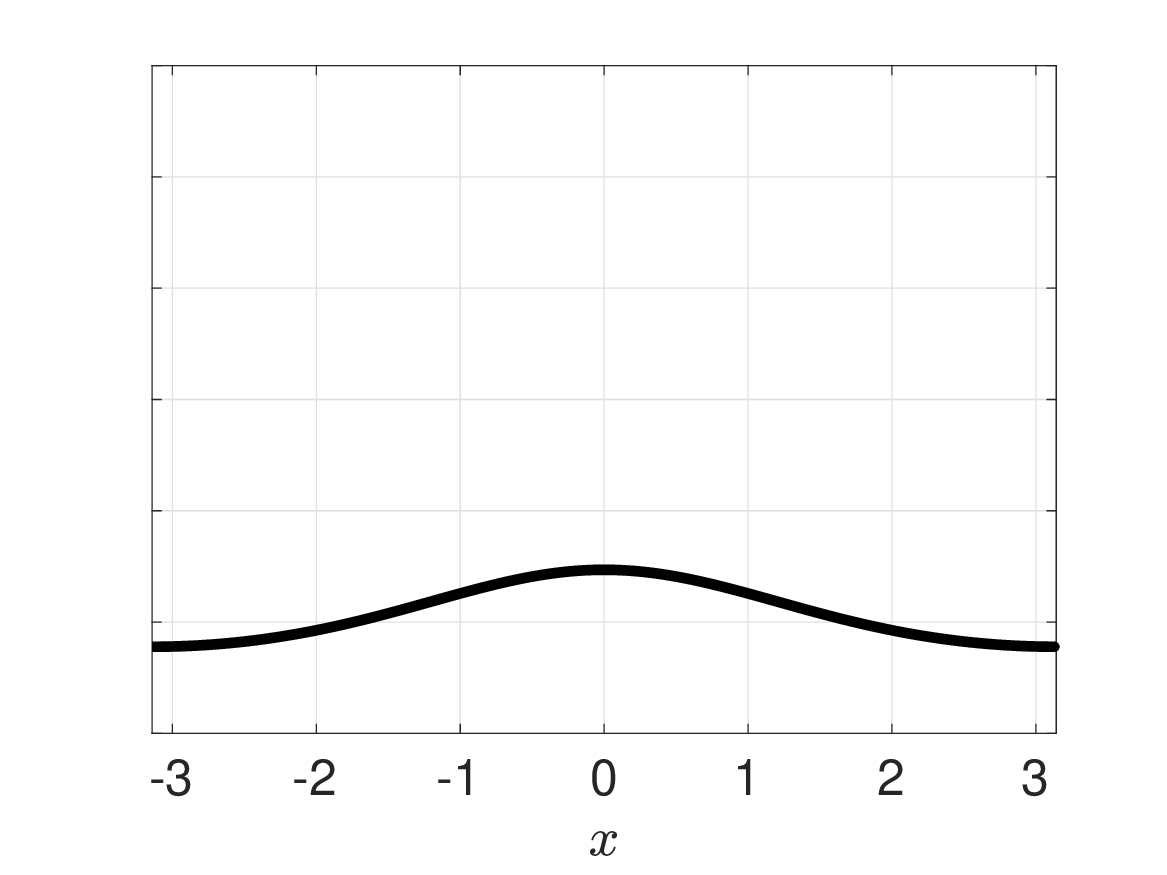}
   % \caption{$t=0.503$}
  \end{subfigure}
  \begin{subfigure}[b]{0.24\linewidth}
    \includegraphics[width=1.1\linewidth]{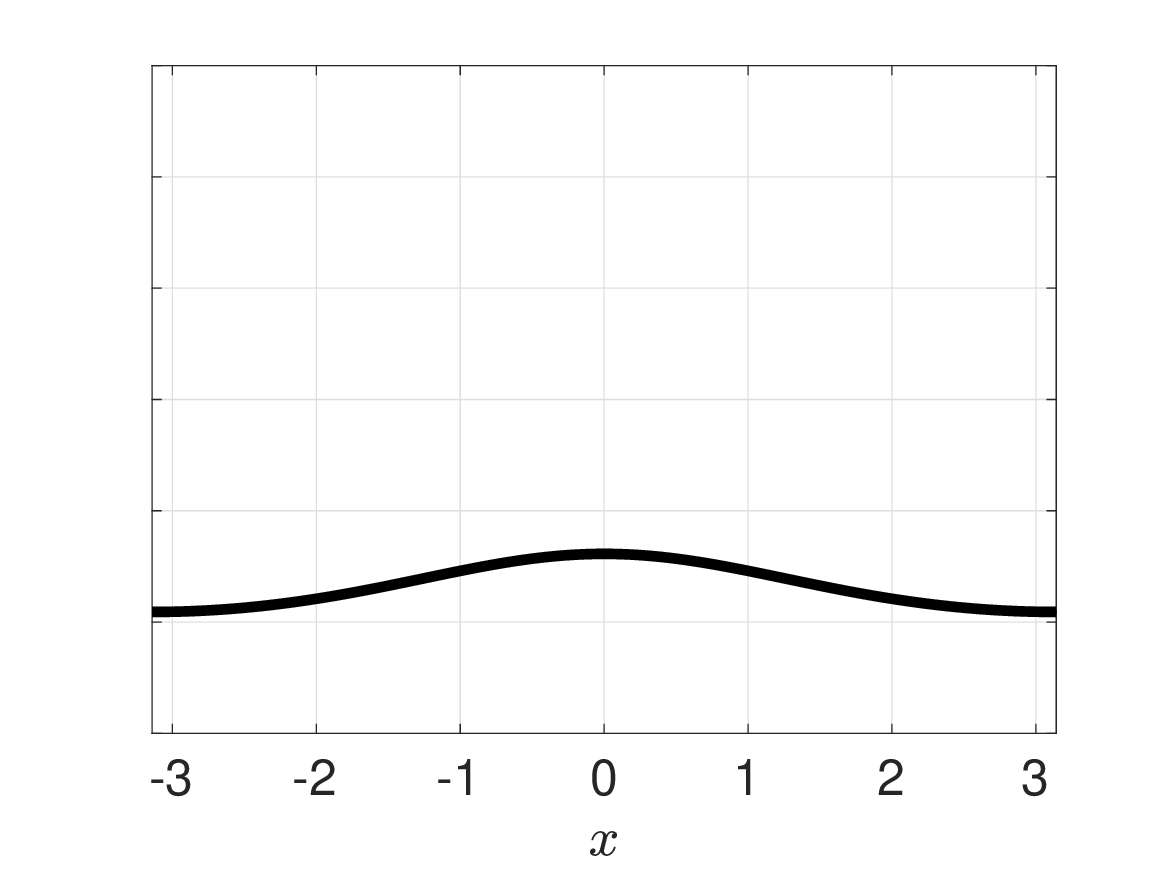}
    %\caption{$t=2.05$}
  \end{subfigure}
  \begin{subfigure}[b]{0.24\linewidth}
    \includegraphics[width=1.1\linewidth]{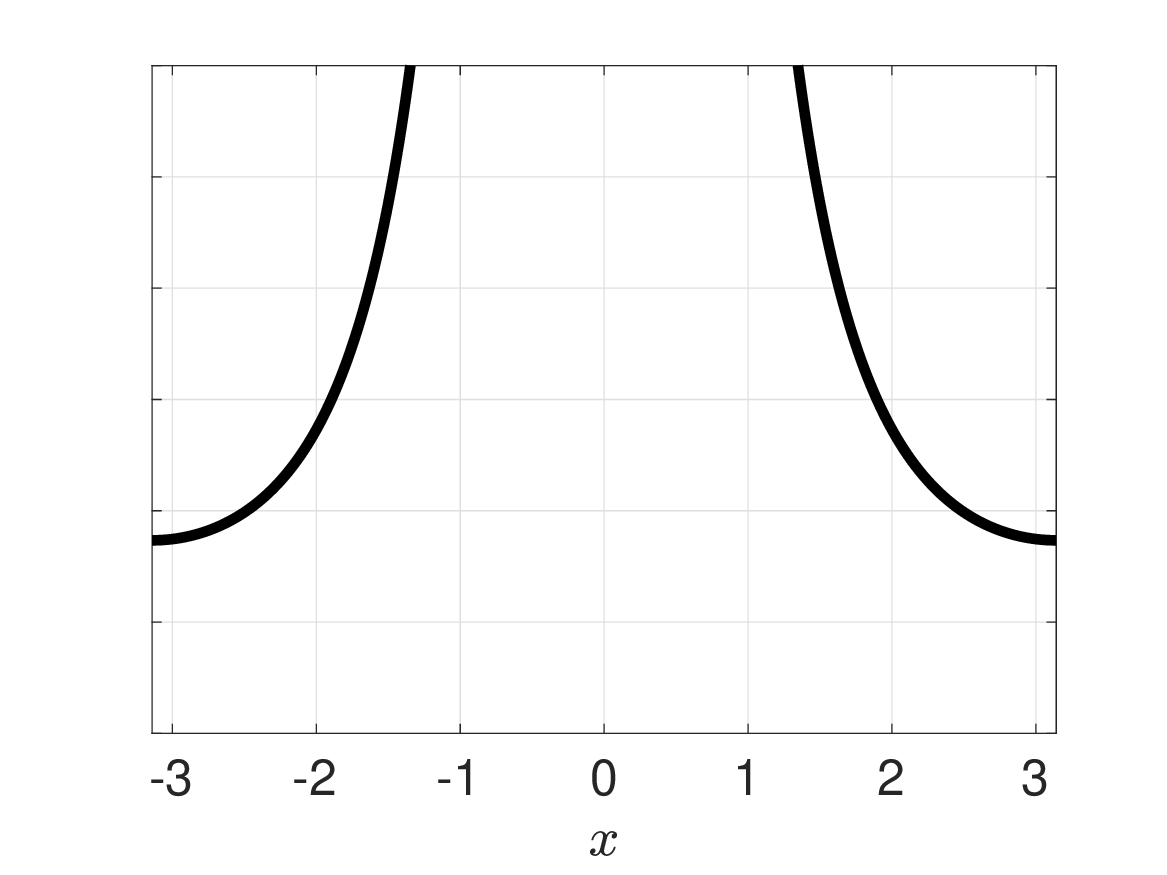}
   % \caption{$t=3.17$}
  \end{subfigure}
  \begin{subfigure}[b]{0.24\linewidth}
    \includegraphics[width=1.1\linewidth]{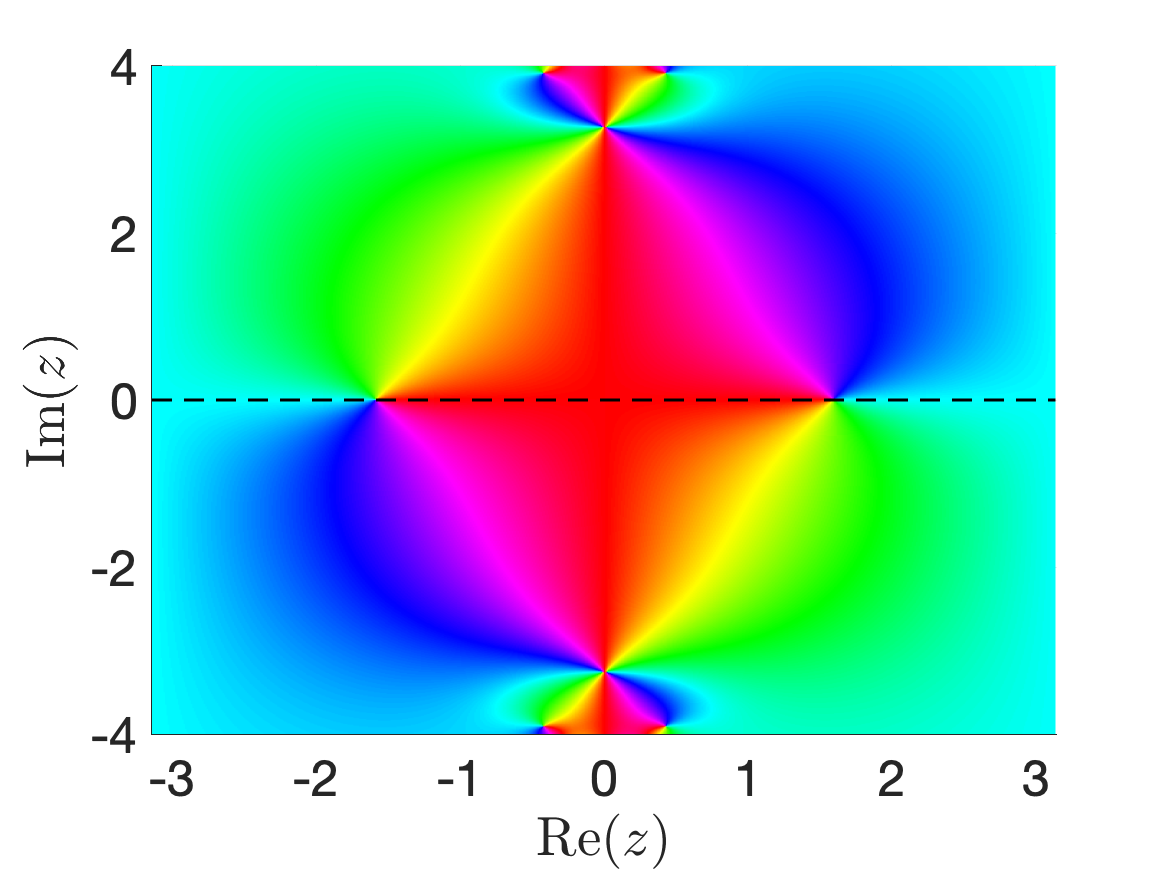}
    \caption{$t=0.15$}
  \end{subfigure}
  \begin{subfigure}[b]{0.24\linewidth}
    \includegraphics[width=1.1\linewidth]{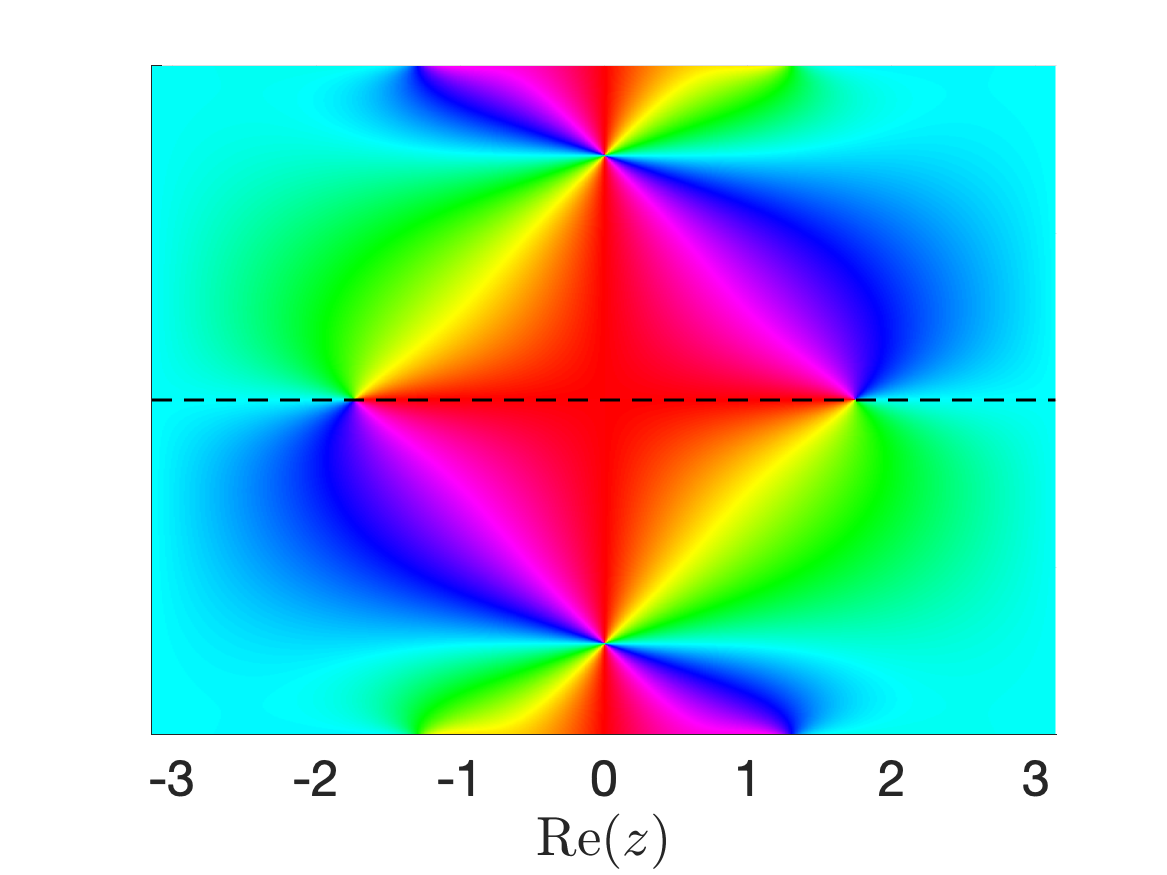}
    \caption{$t=0.503$}
  \end{subfigure}
  \begin{subfigure}[b]{0.24\linewidth}
    \includegraphics[width=1.1\linewidth]{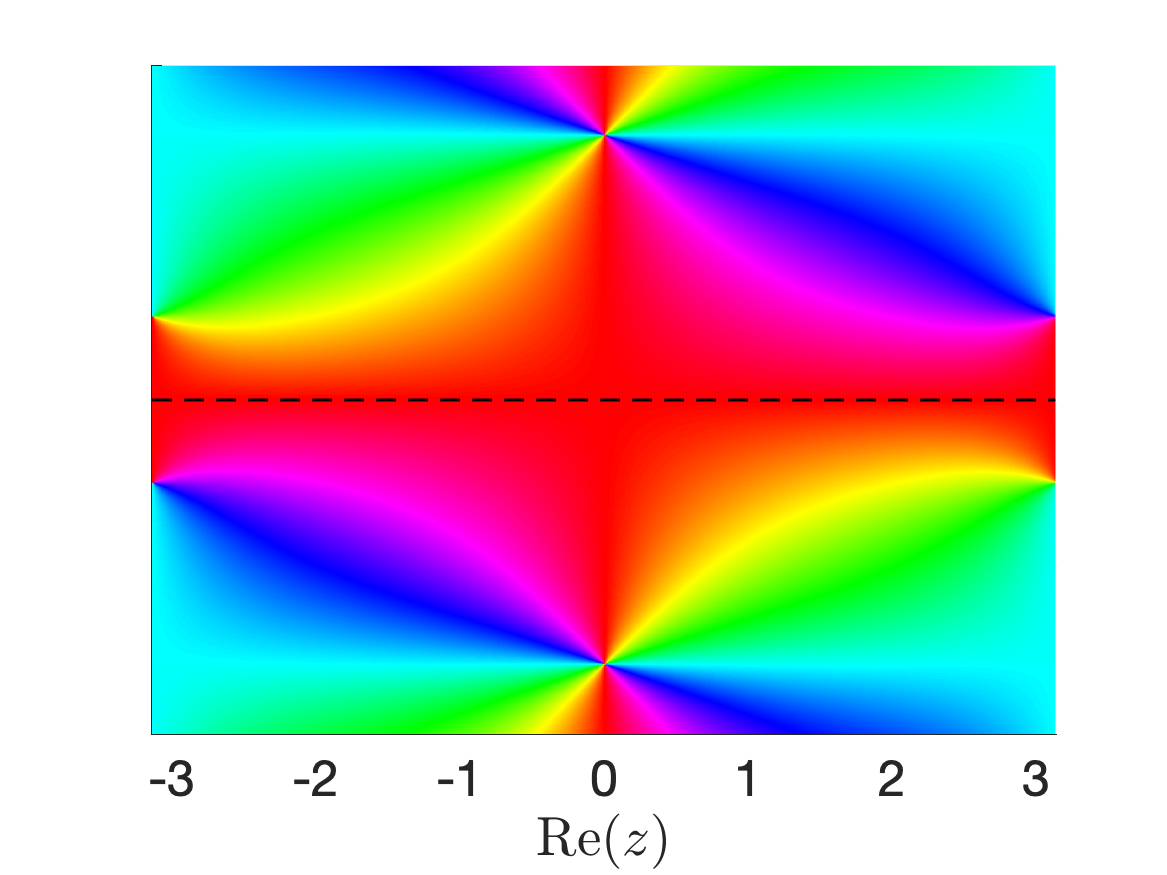}
    \caption{$t=2.05$}
  \end{subfigure}
  \begin{subfigure}[b]{0.24\linewidth}
    \includegraphics[width=1.1\linewidth]{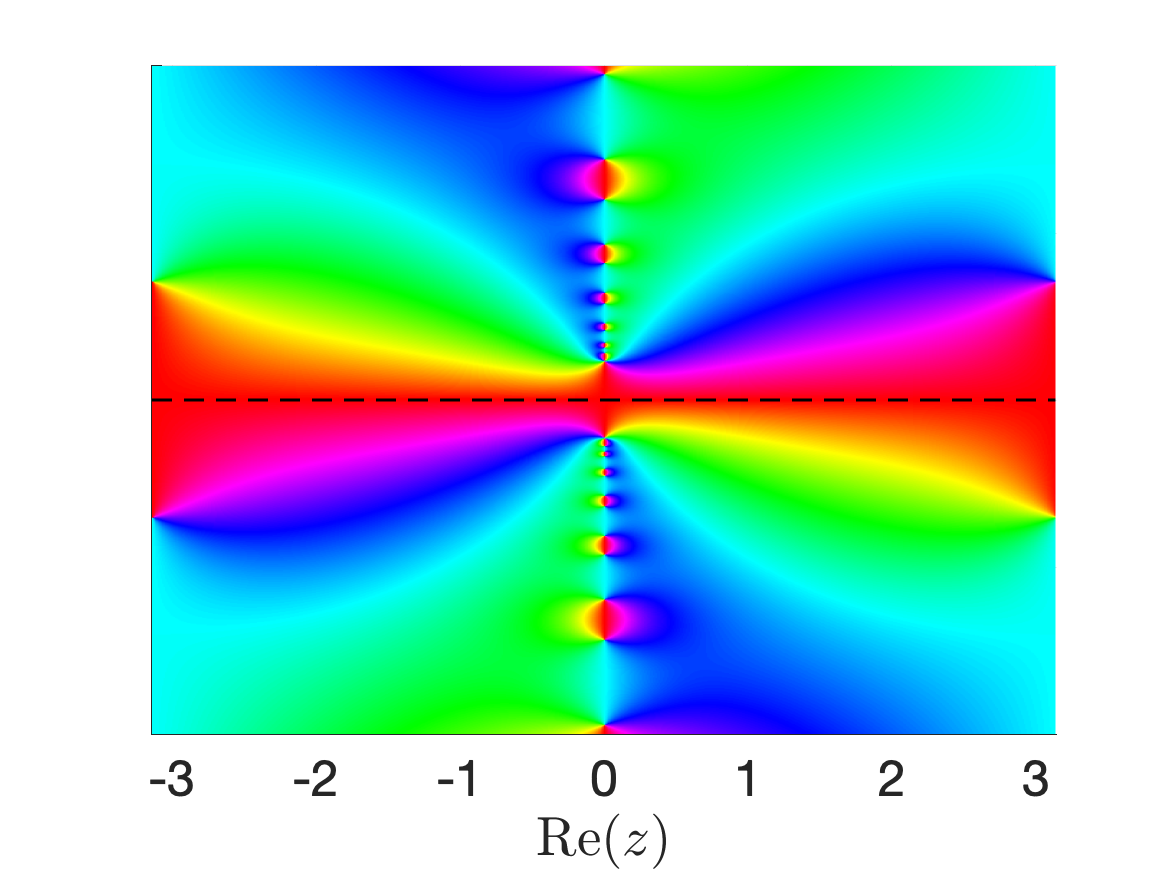}
    \caption{$t=3.17$}
  \end{subfigure}
  \caption{The solution $u(x,t)$ (first row) and  the phase portraits of its numerical analytic continuation $\Phi(z,t)$ as in (\ref{nnapr}) (second row) of the NLH (\ref{heat2}) subject to the initial condition (\ref{inheat}) with $\alpha=1$ for time $t=0.15, 0.503,2.05,3.17$.  The blow up occurs at time $t_b^{(1)}\approx 3.17395$. }
  \label{heatsolfig}
\end{figure}
 
 Next, we extend the numerical solution $u(x,t)$ to the complex plane by our neural network-based method described in Section \ref{le}. The numerical analytic continuation $\Phi(z,t)$ of the solution $u(x,t)$ is computed by (\ref{nnapr}) and presented in Figure \ref{heatsolfig} (second row) in case $\alpha=1$ for time $t=0.15, 0.503,2.05$ and $t=3.17$. As was investigated in \cite{FW24, W03},  singularities of the extended solution of the PDE (\ref{heat2}) subject to the initial condition  (\ref{inheat})
 have the complicated nature. Using local expansion of the solution around the singularity (see formula (9) in \cite{FW24}), it was concluded that the singularities are actually  branch points. Thus, using equations (\ref{spl})-(\ref{spl1}), we compute estimated locations $s_1^{(\pm)}(t)$ of the closest to the real axis singularities (in the upper and lower half-planes) and present their trajectories  with respect to time in Figure \ref{figheatsinwei} (first row) for $\alpha=1,10,1000$.

\begin{figure}[h!]
\centering
\begin{subfigure}{0.32\textwidth}
    \centering
    \includegraphics[width=1.1\linewidth]{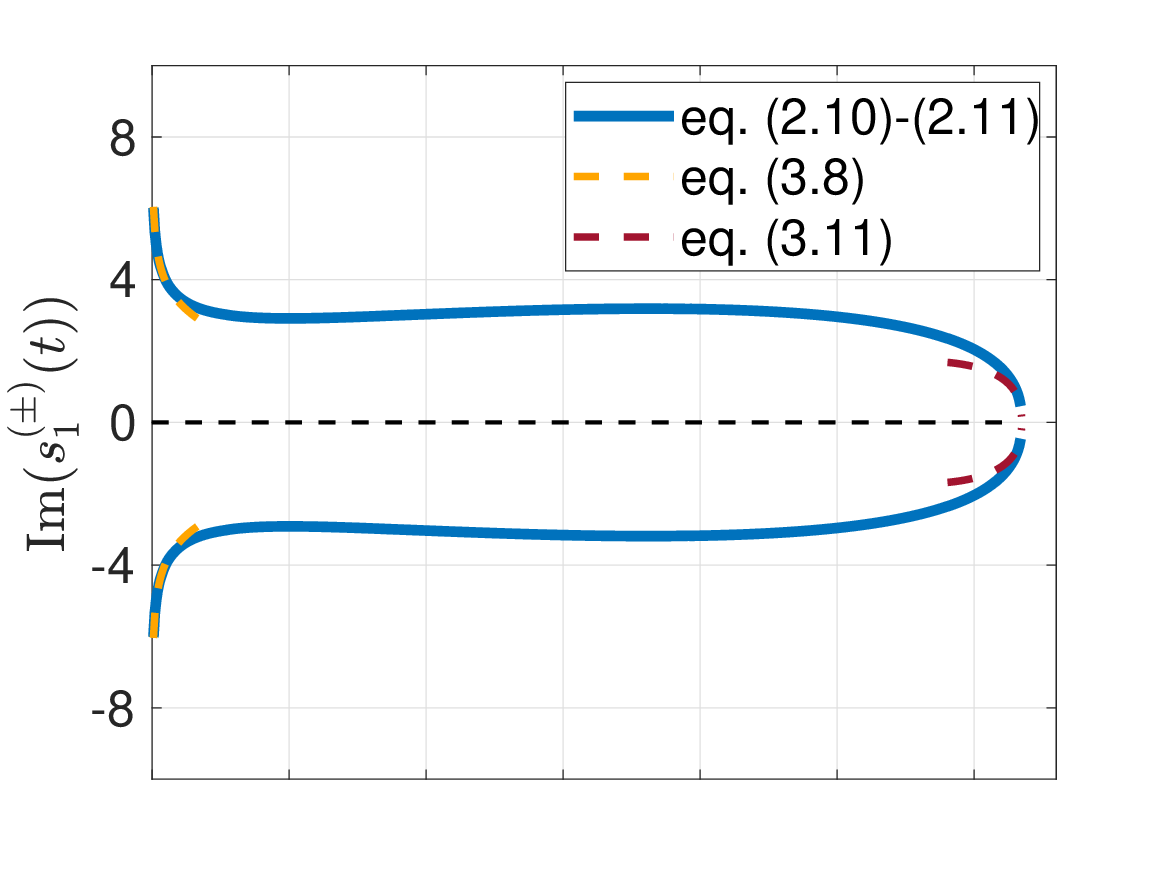}
\end{subfigure}\hfill
\begin{subfigure}{0.32\textwidth}
    \centering
     \includegraphics[width=1.1\linewidth]{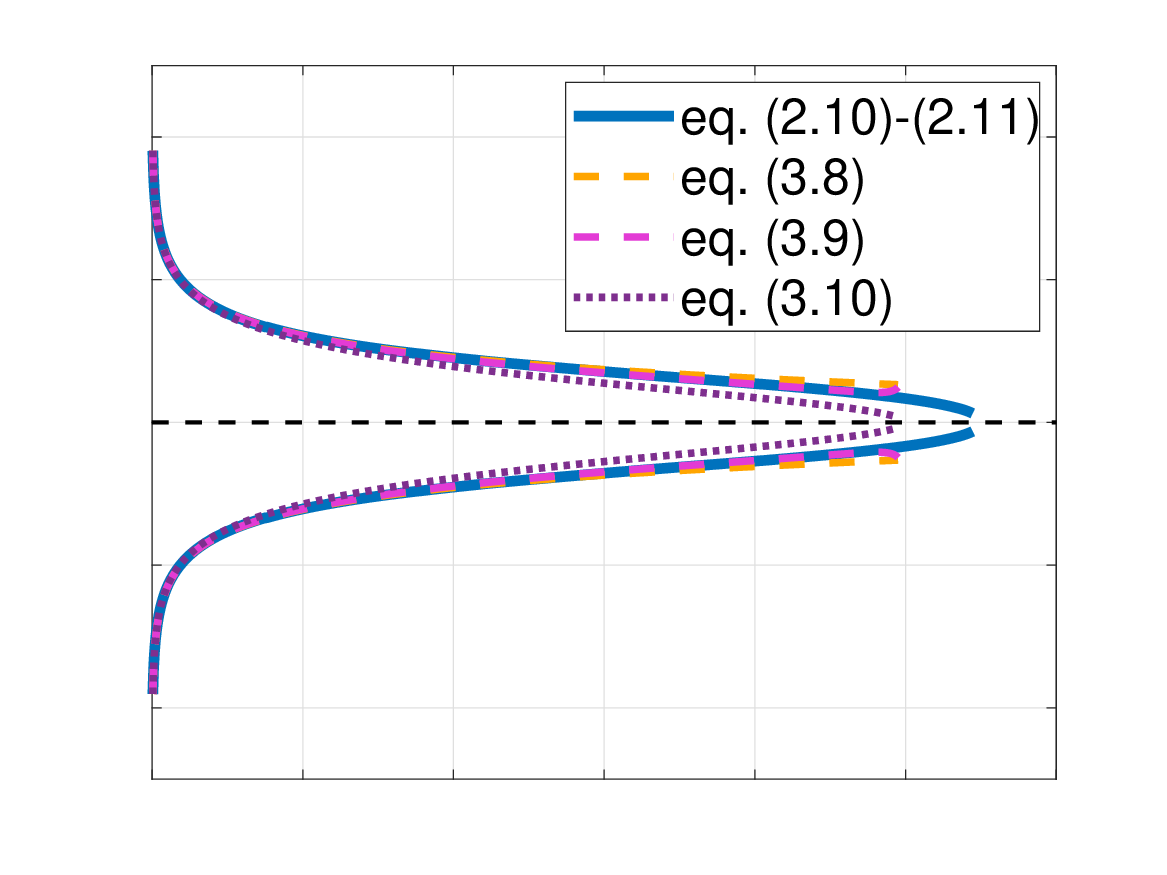}
\end{subfigure}\hfill
\begin{subfigure}{0.32\textwidth}
    \centering
    \includegraphics[width=1.1\linewidth]{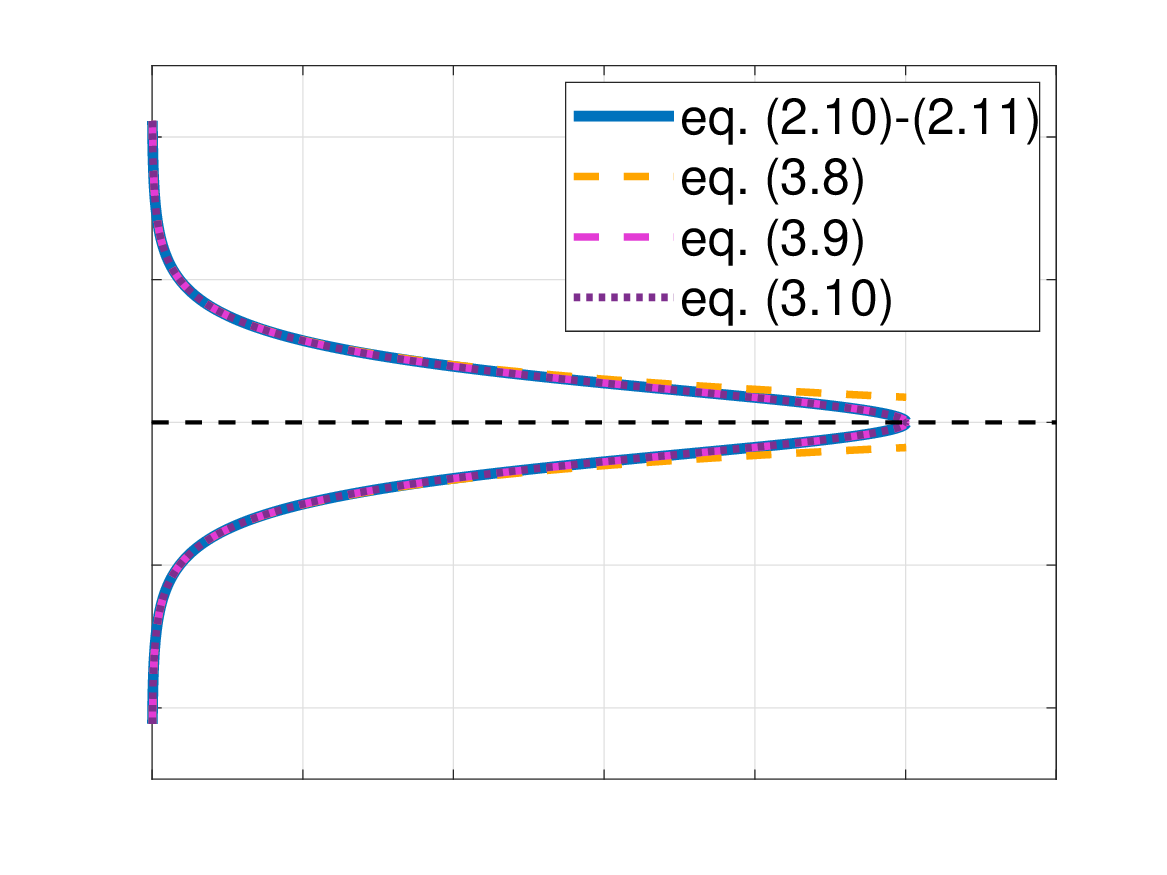}
\end{subfigure}
%\vspace{0.5em}
\begin{subfigure}{0.32\textwidth}
    \centering
    \includegraphics[width=1.1\linewidth]{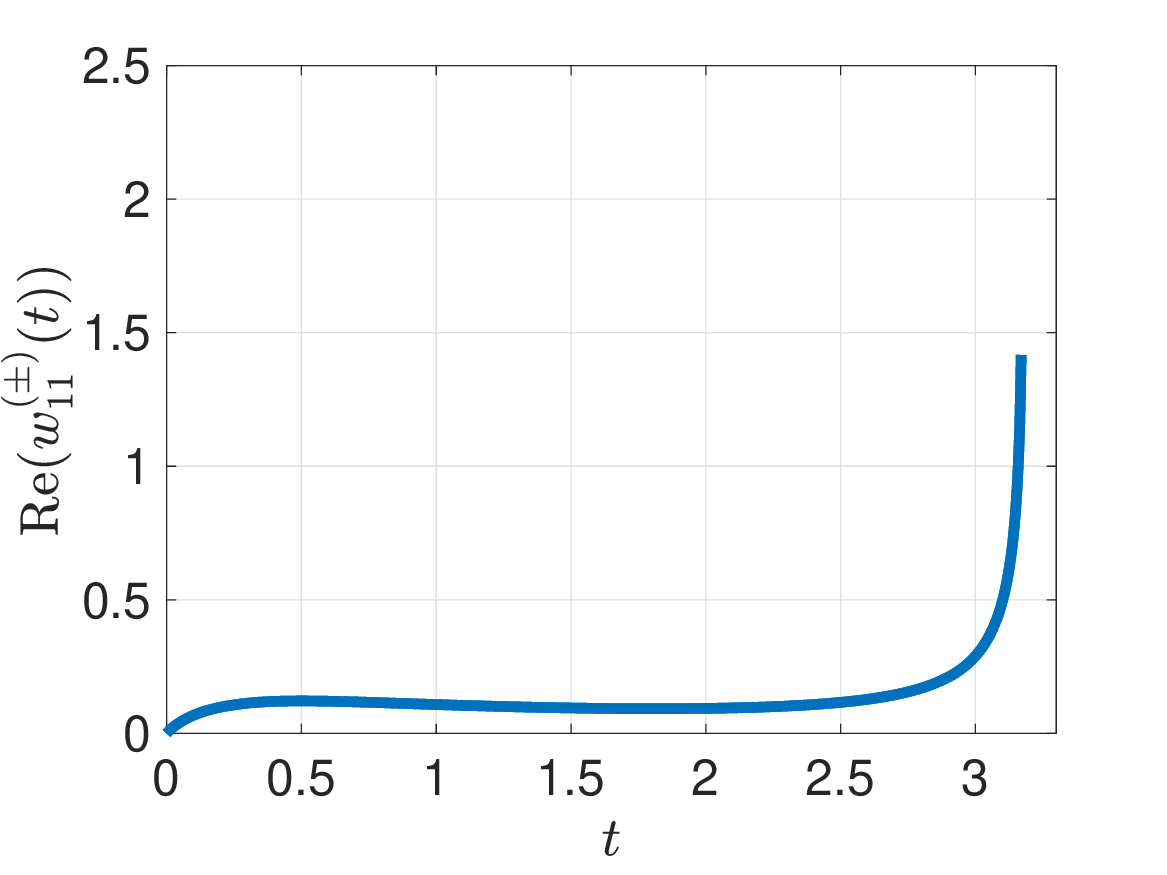}
    \caption{$\alpha=1$}
\end{subfigure}\hfill
\begin{subfigure}{0.32\textwidth}
    \centering
    \includegraphics[width=1.1\linewidth]{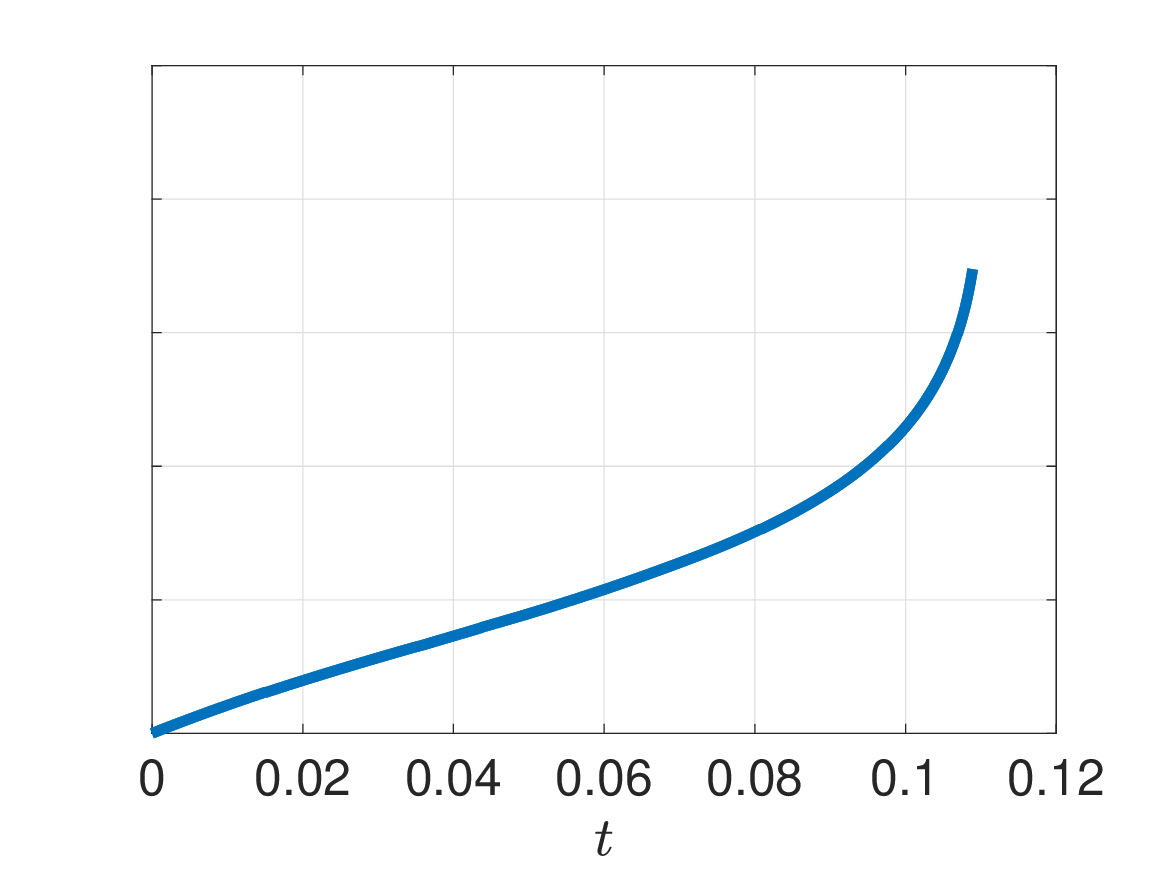}
    \caption{$\alpha=10$}
\end{subfigure}\hfill
\begin{subfigure}{0.32\textwidth}
    \centering
    \includegraphics[width=1.1\linewidth]{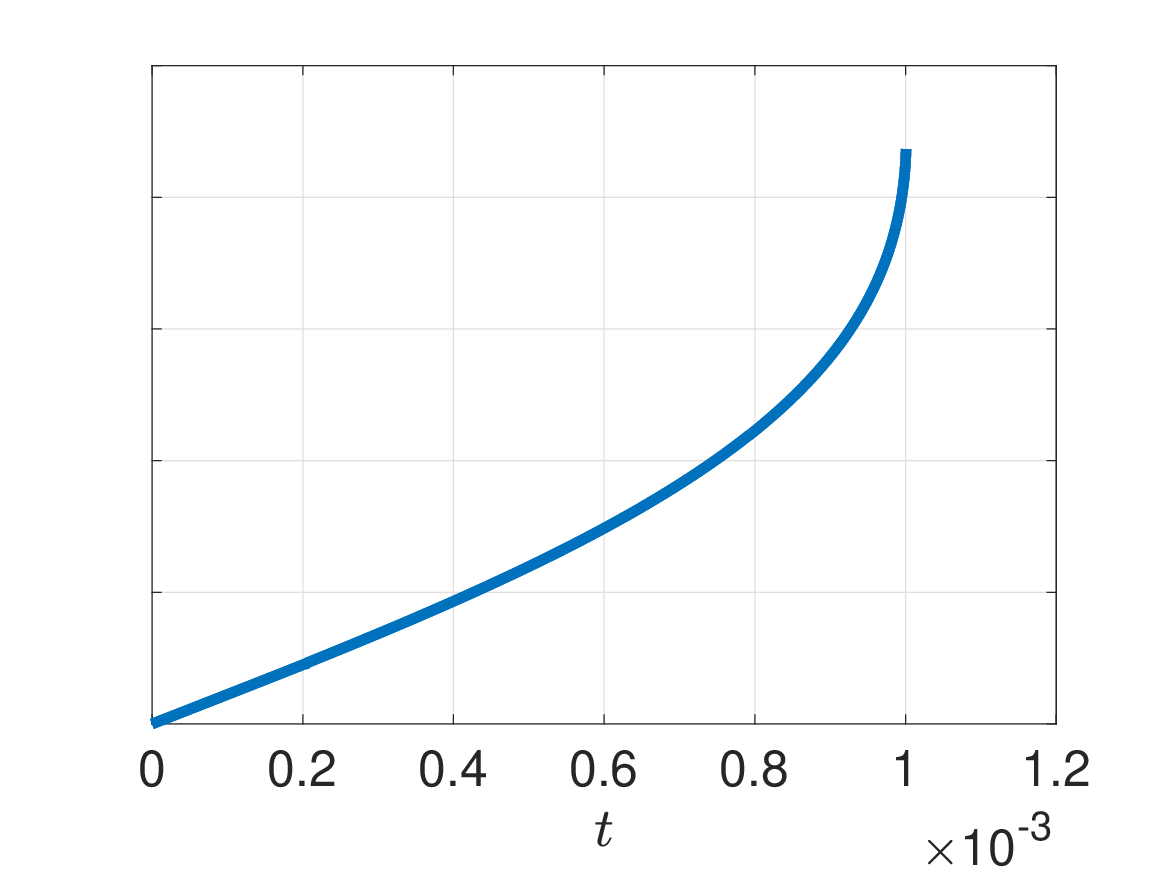}
    \caption{$\alpha=1000$}
\end{subfigure}
 \caption{ Imaginary part of the singularity locations as functions of time  of the extended solution of the NLH (\ref{heat2}) with the initial condition (\ref{inheat})  for  $\alpha=1,10,1000$  computed by equations (\ref{spl})-(\ref{spl1})  and their asymptotics (\ref{asheat})-(\ref{asheat3})   (first row). The blow up occurs for all tree cases at the points $(0,t_b^{(\alpha)})$, where $t_b^{(1)} \approx 3.17395$, $t_b^{(10)}\approx 0.1099821$,   and $t_b^{(1000)} \approx 1.001081 \times 10^{-3}$. Real part  of the corresponding weight locations $w_{11}^{(\pm)}(t)$ (second row) as functions of time was computed with fixed parameters $C_{10}^{(\pm)}(t)=-1$ for each $t \in [0,t_b^{(\alpha)})$.  
  }
  \label{figheatsinwei}
\end{figure}

As in Example \ref{ex31}, we analyze the formation of blow up in the solution $u(x,t)$ by examining the behavior of the singularities $s_1^{(\pm)}(t)$  in the complex plane, in particular how they approach the real axis.
 We start by discussing the case $\alpha=1$.  The singularities $s_1^{(\pm)}(t)$ move along the imaginary axis toward the real axis (Figure \ref{heatsolfig}, second row (a)), until they stop and reverse direction at time $t = 0.503$ (Figure \ref{heatsolfig}, second row (b)). They then move away from the real axis until reaching, at time $t = 2.05$, the position shown in Figure \ref{heatsolfig} (second row (c)), where they turn around once more and rapidly approach the real axis, eventually meeting at the point $(0, t_b^{(\alpha)})$ (Figure \ref{heatsolfig}, second row (d)). The solution profile on the real axis, $u(x,t)$, steepens as the singularities approach the real axis (see Figure \ref{heatsolfig}, first row). The blow up occurs when the two singularities closest to the real axis (in the upper and lower half-planes) collide on it.
As investigated in \cite{FW24}, this type of behavior — in which the singularities, after entering from infinity, reverse direction and move away from the real axis before turning around once again and rapidly approaching it as they meet at the blow up time (see Figure \ref{heatsolfig}, first row (a)) — is observed for small $\alpha$ (roughly $\alpha < 1.2$). 
 For large $\alpha$, in contrast, the singularities approach the real axis monotonically, though not at constant speed (see Figure \ref{heatsolfig}, first row (b) and (c)).
In Figure \ref{heatsolfig} (d), i.e., for time close to the blow up time, the spurious poles and zeros of the extended solution $\Phi(z,t)$  can be seen to align along the branch cuts.
%In Figure \ref{heatsolfig}~(d), i.e., for time close to the blow up time, one can observe that spurious poles and zeros of the extended solution $\Phi(z,t)$ are located along the branch cuts.
It is worth noting that \cite{FW24} employs the so-called quadratic Padé approximation, which is known to provide a more favorable approximation of functions with algebraic branch points.
%Note that in \cite{FW24} the authors use the so-called quadratic Padé approximation, which offers a more favorable approximation for functions with algebraic branch points. 
Note also that \cite{FW24} uses a different color wheel, namely the NIST color wheel.
For $t < t_b^{(\alpha)}$, our method reconstructs the singularities closest to the real axis (in the upper and lower half-planes) as second-order poles. This is illustrated in Figure \ref{heatsolfig} (second row) by the colors wrapping twice around the color wheel when encircling the singularities, as in Example \ref{ex31}. As $t \to 0$ or $t_b^{(\alpha)} \to \infty$, the singularities closest to the real axis are instead reconstructed as simple poles. This numerical result agrees completely with the asymptotic analysis in \cite{FW24}.

Since no explicit solution is known for the NLH \eqref{heat2} with the initial condition \eqref{inheat}, we use the asymptotic estimates from \cite{FW24} to validate the numerical performance of our approximation through equations \eqref{spl}--\eqref{spl1}. The corresponding results are shown in Figure \ref{figheatsinwei} (first row).
  Let $z^{(\pm)}(t)$ be the exact locations of the two  closest to the real axis  singularities and let  $z^{(\pm)}(t)=\pm \mathrm{i} \, \sigma(t)$. Then for all $\alpha>0$, we have \cite{FW24}
\begin{equation}\label{asheat}
    \sigma(t) \sim \log(2/(\alpha t))+2 t \log(1/( t))-1.05695 t,  \quad t \rightarrow 0.
\end{equation}
In the case of large $\alpha$, i.e for $\alpha\rightarrow \infty$ with $t=\mathcal{O}(1/\alpha)$, the following asymptotic holds
%\begin{equation}\label{asheat1}
\begin{align}
    \sigma(t) \sim \mathrm{acosh}&(1/(\alpha t)) + t \sqrt{1-\alpha^2 t^2}  \notag \\
   & \times \left( 2\log(\alpha) - \frac{1-2\alpha^2 t^2}{1-\alpha^2 t^2} - 2\log(\alpha t) - 2  \log(1- \alpha^2 t^2)- 0.05695\right), \label{asheat1}
\end{align}
 %   \sigma(t) \sim \mathrm{acosh}(1/(\alpha t))+ t \sqrt{1-\alpha^2 t^2} \left( 2\log(\alpha) - \frac{1-2\alpha^2 t^2}{1-\alpha^2 t^2} - 2\log(\alpha t) - 2  \log(1- \alpha^2 t^2)- 0.05695\right)
%\end{equation}
with the leading order approximation 
\begin{equation}\label{asheat2}
    \sigma(t) \sim \mathrm{acosh}(1/(\alpha t)), \quad \alpha\rightarrow \infty, \ t=\mathcal{O}(1/\alpha).
\end{equation}
Additionally, for the case $\alpha=1$ it was shown that \cite{W03}
\begin{equation}\label{asheat3}
    \sigma(t) \sim \sqrt{8(t_b^{(1)}-t)|\log(t_b^{(1)}-t)|}, \quad t\rightarrow t_b^{(1)}.
\end{equation}
%In (\ref{asheat3}), we have to assume that blow up time is known and we use the numerical estimates. 

\begin{figure}[h!]
  \centering
  \begin{subfigure}[b]{0.45\linewidth}
    \includegraphics[width=1\linewidth]{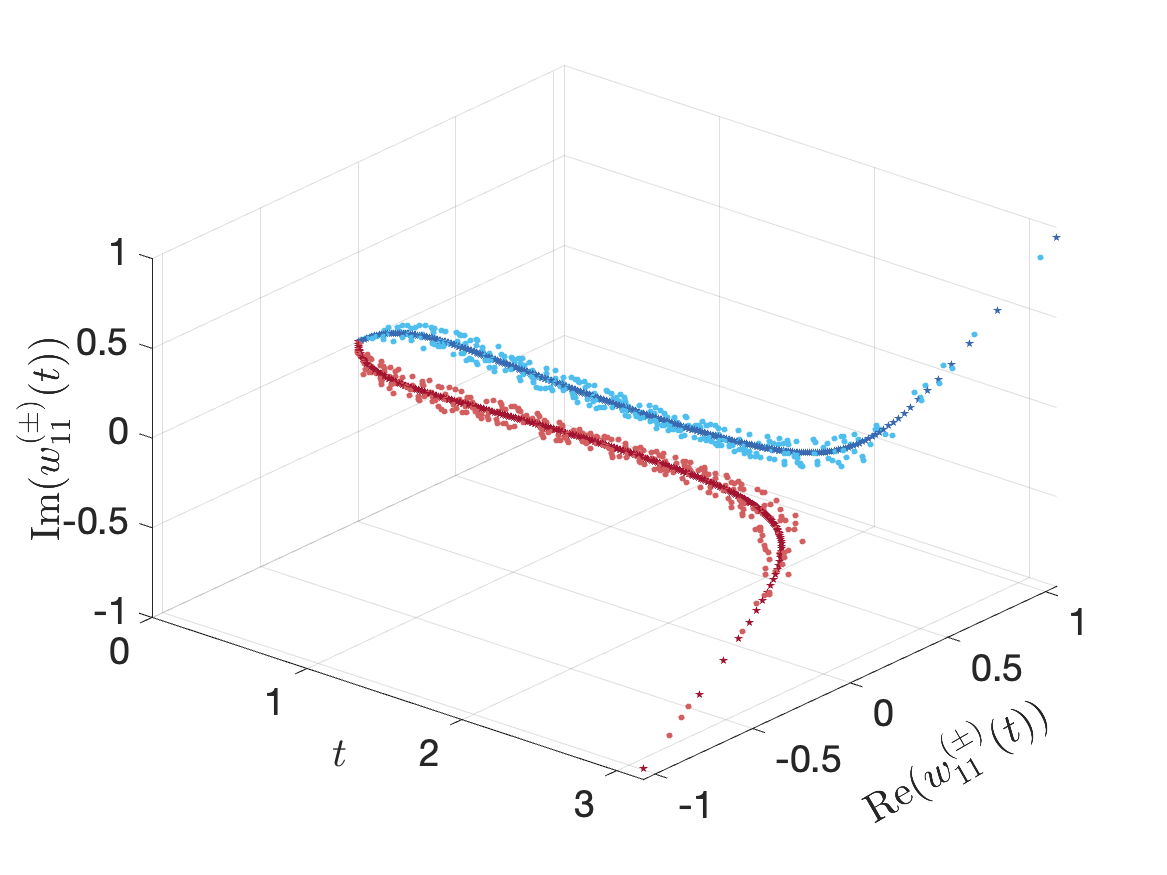}
   % \caption{$t=0.15$}
  \end{subfigure}
  \begin{subfigure}[b]{0.45\linewidth}
    \includegraphics[width=1\linewidth]{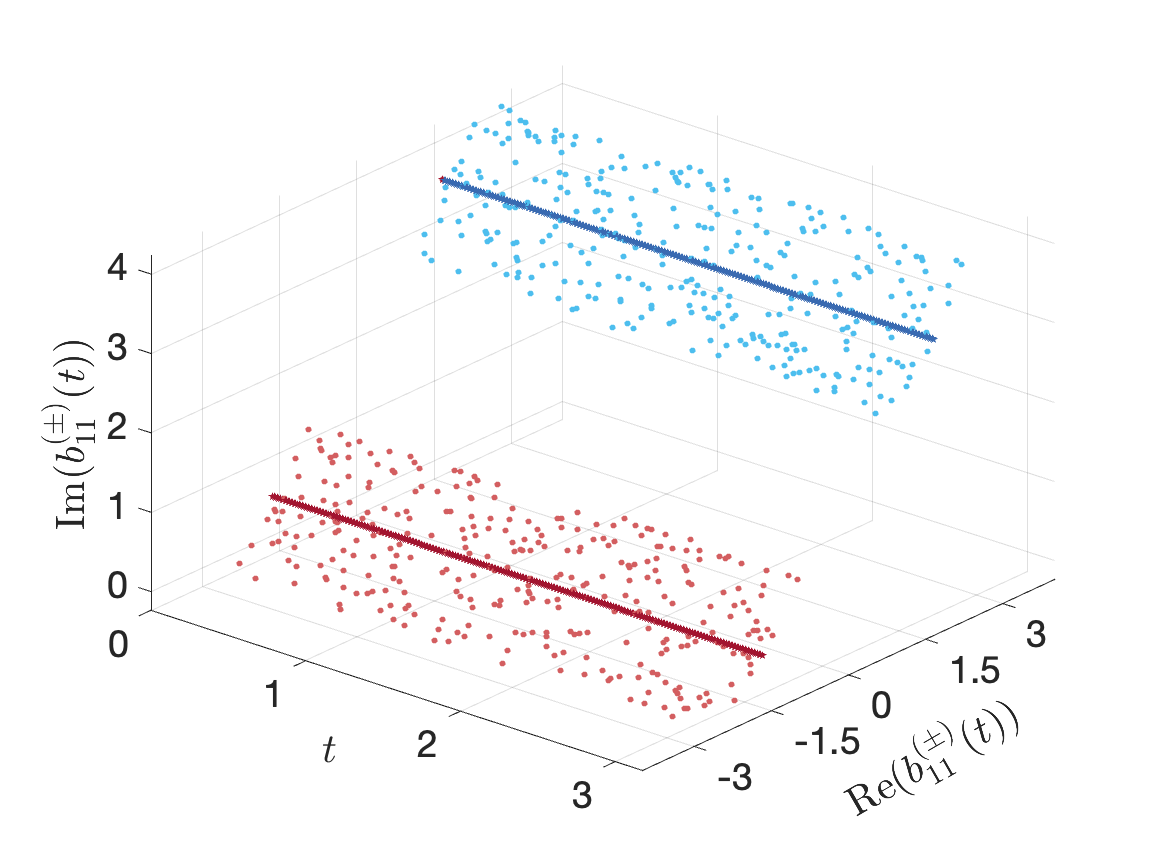}
   % \caption{$t=0.503$}
  \end{subfigure}
  \caption{   Weights $w_{11}^{(\pm)}(t)$ (left) and biases  $b_{11}^{(\pm)}(t)$ (right) of the hidden layers of the components $\Phi^{(\pm)}$  used in (\ref{spl})-(\ref{spl1}) for computation of singularities  $s_{1}^{(\pm)}(t)$ of the extended solution of the NLH (\ref{heat2}) subject to the initial condition (\ref{inheat}) with $\alpha=1$ for time interval $[0,t_b^{(1)})$.
Dark blue: $w_{11}^{(-)}(t)$ and $b_{11}^{(-)}(t)$ computed with fixed parameters $C_0^{(-)}(t)=-1$,  light blue: $w_{11}^{(-)}(t)$ and $b_{11}^{(-)}(t)$ computed with random choice of parameters  $C_{10}^{(-)}(t) \in [-1.5,-0.5]$.  Dark red: $w_{11}^{(+)}(t)$ and $b_{11}^{(+)}(t)$ computed with fixed  $C_0^{(+)}(t)=1$,  light red: $w_{11}^{(+)}(t)$ and $b_{11}^{(+)}(t)$ computed using random choice of parameters  $C_{10}^{(+)}(t) \in [0.5,1.5]$. In all experiments, we set the pole of the activation functions to $z_0=-2\mathrm{i}$. }
  \label{figwbhh}
\end{figure}

Finally, we study trajectories of weights $w_{11}^{(\pm)}(t)$ and biases $b_{11}^{(\pm)}(t)$ as functions of time in the complex plane. According to equations \eqref{spl}--\eqref{spl1}, their behavior determines the behavior of the singularities $s_1^{(\pm)}(t)$.
Let us consider the case $\alpha = 1$, though for other $\alpha > 0$ the results will be similar.
First, we  fix  parameters $C_{10}^{(\pm)}(t)=-1$ for all $t \in [0,t_b^{(1)})$.
According to our numerical computations, weights  are real, positive and satisfy the condition $w_{11}^{(-)}(t)=w_{11}^{(+)}(t)$, and  biases are constant with the values $b_{11}^{(\pm)}(t)= 4.236067977499788$. If a function $w_{11}^{(-)}(t)$ increases (decreases), taking into account that $w_{11}^{(-)}(t)>0$ and $b_{11}^{(-)}(t)+z_0>0$ is constant,  we obtain by \eqref{spl1} that  $ \mathrm{e}^{\mathrm{i} s_{1}^{(-)}(t)}=\mathrm{e}^{- \mathrm{Im} \, s_{1}^{(-)}(t)}$ also increases (decreases). Since $\mathrm{Im} \, s_{1}^{(-)}(t)>0$, we conclude that $\mathrm{Im} \, s_{1}^{(-)}(t)$ must decreases (increases).  Analogously, applying \eqref{spl}, we can derive that if a function $w_{11}^{(+)}(t)$ increases (decreases) and additionally the following  properties hold, $w_{11}^{(+)}(t)>0$, $b_{11}^{(+)}(t)+z_0>0$ is constant, and $\mathrm{Im} \, s_{1}^{(+)}(t)<0$, then  $\mathrm{Im} \, s_{1}^{(+)}(t)$ increases (decreases). This analysis completely coincides with numerical results presented in Figure \ref{figheatsinwei}.   Note that if we choose different values of $C_{10}^{(-)}(t)$ and $C_{10}^{(+)}(t)$ or even choose them randomly, and a complex pole of the activation function $z_0$, then weights and biases of the hidden layers will be different ($w_{11}^{(-)}(t)\neq w_{11}^{(+)}(t)$ and $b_{11}^{(-)}(t) \neq b_{11}^{(+)}(t)$) and complex, but the general behavior will be unchanged (see Figure \ref{figwbhh}).  

\end{example}

   In the proposed framework, the formation of blow up is characterized by the evolution of complex singularities detected through numerical analytic continuation. The neural network $\Phi$ in (\ref{nnapr}) extends the solution $u(x,t)$ into the complex plane, enabling the tracking of singularity trajectories. The proximity of the nearest singularities to the real axis provides an accurate indicator of the onset of blow up. It occurs when one or more singularities collide with the real axis at a finite time. In Example \ref{ex31}, the initial condition (\ref{heatsol}) already possesses complex singularities, whereas in Example \ref{ex32}, the initial condition (\ref{inheat}) becomes an entire function upon complexification, and it is the nonlinear term in equation (\ref{heat2}) that develops singularities in finite time.

%In Example  \ref{ex31}, the initial condition (\ref{heatsol}) has already complex singularities, while in Example   \ref{ex32}, the initial condition (\ref{inheat}) is an entire function after complexification and the nonlinear part of the equation (\ref{heat2}) develops singularities in finite time. 

%\newpage

\subsection{Nonlinear Burgers equation}
\label{secburgers}

 We now study  the formation of a shock. To this end we analyze the nonlinear Burgers equation.

\begin{example}
Let the nonlinear Burgers equation be given by 
\begin{equation}\label{burdis}
    u_t+u u_x= \nu \, \mathrm{e}^{\mathrm{i}  \theta}  u_{xx}, \quad \nu, \theta \geq 0,
\end{equation}
with the initial condition 
\begin{equation}\label{initbur}
    u(x,0)=-\sin (x),
\end{equation}
for  $ x \in [-\pi,\pi)$ and $t\geq 0$.  We investigate the extended solution $u(z,t)$ for the domain $z \in D=[-\pi, \pi)\times  \mathrm{i} \, (-1.8, 1.8)$ and  $t\geq 0$. First, employing $2n$  function values  $u\left( \frac{ \pi j}{ n} , 0 \right)$, $j=-n,\dots, n-1$, of the initial condition (\ref{initbur}) with $n=200$, we construct approximation of the solution $u(x,t)$ as in (\ref{appsr}) by the Fourier spectral method.
  The system of ODEs  (\ref{odes}) has the form    
%\begin{equation}\label{odesbergers}
$$
 \frac{\mathrm{d} \hat{u}_k(t)}{\mathrm{d} t } + \frac{1}{2} \mathrm{i} \, k \sum\limits_{\ell=-n}^n \hat{u}_\ell(t) \, \hat{u}_{k-\ell}(t)+ \nu \, \mathrm{e}^{\mathrm{i} \theta} k^2 \hat{u}_k(t)=0  , \quad k=-n,\dots,n.
 $$
% \end{equation}
Next, we compute the numerical analytic continuation $\Phi(z,t)$ of $u(x,t)$ by (\ref{nnapr}) and detect their complex singularities by (\ref{spl})-(\ref{spl1}). 
 
\begin{figure}[h!]
    \centering
  \begin{subfigure}[b]{0.24\linewidth}
    \includegraphics[width=1.1\linewidth]{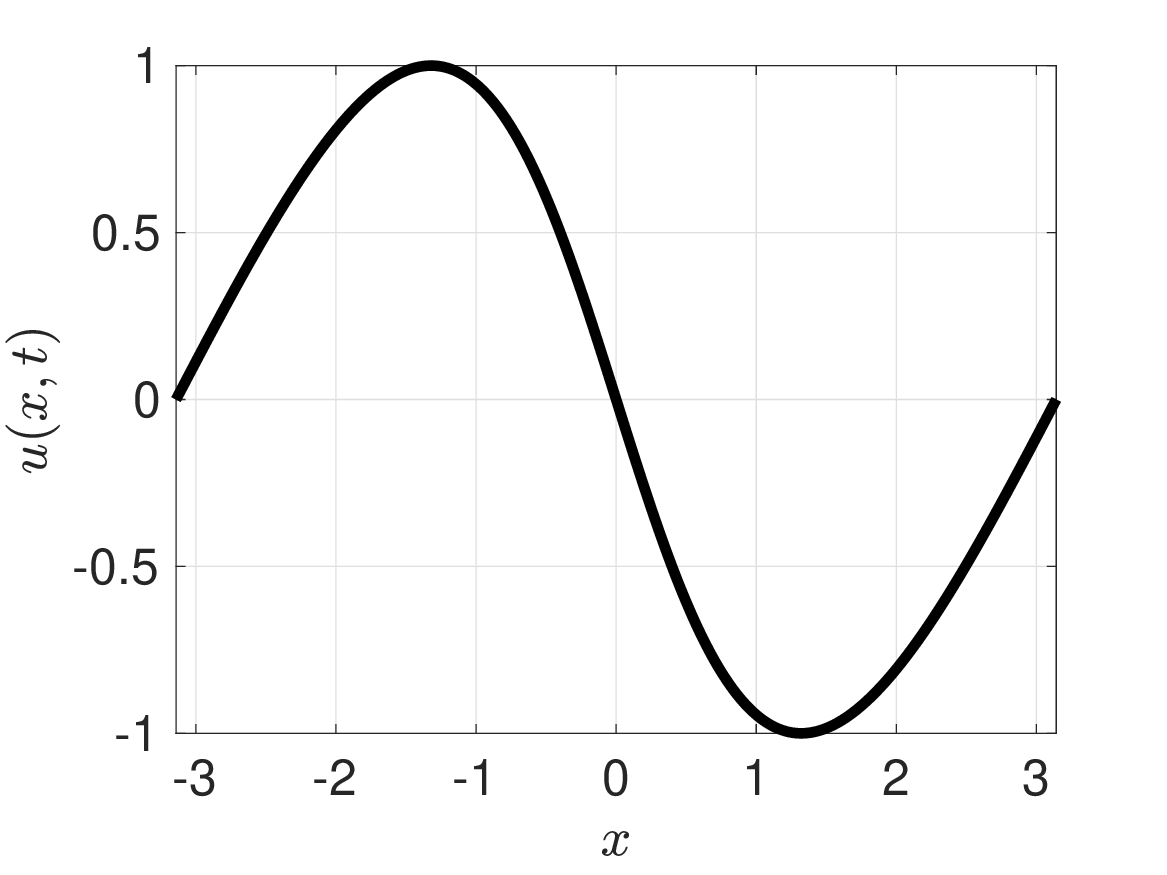}
  %  \caption{$t=0.25$}
  \end{subfigure}
  \begin{subfigure}[b]{0.24\linewidth}
    \includegraphics[width=1.1\linewidth]{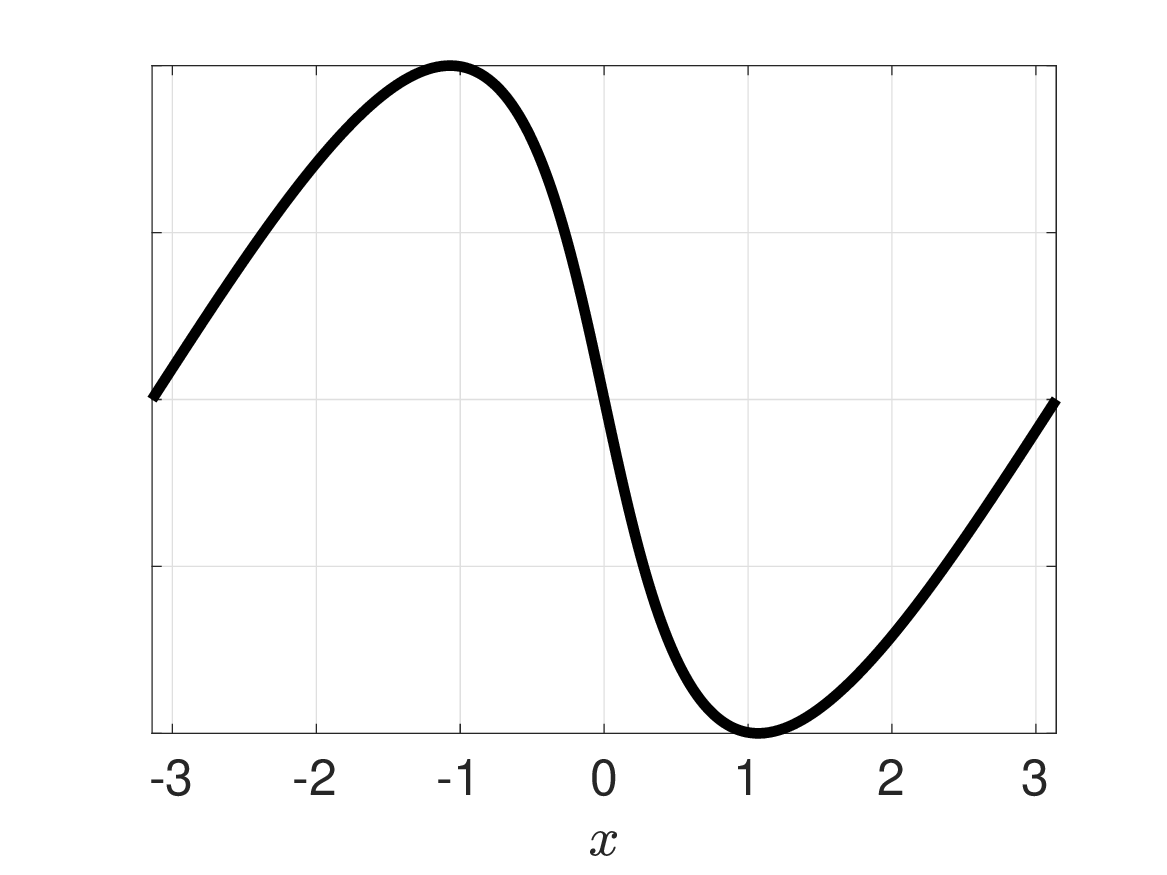}
   % \caption{$t=0.5$}
  \end{subfigure}
    \begin{subfigure}[b]{0.24\linewidth}
    \includegraphics[width=1.1\linewidth]{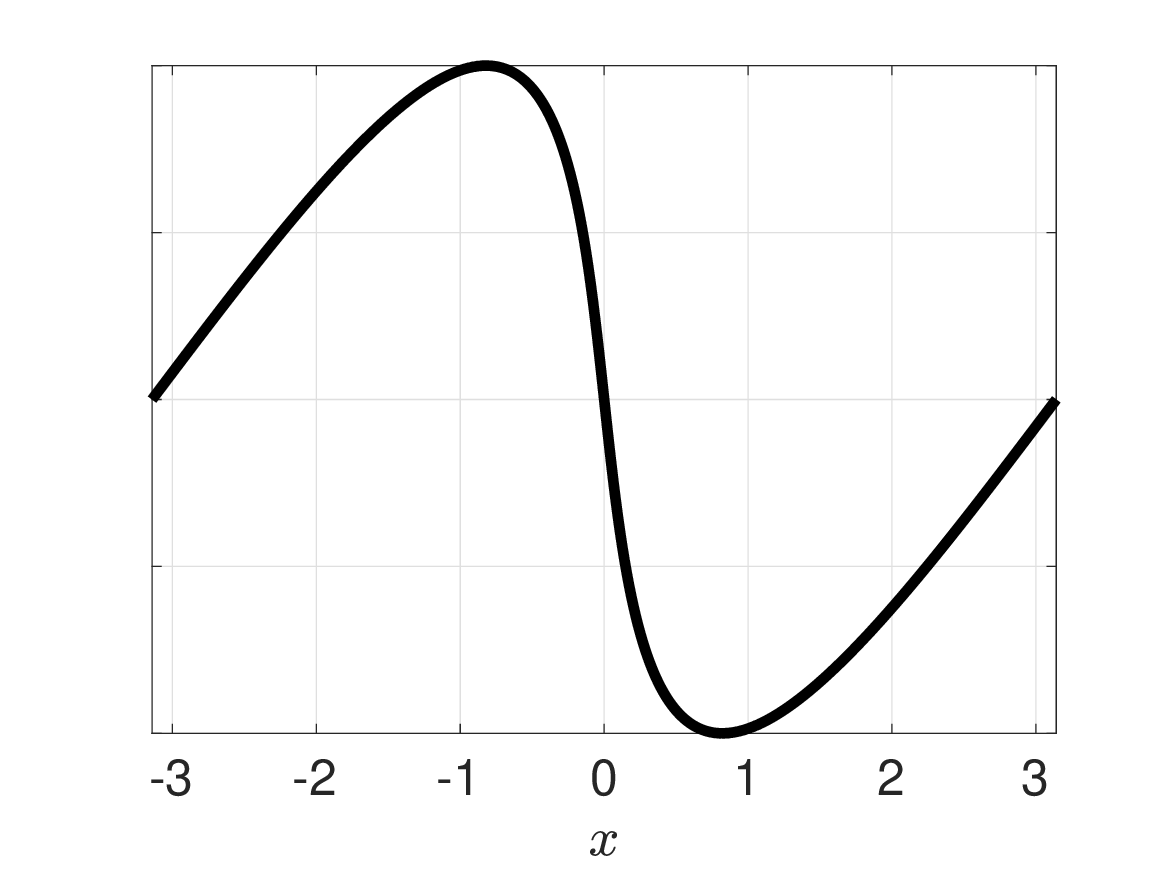}
  %  \caption{$t=0.5$}
  \end{subfigure}
  \begin{subfigure}[b]{0.24\linewidth}
    \includegraphics[width=1.1\linewidth]{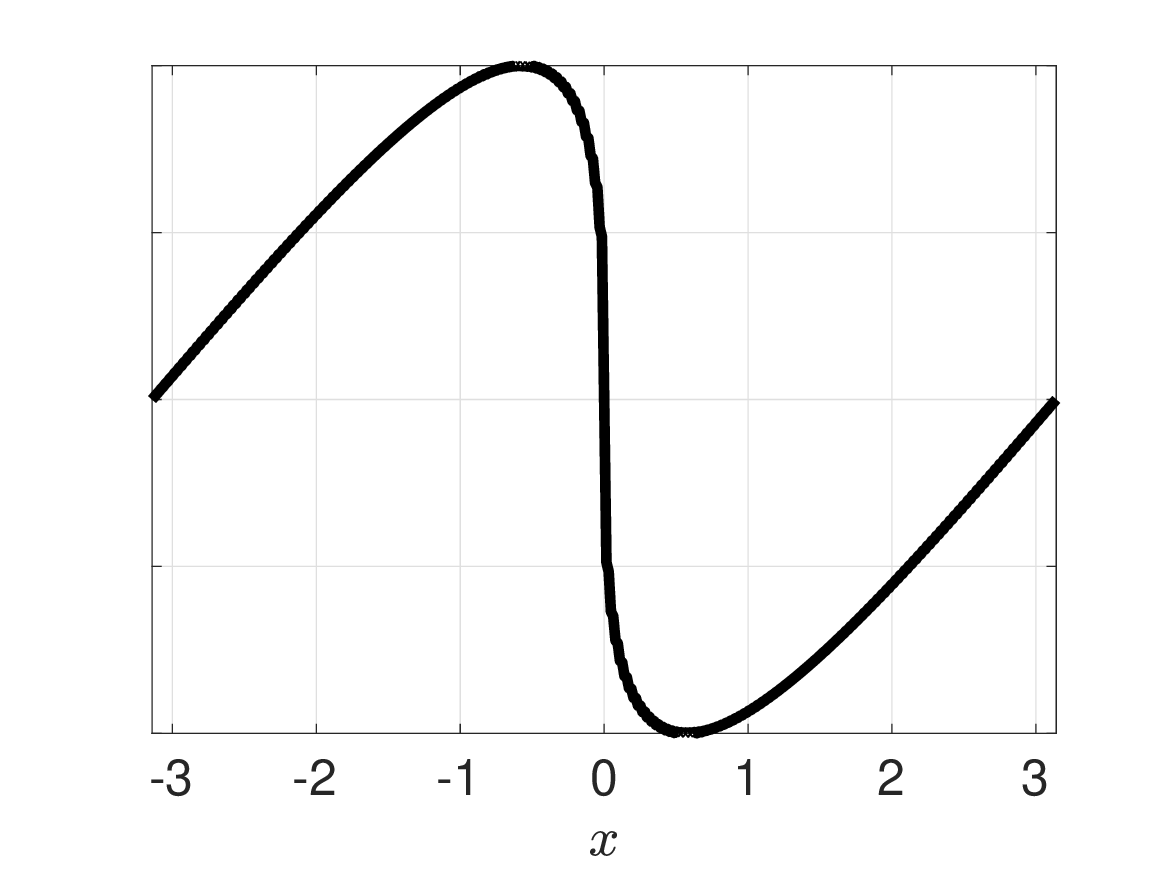}
  %  \caption{$t_s=1$}
  \end{subfigure}
  \begin{subfigure}[b]{0.24\linewidth}
    \includegraphics[width=1.1\linewidth]{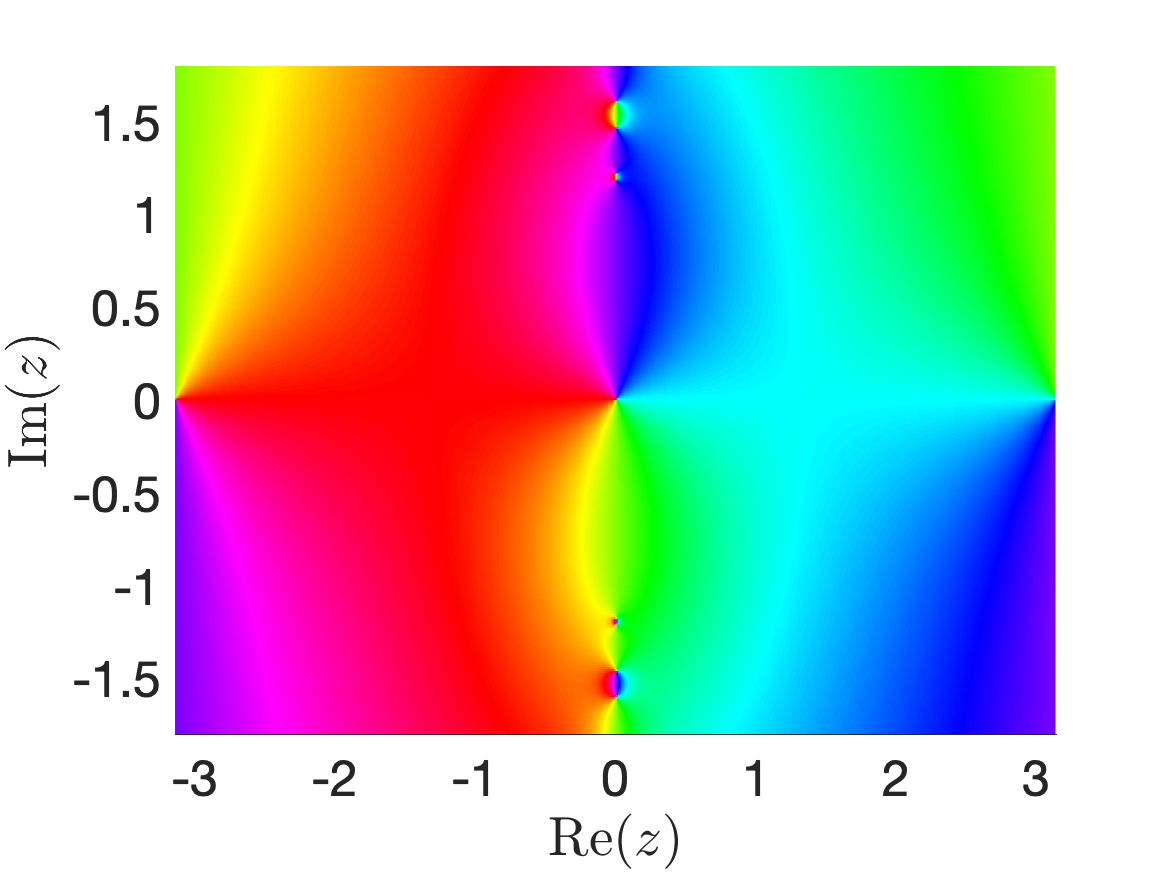}
    \caption{$t=0.25$}
  \end{subfigure}
  \begin{subfigure}[b]{0.24\linewidth}
    \includegraphics[width=1.1\linewidth]{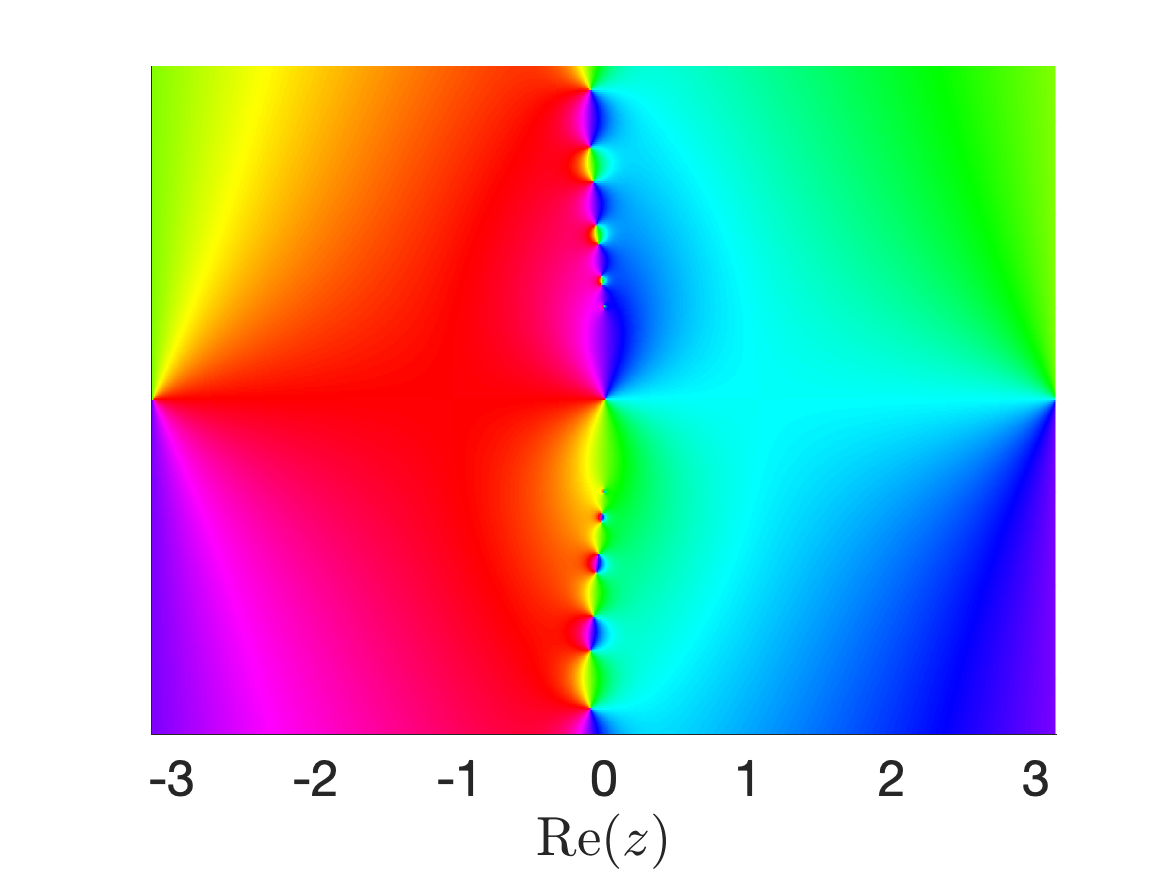}
    \caption{$t=0.5$}
  \end{subfigure}
    \begin{subfigure}[b]{0.24\linewidth}
    \includegraphics[width=1.1\linewidth]{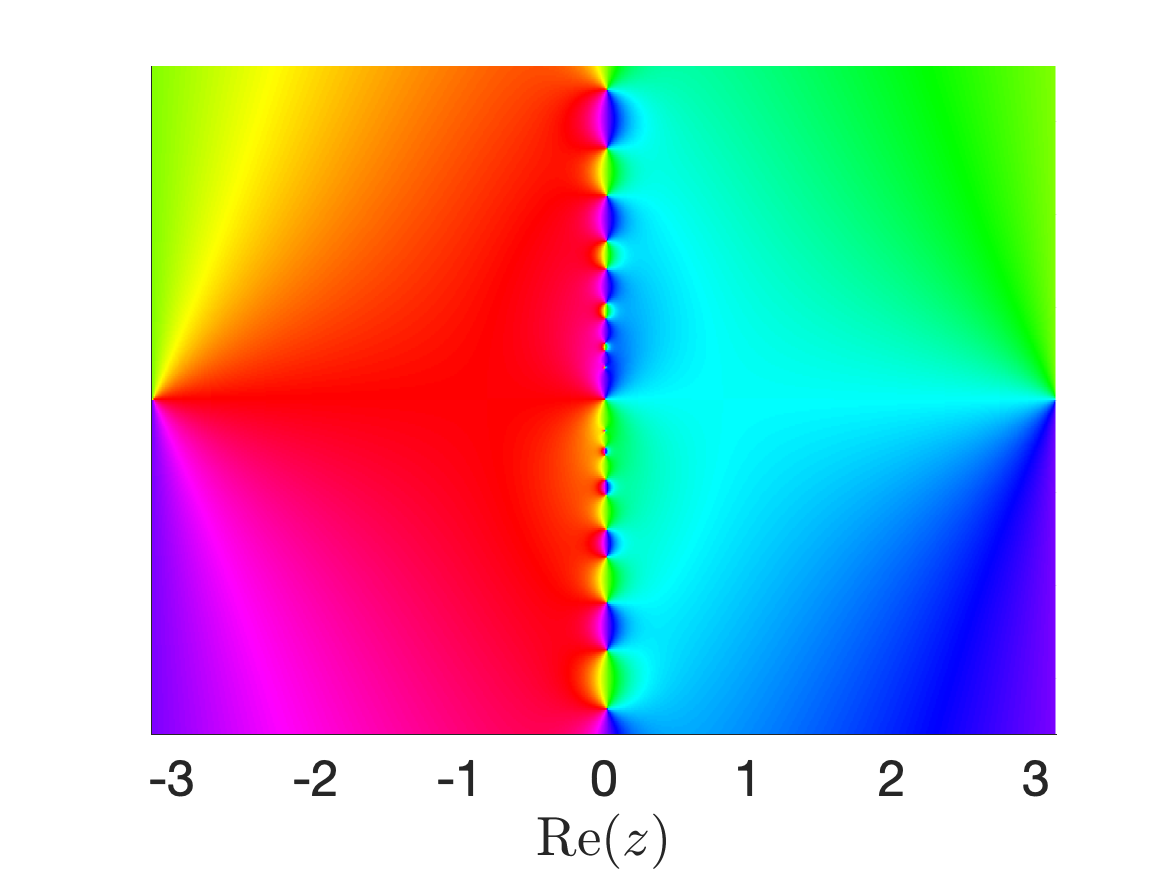}
    \caption{$t=0.75$}
  \end{subfigure}
  \begin{subfigure}[b]{0.24\linewidth}
    \includegraphics[width=1.1\linewidth]{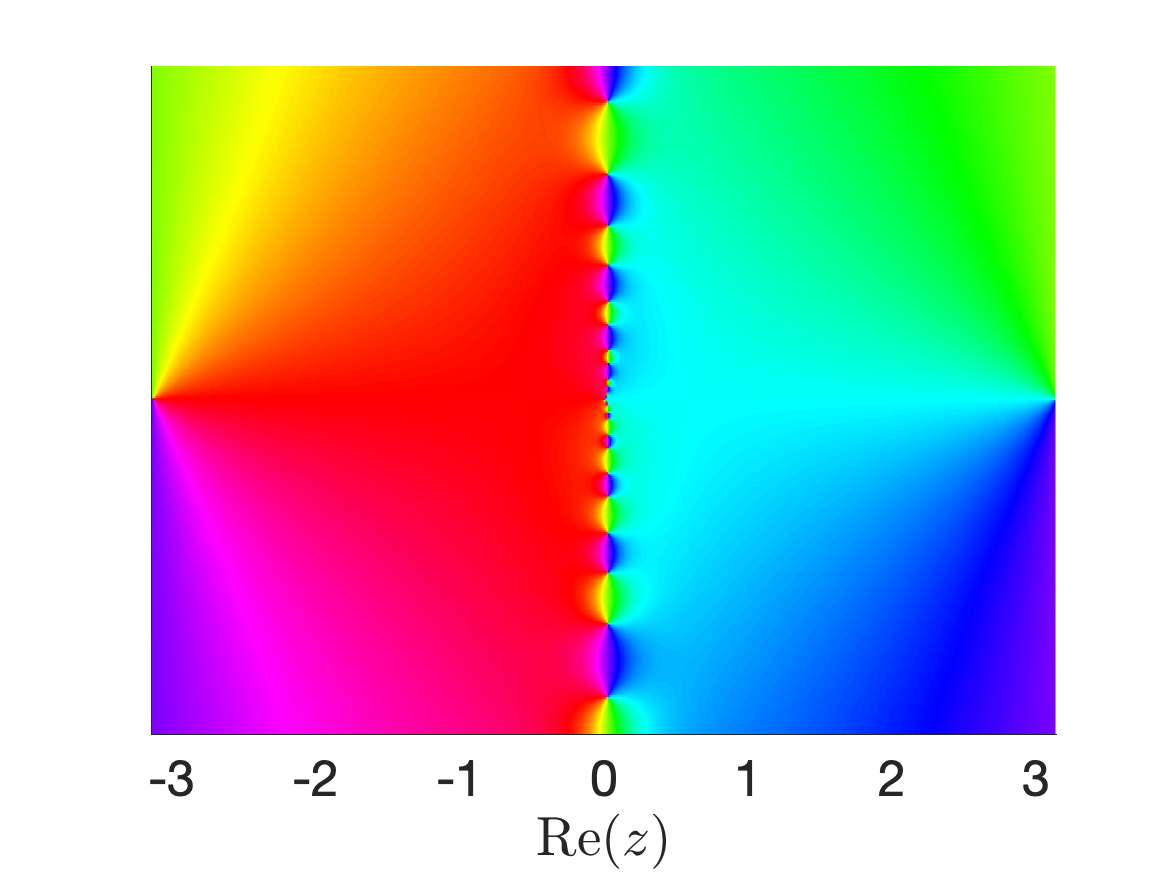}
    \caption{$t_s=1$}
  \end{subfigure}
  \caption{ The solution $u(x,t)$ (first row) and  the phase portraits of its extension   $\Phi(z,t)$ as in (\ref{nnapr})  into the complex plane (second row) of the Burgers equation (\ref{burdis})  with the initial condition (\ref{initbur}) for $\nu=0$ (inviscid case) for time $t=0.25, 0.5, 0.75$ and $t_s=1$.  Spurious poles and zeros are located along the  branch cuts. The shock occurs at time $t_s=1$. }
  \label{inburext}
\end{figure}

The case $\nu=0$ in (\ref{burdis})  corresponds to the \textbf{inviscid Burgers equation} that subject to the initial condition (\ref{initbur})    develops a shock at time $t_s=1$. 
 In Figure \ref{inburext} (first row) we observe the typical steepening of the solution $u(x,t)$ for $t=0.25, 0.5, 0.75$ until the shock forms at $t_s=1$.  
 The extended solution $u(z,t)$ has two square root singularities, placed symmetrically on the imaginary axis with respect to the origin. These singularities are  ``born'' at $\pm \mathrm{i} \infty$ for small $t>0$ and travel from there on the imaginary axis toward the real axis, where they meet at time $t_s=1$  when the shock becomes visible on the real axis. 
This behavior was first reported in \cite{BF84, BF89} for a cubic polynomial initial condition and later studied for the trigonometric initial condition (\ref{initbur}) in \cite{W22, CG15}.
 %This behaviour was first reported in  \cite{BF84, BF89} with a cubic polynomial as the initial condition. For the trigonometric initial condition (\ref{initbur}), it was discussed in \cite{W22, CG15}.
 In \cite{W22} it was investigated that the branch points dynamics of the extended solution $u(z,t)$ of the inviscid Burgers equation with the initial condition (\ref{initbur}) is given by
\begin{equation}\label{brpoints}
    \zeta^{(\pm)}(t)=\pm \mathrm{i} \, (\sqrt{1-t^2}- \tanh^{-1}\sqrt{1-t^2}), \quad 0<t <1.
\end{equation}
%In  Figure \ref{inburext} (second row), we present the approximation $\Phi(z,t)$  of the extended solution as in (\ref{nnapr})    for time $t=0.25, 0.5,0.75$ and $t_s=1$. As is typical for rational approximation, our neural network-based approach places spurious poles and zeros along the branch cuts. We would like to mention again  quadratic Padé approximation \cite{FH19}, which is better suited to approximating functions with algebraic branch points.
%In Figure \ref{inburext}, we can also observe how the complex singularities of the extended solution approaching the real axis. 
In  Figure \ref{inburext} (second row), we present the approximation $\Phi(z,t)$  of the extended solution as in (\ref{nnapr})    for time $t=0.25, 0.5,0.75$ and $t_s=1$. As is typical for rational approximation, our neural network-based approach places spurious poles and zeros along the branch cuts. We would like to mention again  quadratic Padé approximation \cite{FH19}, which is better suited to approximating functions with algebraic branch points.
%In Figure \ref{inburext}, we can also observe how the complex singularities of the extended solution approaching the real axis. 

\begin{figure}[h!]
  \centering
  \begin{subfigure}[b]{0.45\linewidth}
    \includegraphics[width=0.92\linewidth]{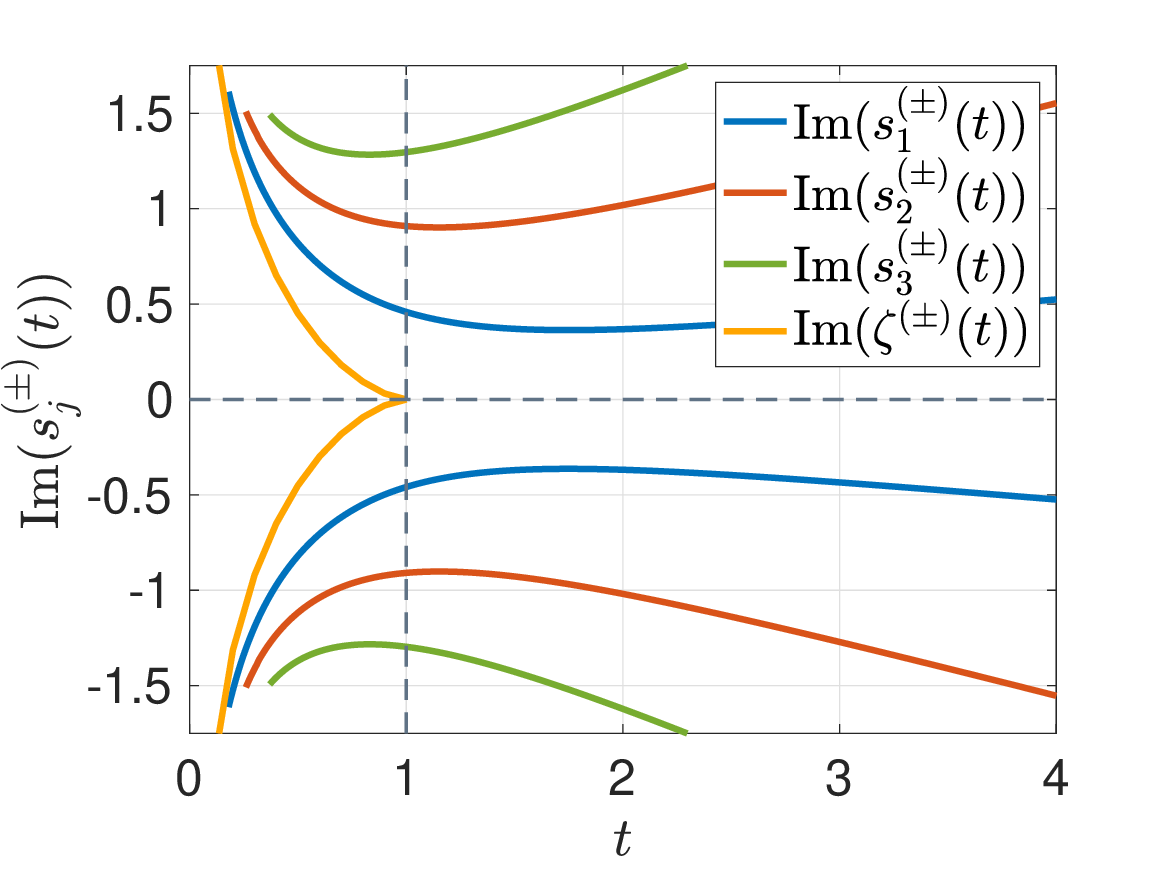}
   % \caption{$t=0$}
  \end{subfigure}
    \begin{subfigure}[b]{0.45\linewidth}
    \includegraphics[width=0.95\linewidth]{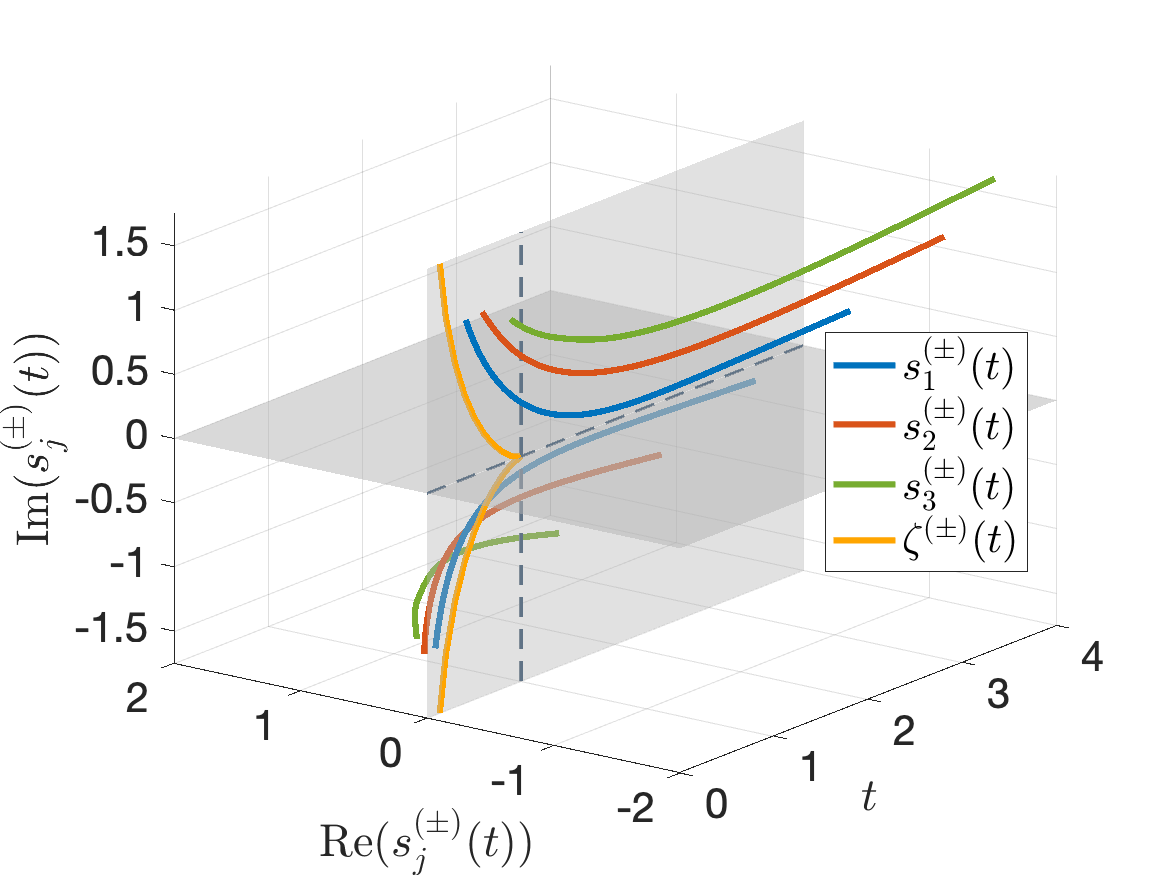}
   % \caption{$t=0$}
  \end{subfigure}
  \begin{subfigure}[b]{0.45\linewidth}
    \includegraphics[width=0.92\linewidth]{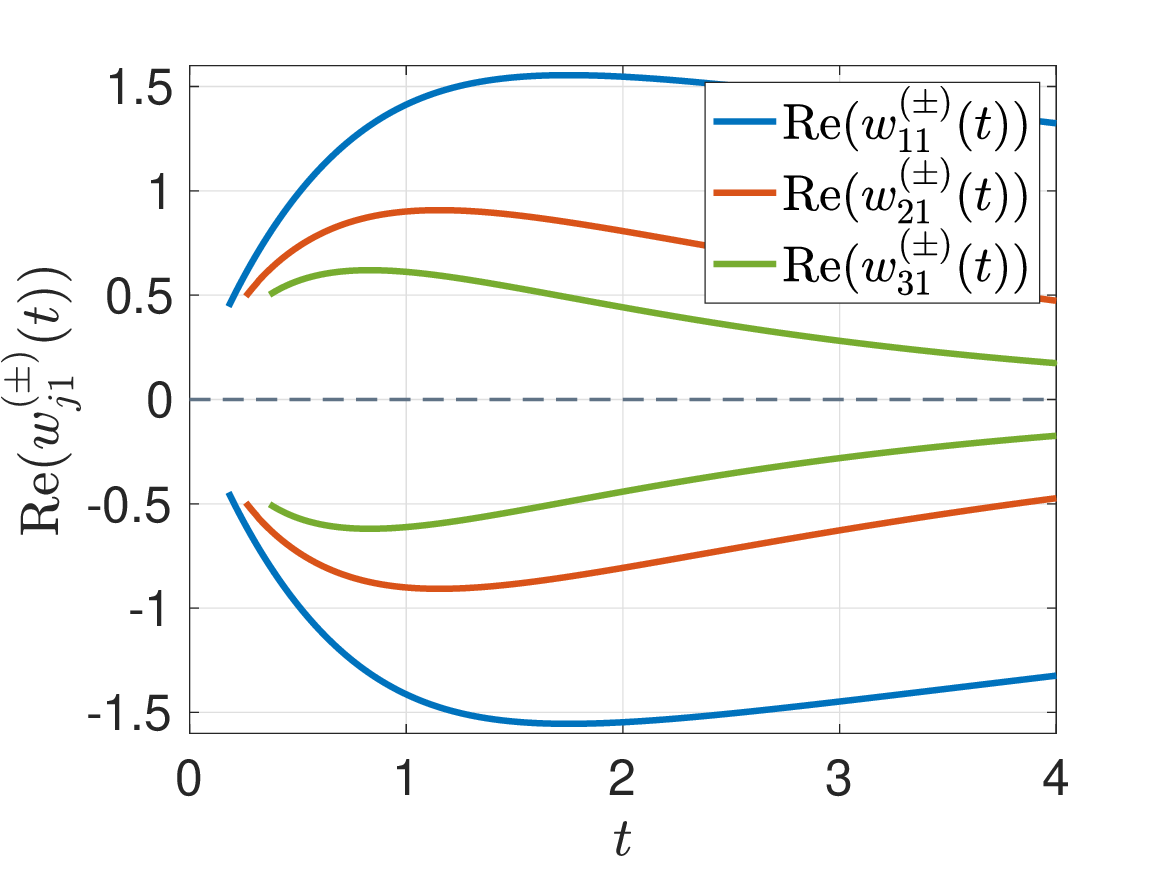}
  %  \caption{$t=0.5$}
  \end{subfigure}
  \begin{subfigure}[b]{0.45\linewidth}
    \includegraphics[width=0.95\linewidth]{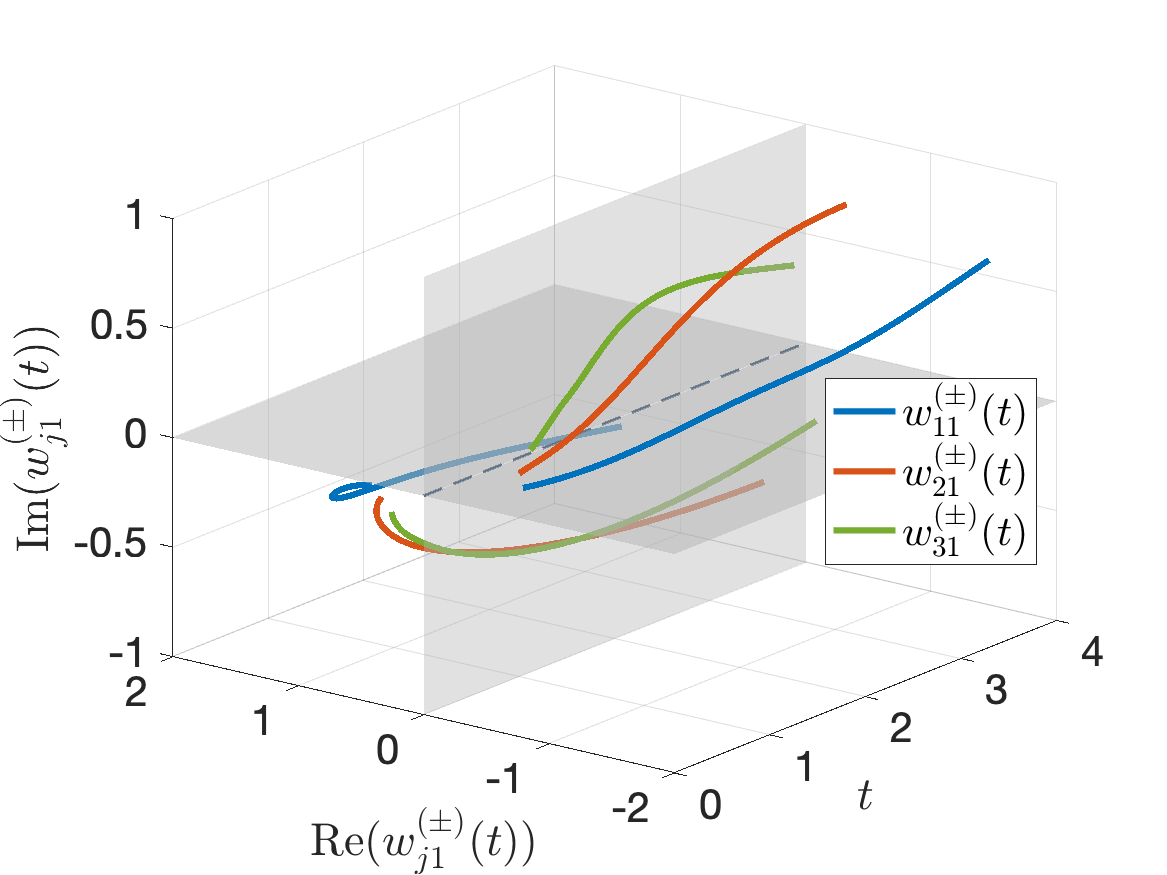}
  %  \caption{$t=0.5$}
  \end{subfigure}
      \begin{subfigure}[b]{0.45\linewidth}
    \includegraphics[width=0.92\linewidth]
    {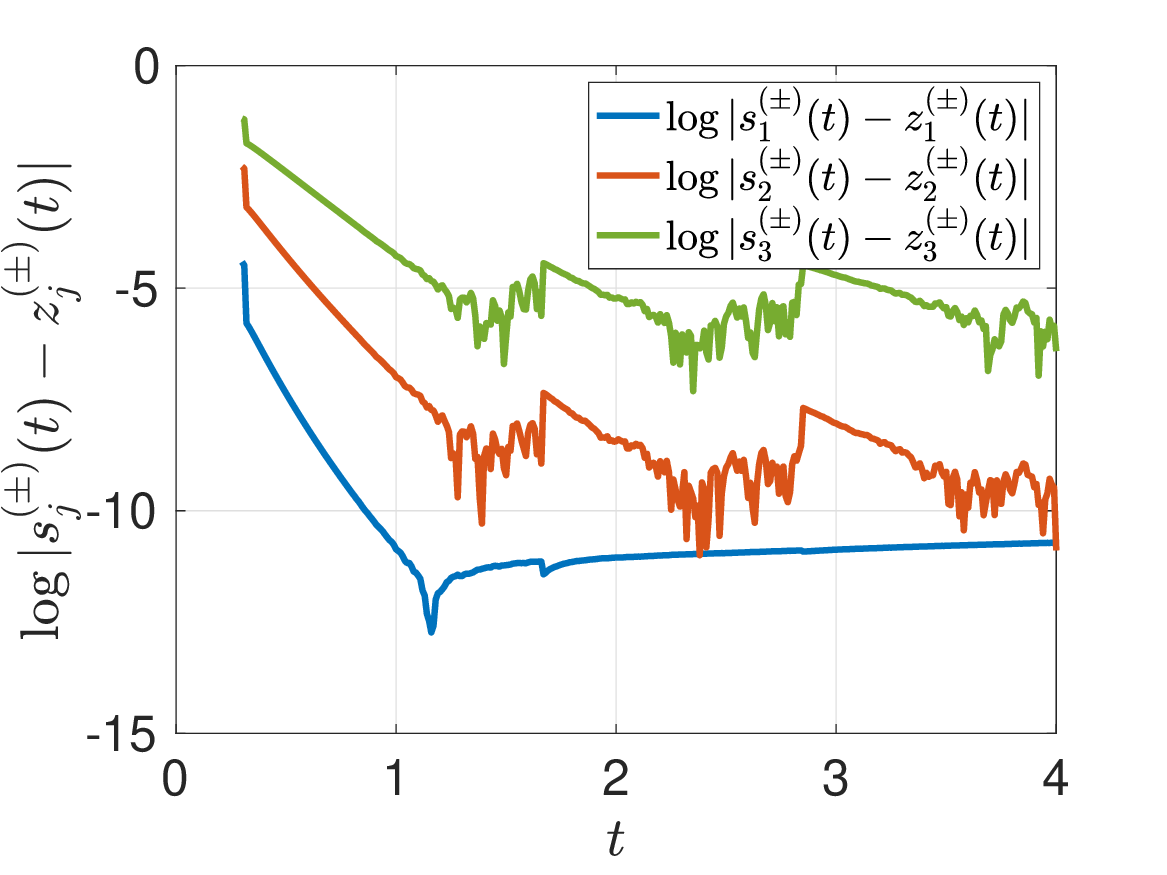}
    \caption{$\theta=0$}
      \end{subfigure}
          \begin{subfigure}[b]{0.45\linewidth}
    \includegraphics[width=0.92\linewidth]
    {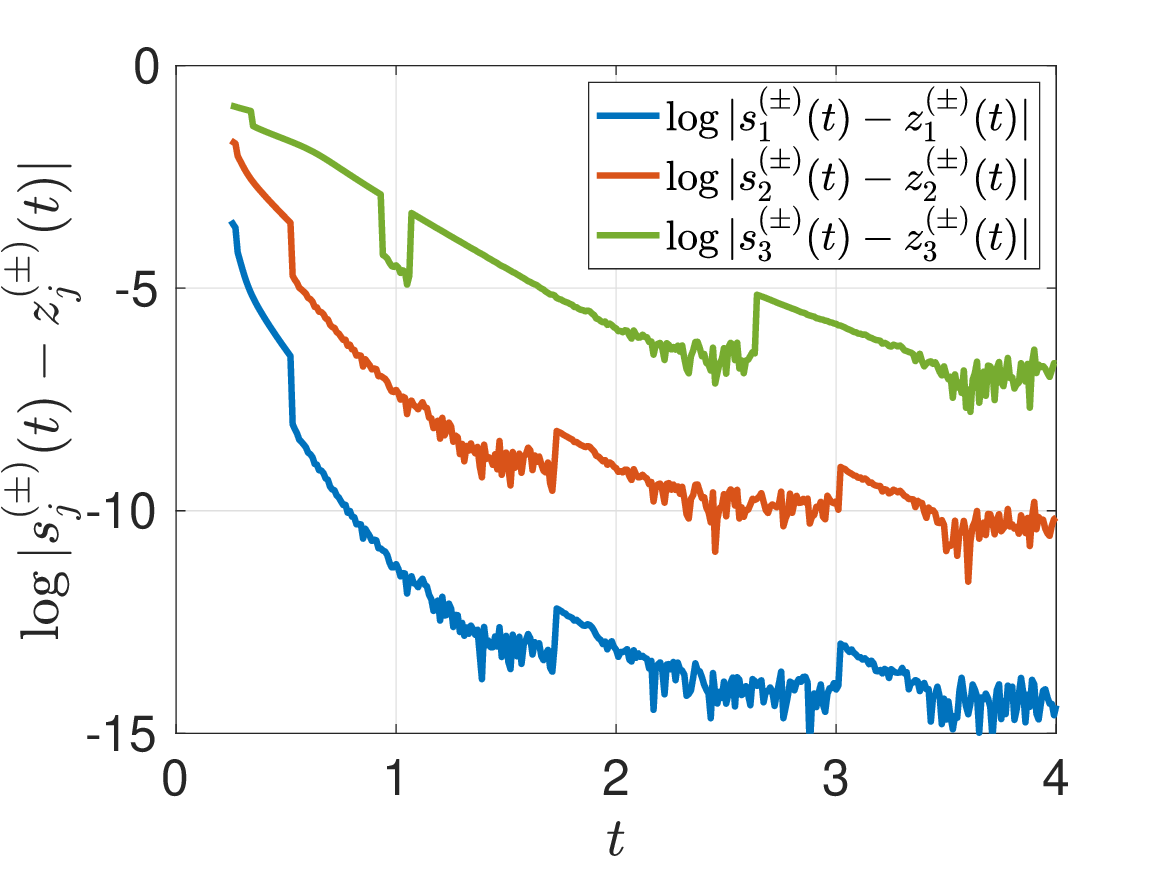}
    \caption{$\theta=\pi/4$}
      \end{subfigure}
  \caption{ Trajectories in the complex plane of the first three complex conjugate pairs of poles $s_j^{(\pm)}(t)$, $j=1,2,3$ as functions of time computed by equations (\ref{spl})-(\ref{spl1}) (first row), trajectories of the  weights $w_{j1}^{(\pm)}(t)$ as functions of time used to compute poles $s_j^{(\pm)}(t)$ for $j=1,2,3$ (second row), and  the corresponding absolute approximation error $|s_j^{(\pm)}(t) - z_j^{(\pm)}(t)|$ (third row) of  the extended solution  of the Burgers equation (\ref{burdis}) subject to the initial condition (\ref{initbur})   for $\nu=0.1$, $\theta=0$ (left) and $\theta=\pi/4$ (right)   for time interval $[0,4]$.   Orange curves show the dynamics of branch points $ \zeta^{(\pm)}(t)$ (for $\nu=0$) as given by (\ref{brpoints}). The pole curves approaches the branch points  curve asymptotically as $t\rightarrow 0+$. 
 % With gray dashed curves we marked the estimated location of singularities computes as poles of  (\ref{solburvis}). 
 }
  \label{fig_bur_pol_w}
\end{figure}

When the parameters $\nu$ and $\theta$ are positive, shock formation does not occur; in the complex-plane interpretation, this corresponds to the singularities never reaching the real axis. Moreover, in this case the Burgers equation \eqref{burdis} with initial condition \eqref{initbur} admits an explicit solution given as a ratio of two series 
\begin{equation}\label{solburvis}
    u(x,t)=-2 \nu \, \mathrm{e}^{\mathrm{i} \theta} \frac{\varphi_x}{\varphi}, \quad \varphi(x,t):=  \sum\limits_{k=-\infty}^\infty (-1)^k I_k \left( \frac{1}{2 \nu \, \mathrm{e}^{\mathrm{i} \theta}}\right) \mathrm{e}^{\mathrm{i} k x-\nu \mathrm{e}^{\mathrm{i} \theta} k^2 t} ,
\end{equation}
where $I_k$ are the modified Bessel functions of the first kind. This solution is derived from the famous Hopf-Cole transform \cite{C51}. As follows from the Painlevé property, the analytic continuation of $u(x,t)$ for each $t>0$ is a meromorphic function in the complex plane with an infinite sequence of conjugate pole pairs. We refer to \cite{VL22} regarding this property for \textbf{the viscous Burgers equation}, i.e. with  $\nu>0$ and $\theta =0$ in (\ref{burdis}), for any initial data. If $\nu>0$ and $\theta >0$ in (\ref{burdis}), we obtain \textbf{the dispersive  Burgers equation} \cite{CG15, SC96}. Numerically, we observe behavior of the extended solution $u(z,t)$ analogous to the real positive-viscosity case, but we have not been able to find a corresponding discussion of the Painlevé property for complex viscosity in the literature. Unfortunately, for finite time $t$ there is no known closed-form formula for all poles $z(t)$ in terms of elementary or standard special functions. The pole locations are defined implicitly by the transcendental equation $\varphi(z,t)=0$, where $\varphi(z,t)$ is defined in (\ref{solburvis}), and can be computed numerically. Employing the root funding, we compute locations of  three closest to the real axis (in the upper and lower half-planes) pairs of complex-conjugated poles of $u(z,t)$ and denote them by $z_j^{(\pm)}(t)$, $j=1,2,3$.  Let also $s_j^{(\pm)}(t)$, $j=1,2,3$ be the estimated locations of those poles computed by equations (\ref{spl})-(\ref{spl1}). In Figure \ref{fig_bur_pol_w}, we present trajectories of $s_j^{(\pm)}(t)$, $j=1,2,3$ (first row) in the complex plane as functions of time for the interval $t \in [0,4]$, as well as the corresponding absolute approximation errors  $|s_j^{(\pm)}(t) - z_j^{(\pm)}(t)|$, $j=1,2,3$ (third row).
%, where $z_j^{(\pm)}$ are computed as poles of $u(z,t)$ in (\ref{solburvis}).
Note that from our numerical experiments (Figure \ref{fig_bur_pol_w}, third row) follows that  poles located more close to the real axis are reconstructed with smaller approximation error, which is explained with the fact that as inputs for the algorithm we use data from the real axis.

 \begin{figure}[h!]
    \centering
      \begin{subfigure}[b]{0.24\linewidth}
    \includegraphics[width=1.1\linewidth]{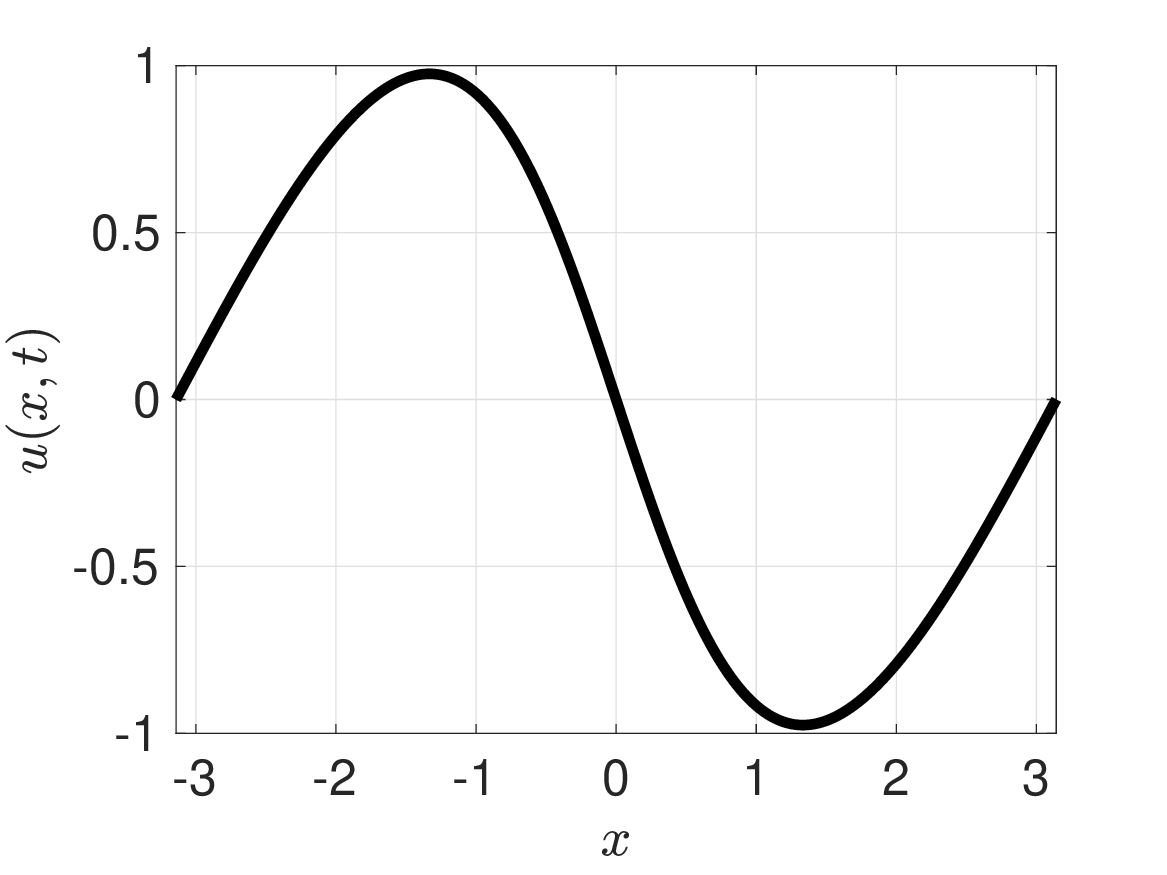}
   \caption{$t=0.25$}
  \end{subfigure}
  \begin{subfigure}[b]{0.24\linewidth}
    \includegraphics[width=1.1\linewidth]{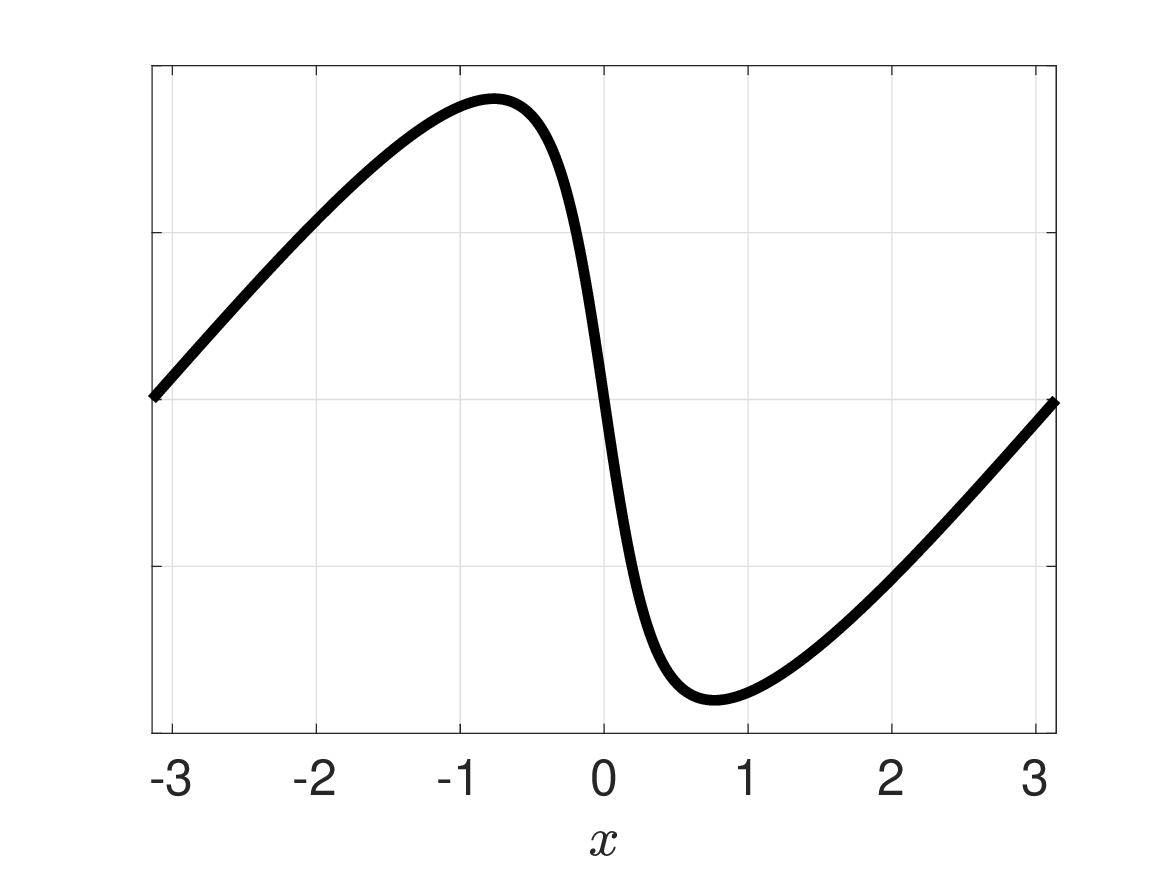}
    \caption{$t_s=1$}
  \end{subfigure}
    \begin{subfigure}[b]{0.24\linewidth}
    \includegraphics[width=1.1\linewidth]{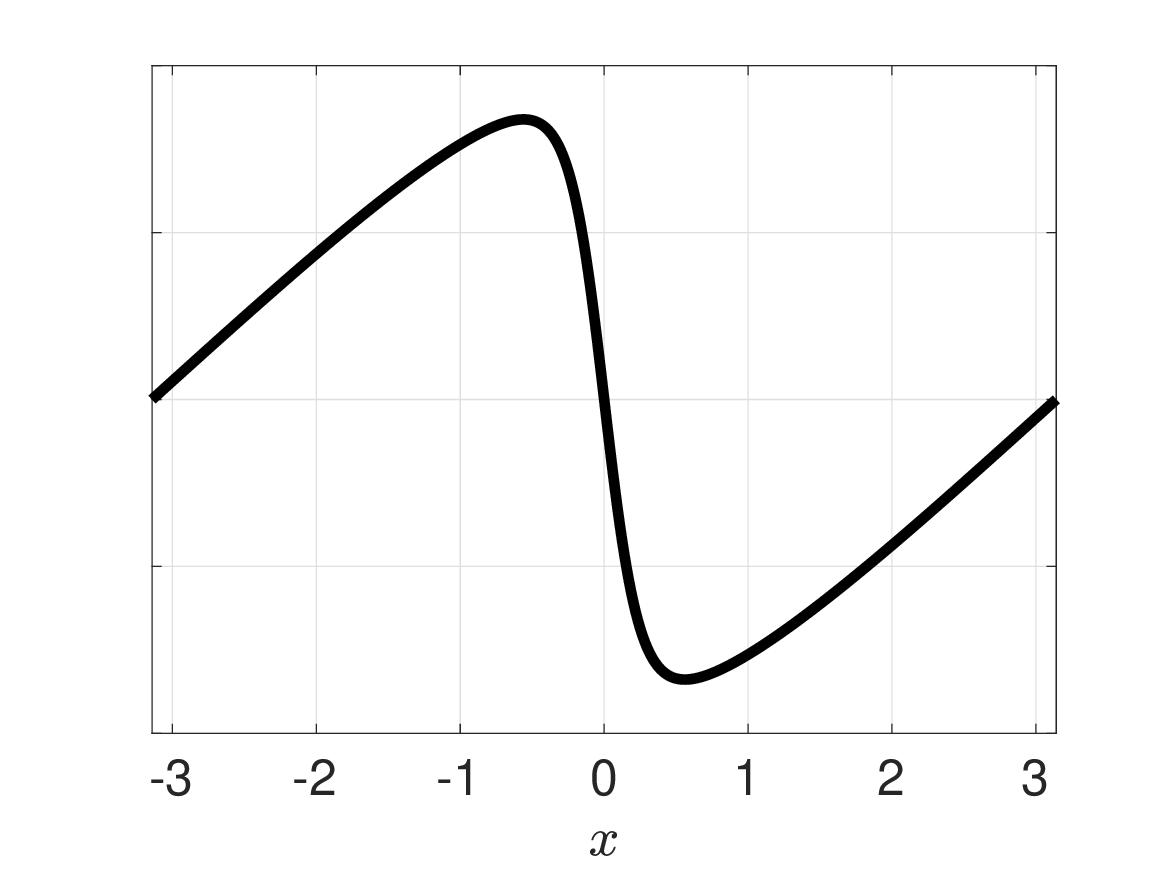}
   \caption{$t=1.5$}
  \end{subfigure}
  \begin{subfigure}[b]{0.24\linewidth}
    \includegraphics[width=1.1\linewidth]{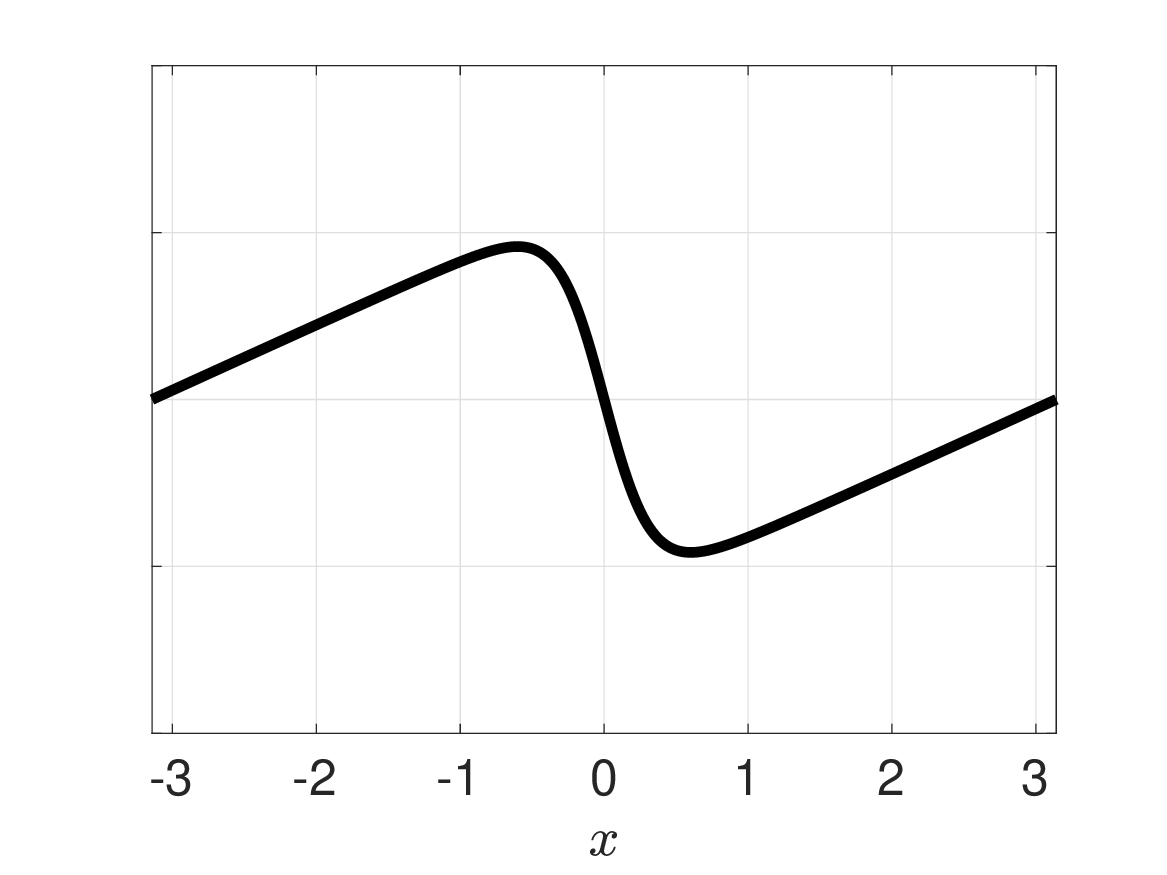}
    \caption{$t=4$}
  \end{subfigure}
  \caption{ The solution $u(x,t)$  
  %and  the phase portraits of its extension   $\Phi(z,t)$ as in (\ref{nnapr})  into the complex plane (second row)
  of the Burgers equation (\ref{burdis}) subject to the initial condition (\ref{initbur}) for $\nu=0.1$ and $\theta=0$ (viscous case) for time $t=0.25$, $t_s=1$ and $t=1.5, 4$.  The shock does not occur. }
  \label{bur_vis_fig}
\end{figure}

Let us first consider the viscous Burgers equation. Similarly as in Example \ref{ex32}, the initial condition \eqref{initbur} is an entire function upon complexification and the poles of the analytic continuation of $u(x,t)$ are ``born'' at infinity for small $t > 0$. They travel toward the real axis along the imaginary axis at a slower rate than in the inviscid case $\nu = 0$ (see Figure \ref{fig_bur_pol_w}, first row (a)). As the solution steepens (Figure \ref{bur_vis_fig}), the singularities approach the real axis, but never reach it, and then eventually  recede back toward $\pm \mathrm{i}\infty$ (Figure \ref{fig_bur_pol_w}, first row (a)), so that the shock loses amplitude. For $\theta > 0$, the poles no longer remain on the imaginary axis. They retain the same qualitative behavior, but travel toward the real axis not along the imaginary axis, but along a slanted line in the complex plane whose direction is set by the parameter $\theta$ (see Figure \ref{fig_bur_pol_w}, first row (b)). In Figure \ref{burviscsol}, we observe three pairs of complex-conjugate poles $s_j^{(\pm)}(t)$, $j = 1, 2, 3$, approaching the real axis along the imaginary axis for $\theta = 0$ (first row (a)), and along a ``rotated'' imaginary axis for $\theta = \pi/4$ (first row (b)), together with the corresponding pole-scaled approximation error \eqref{errorrel} (second row) at time $t_s = 1$.

\begin{figure}[h!]
  \centering
  \begin{subfigure}[b]{0.45\linewidth}
   \includegraphics[width=1\linewidth]{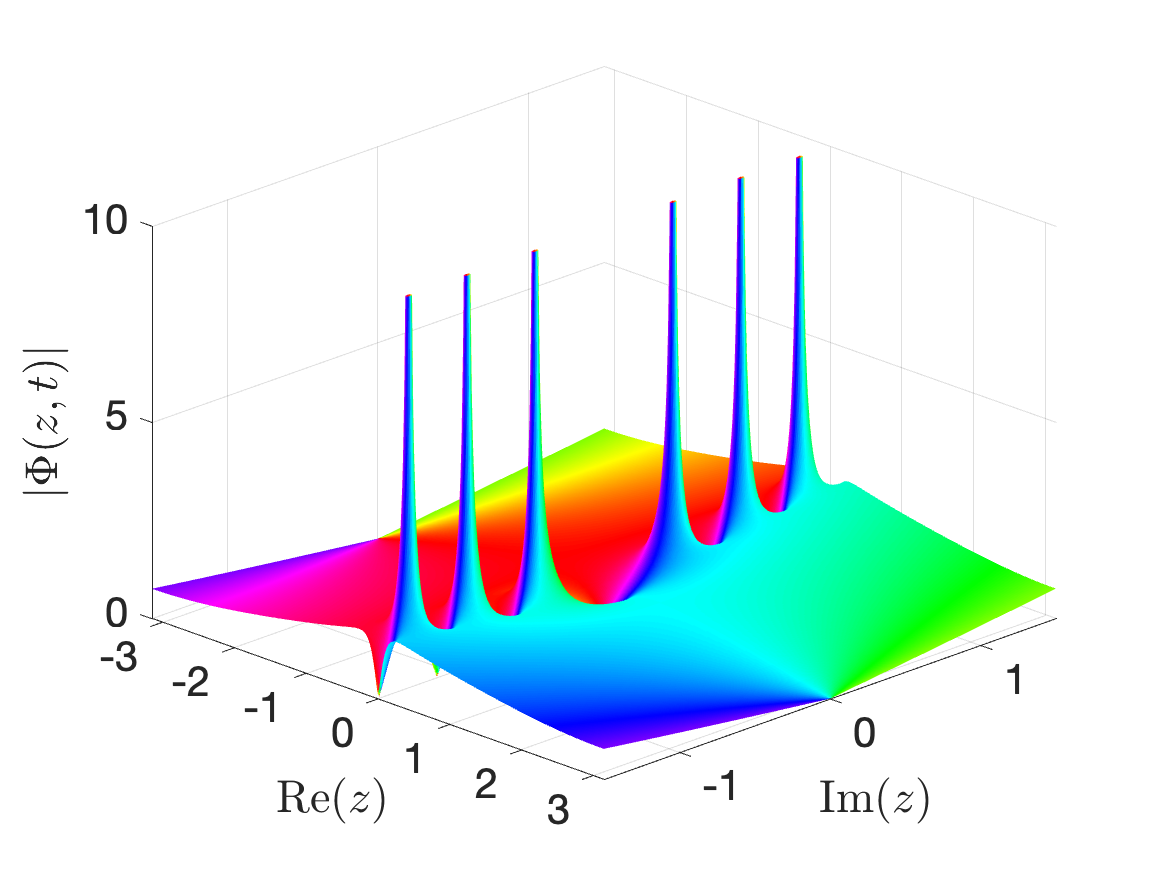}
    %\caption{$t=0.25$}
  \end{subfigure}
  \begin{subfigure}[b]{0.45\linewidth}
   \includegraphics[width=1\linewidth]{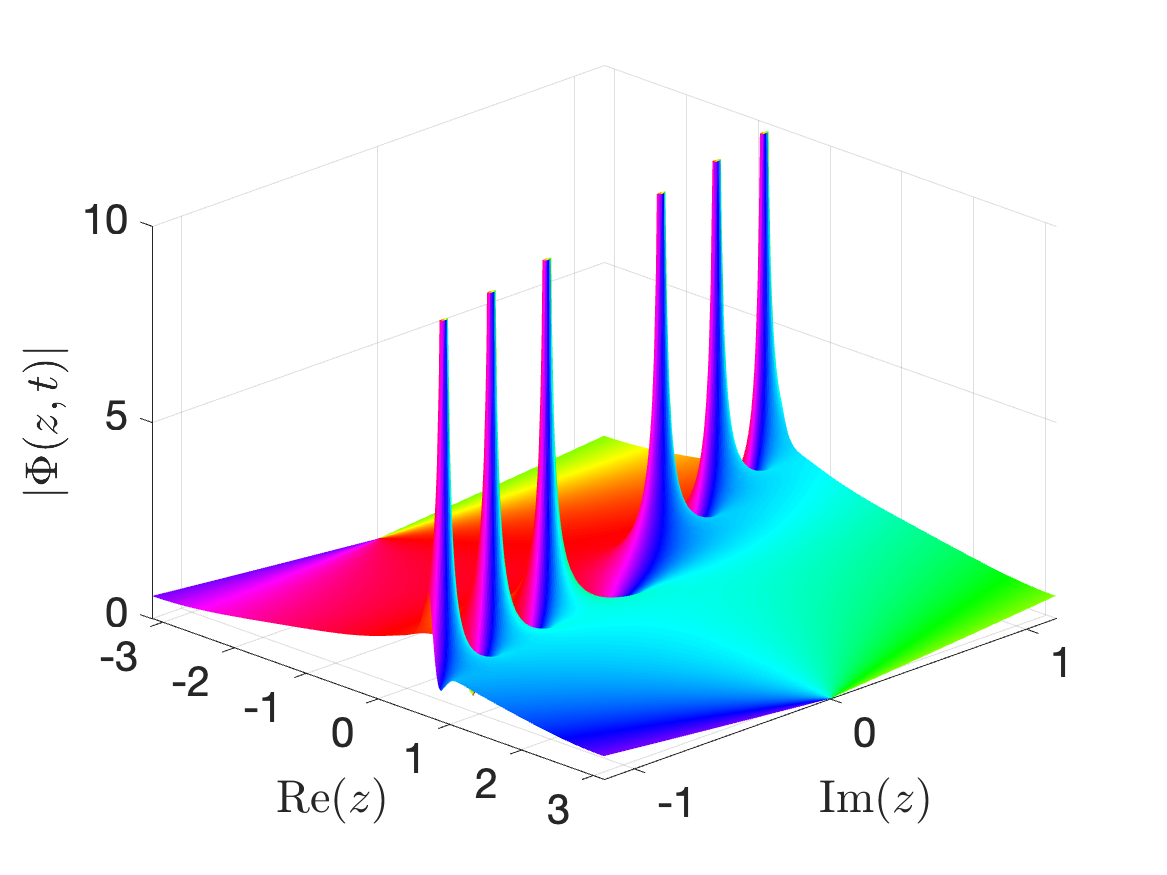}
   % \caption{$t_s=1$}
  \end{subfigure}
  \begin{subfigure}[b]{0.45\linewidth}
     \includegraphics[width=0.91\linewidth]{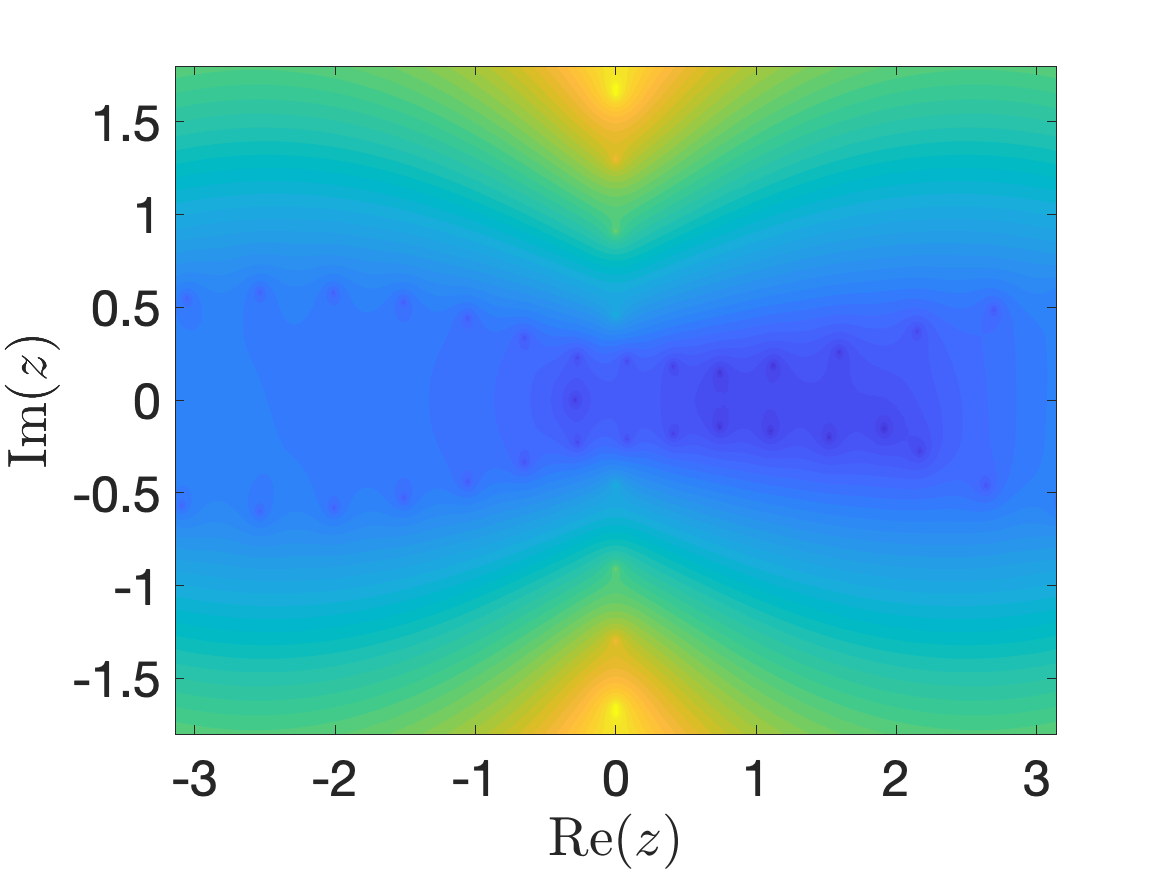}
    \caption{$\theta=0$}
  \end{subfigure}
  \begin{subfigure}[b]{0.45\linewidth}
  \includegraphics[width=0.91\linewidth]{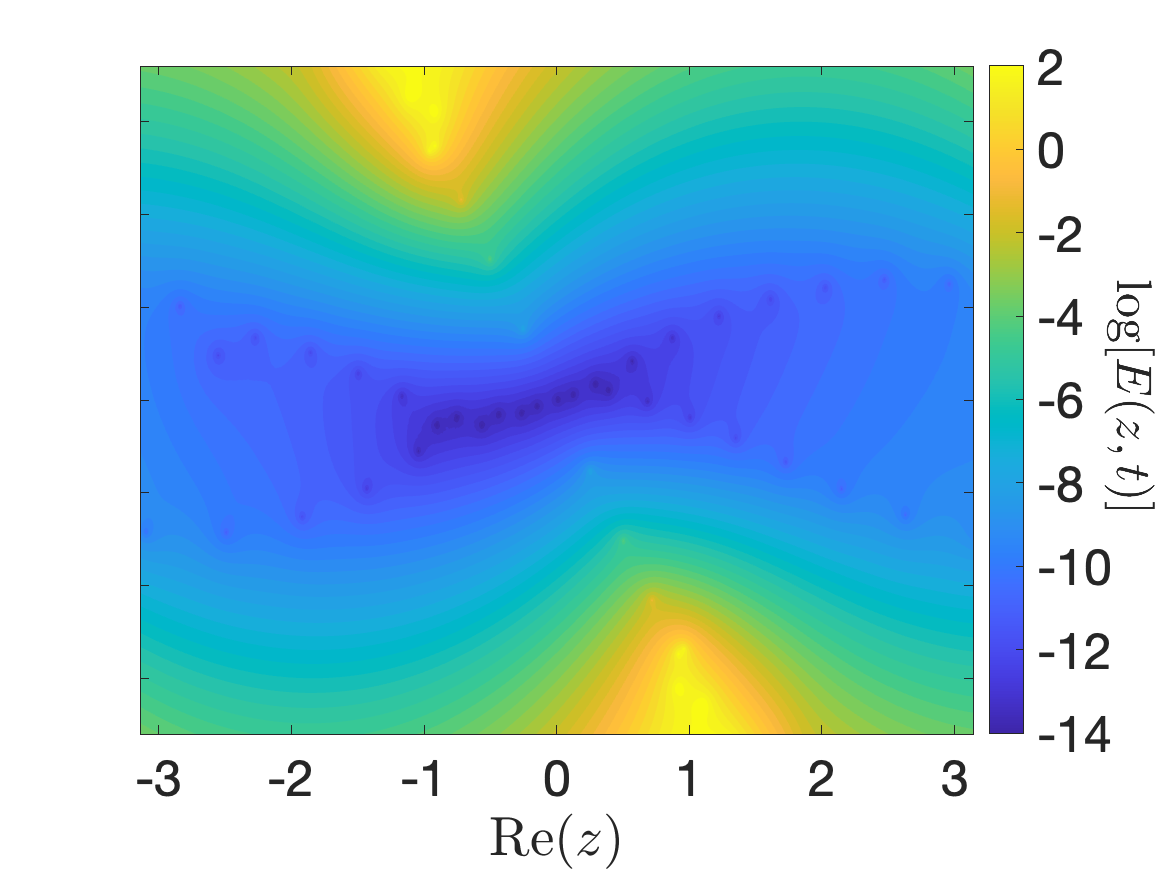}
   \caption{$\theta=\pi/4$}
  \end{subfigure}
  \caption{The analytic landscapes of the numerical analytic continuation $\Phi(z,t)$ as in (\ref{nnapr}) of the solution $u(x,t)$ (first row) and the corresponding pole-scaled approximation error (\ref{errorrel}) in the logarithmic scale  (second row) of the Burgers equation (\ref{burdis}) subject to the initial condition (\ref{initbur})   for $\nu=0.1$ and $\theta=0$ (left)  and   $\theta=\pi/4$  (right) at time  $t_s=1$. The shock does not occur. }
  \label{burviscsol}
\end{figure}

Finally, we study dynamics of weights $w_{j1}^{(\pm)}(t)$ and biases $b_{j1}^{(\pm)}(t)$ of the hidden layers of the components $\Phi^{(\pm)}$ that we used in (\ref{spl})-(\ref{spl1}) to compute singularities $s_{j}^{(\pm)}(t)$ for $j=1,2,3$. We choose again $\nu=0.1$ and consider the case $\theta=0$.
Let us first fix parameters  $C_{j0}^{(-)}(t)=1$ and $C_{j0}^{(+)}(t)=-1$ for  $j=1,2,3$ and each $t \in [0,4]$. 
According to our computations,   biases are constant with the values $b_{j1}^{(-)}(t)=-0.236067977499790$, $b_{j1}^{(+)}(t)=4.236067977499788$,  and   weights  satisfy the properties $w_{j1}^{(-)}(t)=-w_{j1}^{(+)}(t)$, $w_{j1}^{(-)}(t)<0$ and $w_{j1}^{(+)}(t)>0$ for $j=1,2,3$. Analyzing formulas (\ref{spl})-(\ref{spl1}), we conclude that poles must be purely imaginary and complex conjugated.   Moreover, since $b_{j1}^{(-)}(t)+z_0<0$ ($b_{j1}^{(+)}(t)+z_0>0$), and $w_{j1}^{(-)}(t)$ is negative and decreasing ($w_{j1}^{(+)}(t)$ is positive and increasing), we obtain that  $\mathrm{Im}(s_j^{(-)} (t))$ decreases ($\mathrm{Im}(s_j^{(-)} (t))$ increases), which completely coincides with our numerical experiments (Figure \ref{fig_bur_pol_w}, first and second rows (a)).
For other values of the parameters $C_{j0}^{(\pm)}(t)$ — for instance, choosing them at random as $C_{j0}^{(\pm)}(t) \in [-1.5,-0.5] + [-1.5,-0.5]\mathrm{i}$ for $j = 1,2,3$ — the weights $w_{j1}^{(\pm)}(t)$ and biases $b_{j1}^{(\pm)}(t)$ will in general be complex, but the overall behavior remains unchanged (see Figure \ref{figwbb}). For $\theta = \pi/4$ and the choice of parameters $C_{j0}^{(-)}(t) = 1$ and $C_{j0}^{(+)}(t) = -1$ for $j = 1,2,3$, we obtain the same values of the biases, $b_{j1}^{(-)}(t) = -0.236067977499790$ and $b_{j1}^{(+)}(t) = 4.236067977499788$, as for $\theta = 0$. The weights again satisfy $w_{j1}^{(-)}(t) = -w_{j1}^{(+)}(t)$, but they are now complex (see Figure \ref{fig_bur_pol_w}, second row (b)).

\begin{figure}[h!]
  \centering
  \begin{subfigure}[b]{0.45\linewidth}
    \includegraphics[width=1\linewidth]{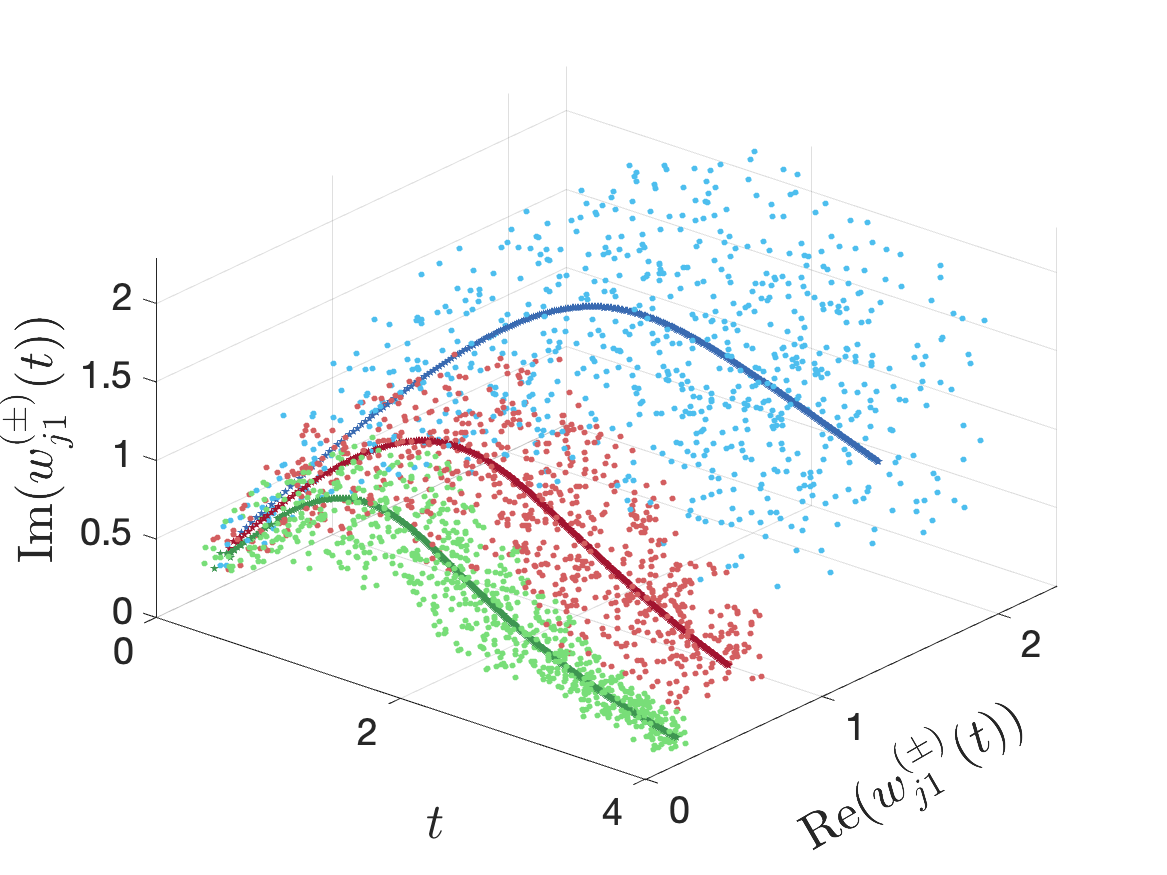}
   % \caption{$t=0.15$}
  \end{subfigure}
  \begin{subfigure}[b]{0.45\linewidth}
    \includegraphics[width=1\linewidth]{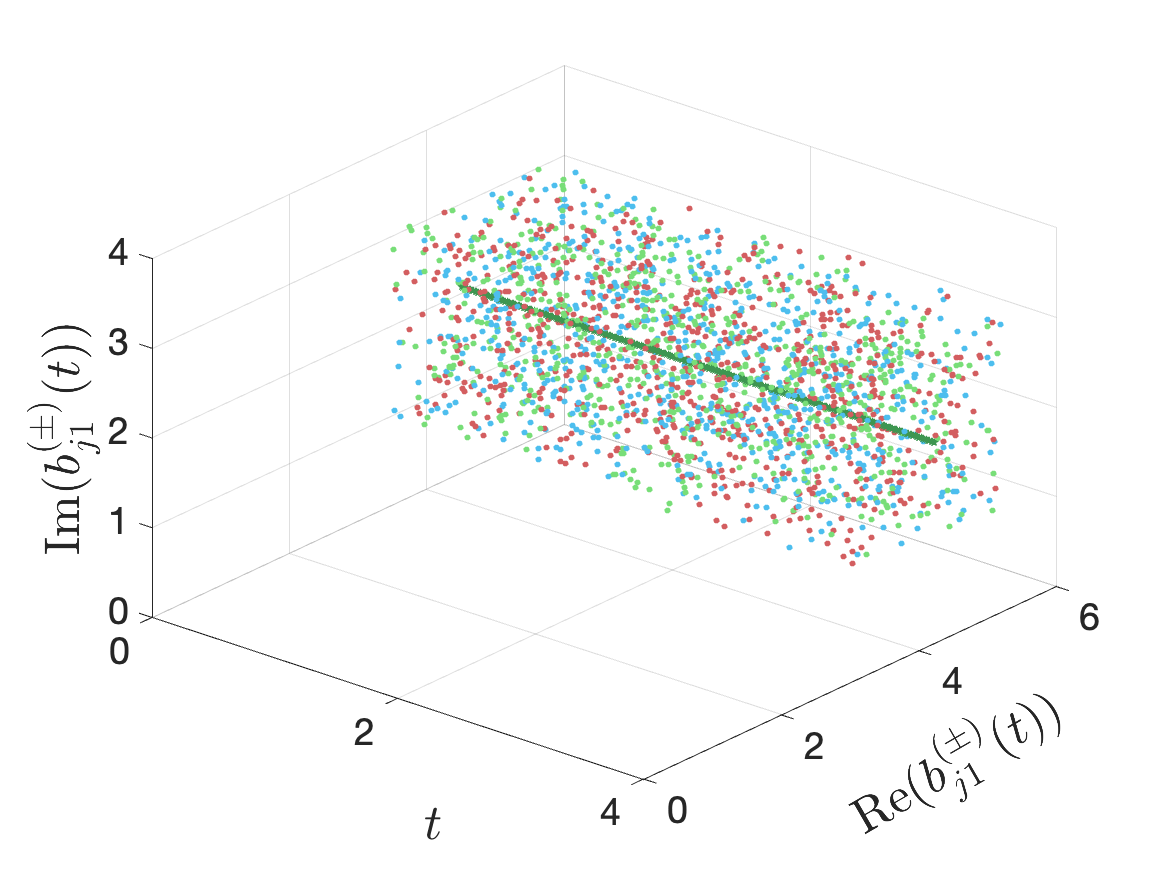}
   % \caption{$t=0.503$}
  \end{subfigure}
  \caption{
    Locations of weights $w_{j1}^{(\pm)}(t)$ (left) and biases  $b_{j1}^{(\pm)}(t)$ (right) used in (\ref{spl})-(\ref{spl1}) for computation of singularities  $s_{j}^{(\pm)}(t)$ for $j=1,2,3$ of the extended solution  of the Burgers equation (\ref{burdis}) subject to the initial condition (\ref{initbur})   for $\nu=0.1$ and  $\theta=0$  in the time interval $[0,4]$. 
%Dark blue -- $w_{11}^{(-)}(t)$ and $b_{11}^{(-)}(t)$ computed for $C_0^{(-)}(t)=-1$ and $z_0=-2$,  light blue -- $w_{11}^{(-)}(t)$ and $b_{11}^{(-)}(t)$ computed for random choice of parameters  $C_{10}^{(-)}(t) \in [-1.5,-0.5]$ and $z_0=-2\mathrm{i}$.  
Dark blue: $w_{11}^{(\pm)}(t)$ and $b_{11}^{(\pm)}(t)$, dark red: $w_{21}^{(\pm)}(t)$ and $b_{21}^{(\pm)}(t)$, dark green: $w_{31}^{(\pm)}(t)$ and $b_{31}^{(\pm)}(t)$,
  %, dark purple -- $w_{14}^{(\pm)}(t)$ and $b_{14}^{(\pm)}(t)$, 
  computed with fixed parameters $C_{j0}^{(\pm)}(t)=-1-\mathrm{i}$;  light blue: $w_{11}^{(\pm)}(t)$ and $b_{11}^{(\pm)}(t)$, light red: $w_{21}^{(\pm)}(t)$ and $b_{21}^{(\pm)}(t)$, light green: $w_{31}^{(\pm)}(t)$ and $b_{31}^{(\pm)}(t)$,
  %, light purple -- $w_{14}^{(\pm)}(t)$ and $b_{14}^{(\pm)}(t)$,
   computed with  random choice of parameters  $C_{j0}^{(\pm)}(t) \in [-1.5,-0.5]+ \mathrm{i} \,  [-1.5,-0.5]$. }
  \label{figwbb}
\end{figure}

\end{example}

As is well known, shock formation is associated with a finite-time loss of regularity of the solution.   From the viewpoint of numerical analytic continuation, this phenomenon is characterized by the motion of complex singularities toward the real axis.  
%As the nearest singularity approaches the real axis, the width of the analyticity strip decreases. 
As the complex singularities approach and eventually reach the real axis, the solution loses analyticity, resulting in shock formation.  Unlike blow-up, where the solution itself becomes unbounded, shock formation is characterized by the boundedness of the solution together with the blow-up of its spatial gradient. This distinction is reflected in the type of complex singularities approaching the real axis. Blow-up is typically associated with unbounded singularities, such as poles or logarithmic branch points, whereas shock formation is caused by bounded singularities, most notably algebraic branch points.

In the case of the nonlinear Burgers equation (\ref{burdis}), viscosity $\nu$ regulates the formation of a shock.
From the perspective of numerical analytic continuation, it means that viscosity controls the evolution of the complex singularities.
In the absence of viscosity, the nearest complex singularities approach the real axis and eventually collapse on it, initiating shock formation.
 The positive viscosity slows the motion of the singularities toward the real axis,  preventing shock formation.

%The collision of the singularities with the real axis marks the onset of the gradient catastrophe and the formation of a shock.

\subsection{Nonlinear Schrödinger equation (NLS)}
\label{secnls}

In this subsection, we investigate the emergence of rogue waves through the analysis of solutions to the NLS.

\begin{example}
\label{nlsex}
In our last example, we investigate the NLS of the following form 
\begin{equation}\label{nls}
  \mathrm{i}  u_t+ u_{xx}+2 |u|^2 u=0,
\end{equation}
with  $ x \in [-\pi,\pi)$, $t\in \rr$ and the $2\pi$-periodic Akhmediev breather as  the initial condition \cite{GS18}. The focusing NLS equation \eqref{nls} is the simplest universal model describing the modulational (or Benjamin–Feir) instability of quasi-monochromatic waves in weakly nonlinear media, widely regarded as the main physical mechanism behind the appearance of rogue waves in nature. We refer to \cite{GS18} for more details on this phenomenon, and in this paper we offer yet another interpretation: using the dynamics of the poles of the extended solutions in the complex plane, we study the formation of rogue waves and their properties.

The $2\pi$-periodic Akhmediev breather is given by 
(see, for example, \cite{GS18})
\begin{equation}\label{solnls}
  u(x,t)=    \left( \frac{\cosh(\sqrt{3} t+2\pi\mathrm{i}/3)+\sqrt{3} \cos(x)/2}{\cosh(\sqrt{3} t) -\sqrt{3} \cos(x)/2 } \right) \mathrm{e}^{2 \, \mathrm{i}  t+2 \, \mathrm{i} \, \pi/3}
\end{equation}
and is $2\pi$-periodic in space variable $x$ and localized in time  variable $t$.
 The Akhmediev breather (\ref{solnls}) is first-order rogue wave and its maximum height in the physical space occurs at $x=t=0$ (see Figure \ref{fig_nls_2}).

\begin{figure}[h!]
  \centering
      \begin{subfigure}[b]{0.24\linewidth}
    \includegraphics[width=1.1\linewidth]{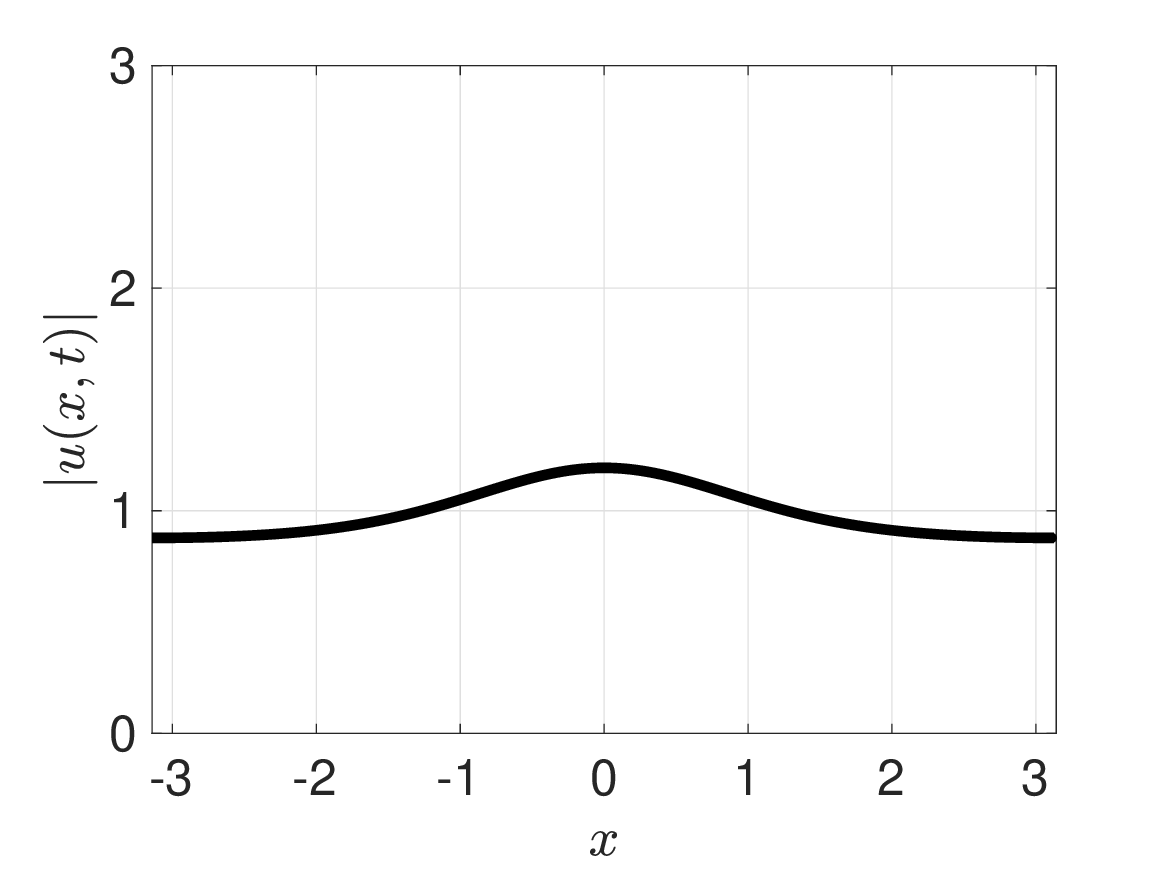}
    \caption{$t=-1$}
  \end{subfigure}
    \begin{subfigure}[b]{0.24\linewidth}
    \includegraphics[width=1.1\linewidth]{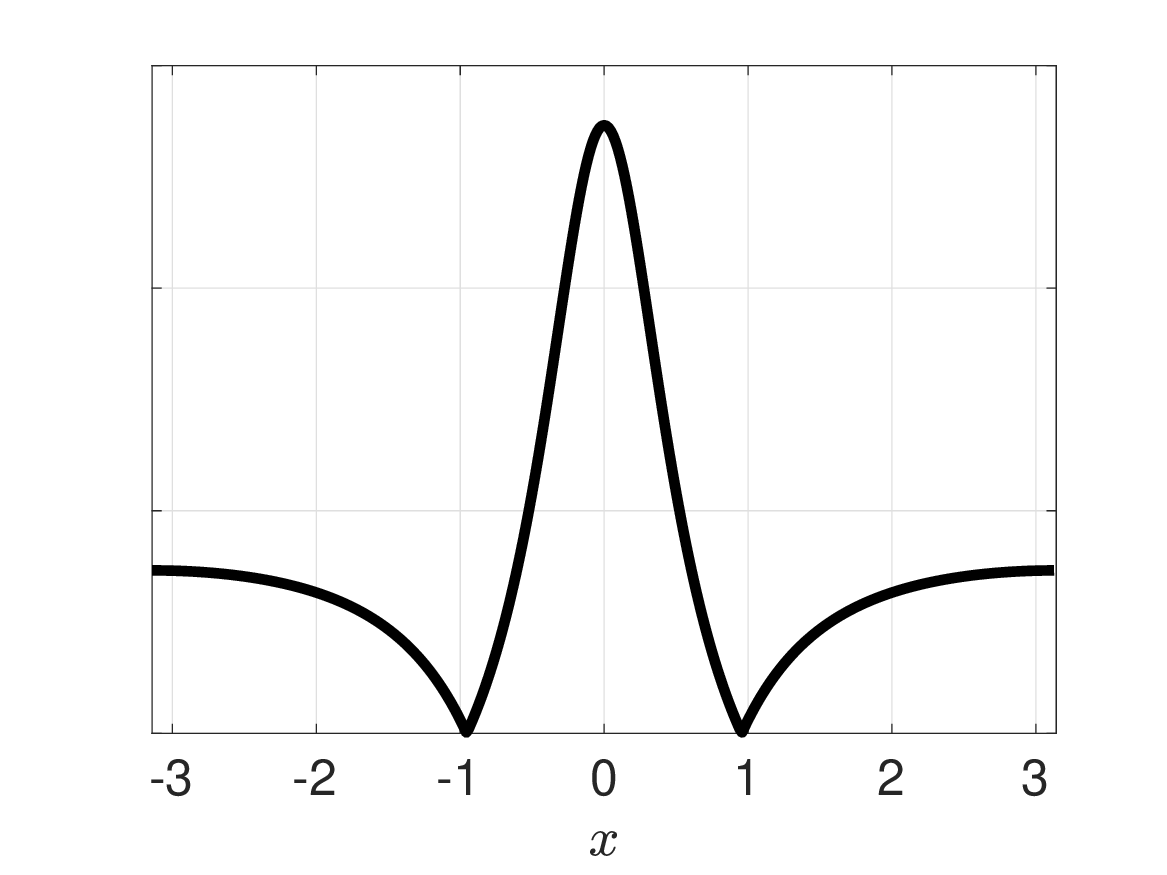}
    \caption{$t=0$}
  \end{subfigure}
  \begin{subfigure}[b]{0.24\linewidth}
    \includegraphics[width=1.1\linewidth]{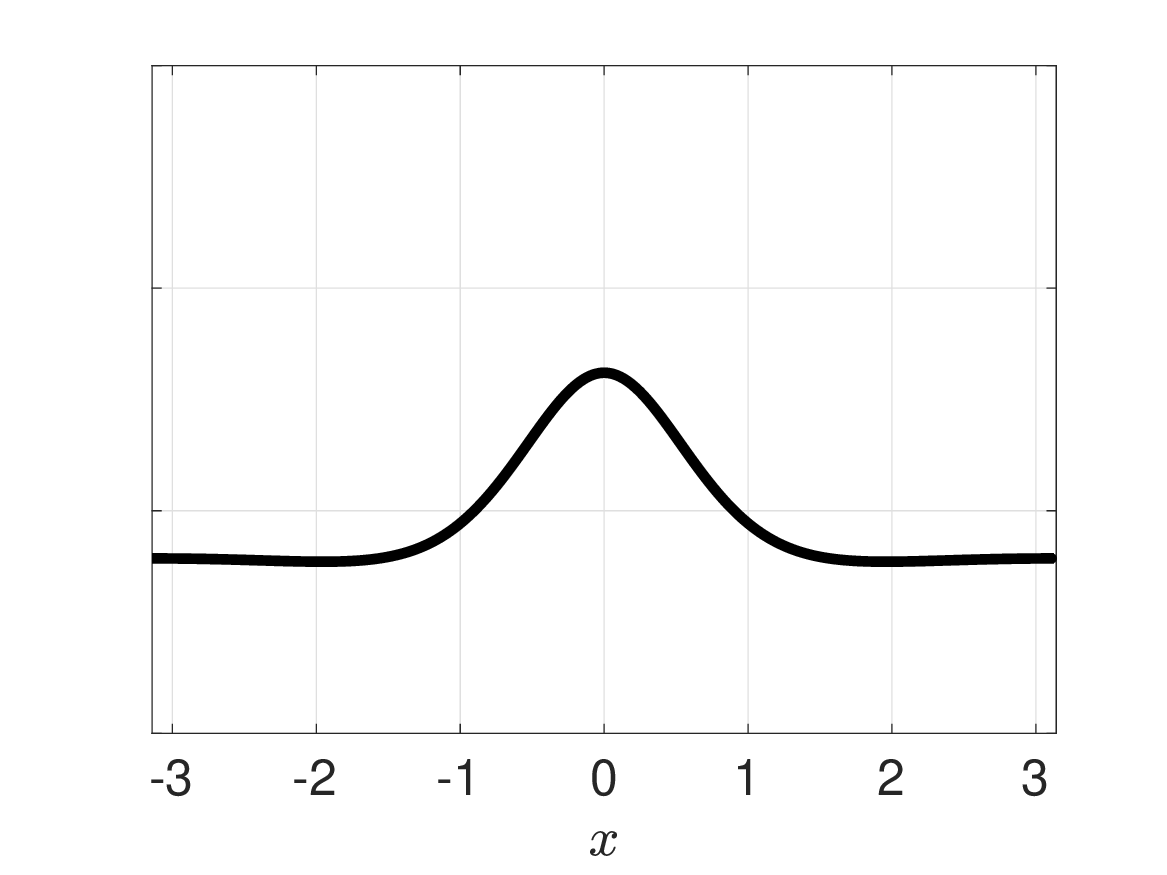}
    \caption{$t=0.5$}
  \end{subfigure}
  \begin{subfigure}[b]{0.24\linewidth}
    \includegraphics[width=1.1\linewidth]{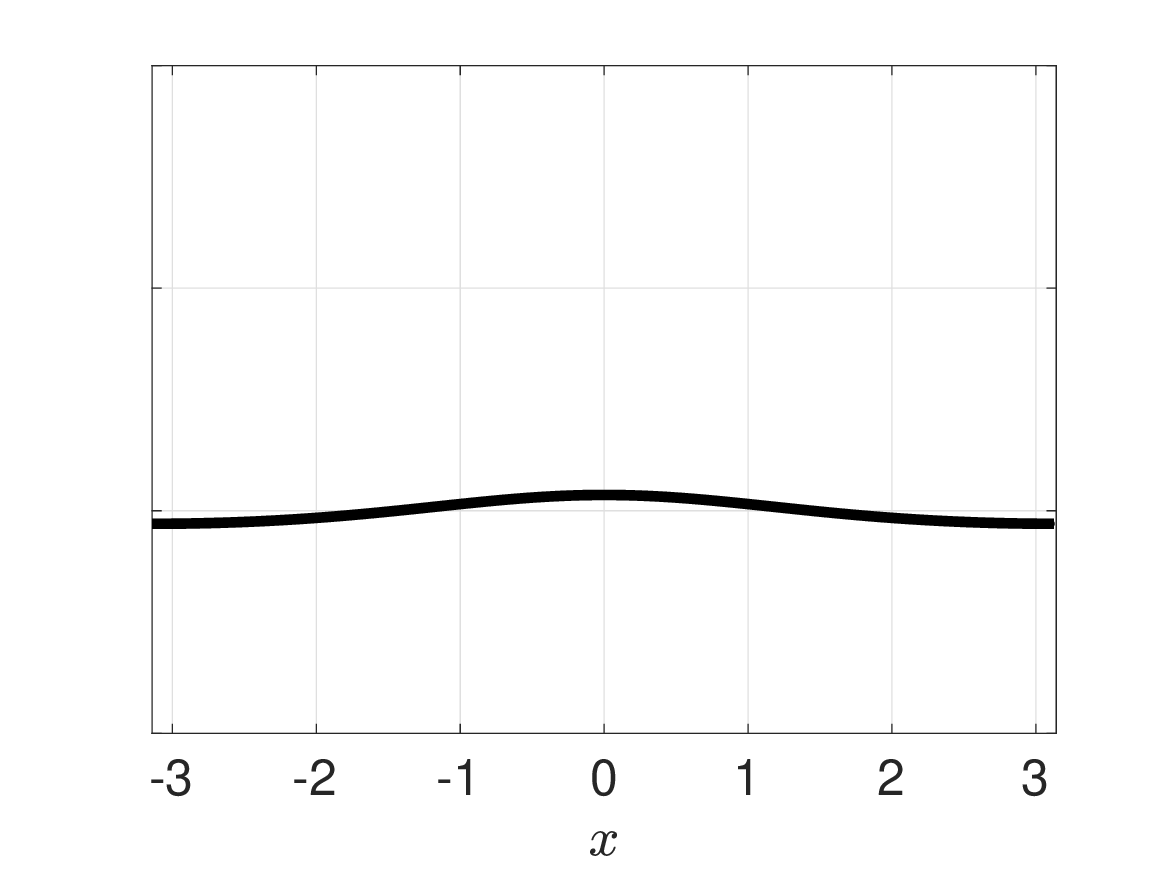}
    \caption{$t=1.5$}
  \end{subfigure}
  \caption{The solution $u(x,t)$ of the NLS (\ref{nls}) subject to the initial condition $u(x,0)$ as in  (\ref{solnls})  for time $t=-1,0,0.5,1.5$. The rogue wave achieves its maximum height at the point $(x,t)=(0,0)$. }
  \label{fig_nls_2}
\end{figure}

We investigate the analytic continuation $u(z,t)$ in the domain $z \in D= [-\pi, \pi)\times  \mathrm{i} \, (-5, 5)$ for time interval $t \in [-2,2]$. In the first step of our method, we solve the  NLS equation (\ref{nls}) with the initial condition $u(x,0)$ as in (\ref{solnls}) by   the split-step Fourier method employing  samples $u\left( \frac{ \pi j}{ n} ,0 \right)$, $j=-n,\dots, n-1$ with $n=200$.
In the second step, we construct an approximation of the extended solution  $\Phi(z,t)$ into the complex plane through  (\ref{nnapr}) and compute estimated locations of its complex singularities via (\ref{spl})-(\ref{spl1}).   

\begin{figure}[h!]
    \centering  
     \begin{subfigure}[b]{0.45\linewidth}
    \includegraphics[width=0.92\linewidth]{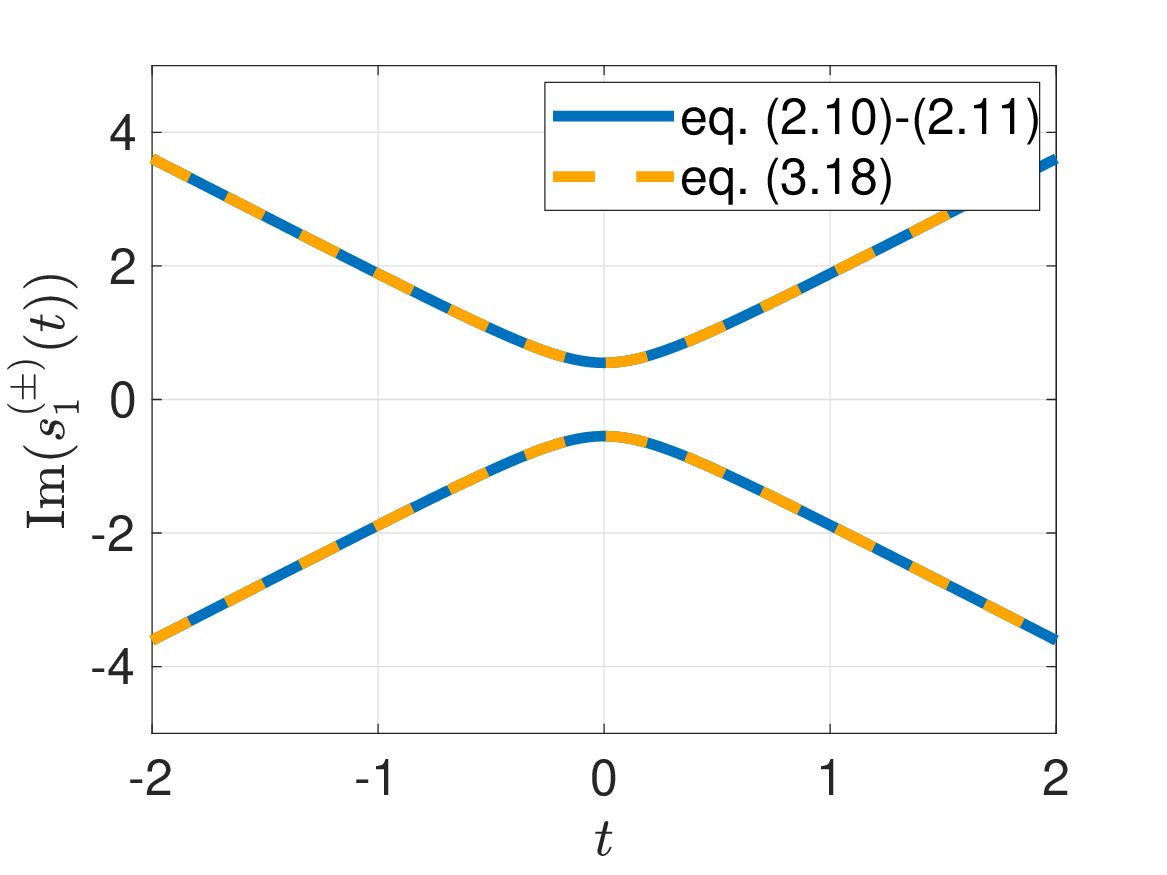}
   %\caption{$t=0.5$}
  \end{subfigure}
     \begin{subfigure}[b]{0.45\linewidth}
     \includegraphics[width=0.92\linewidth]{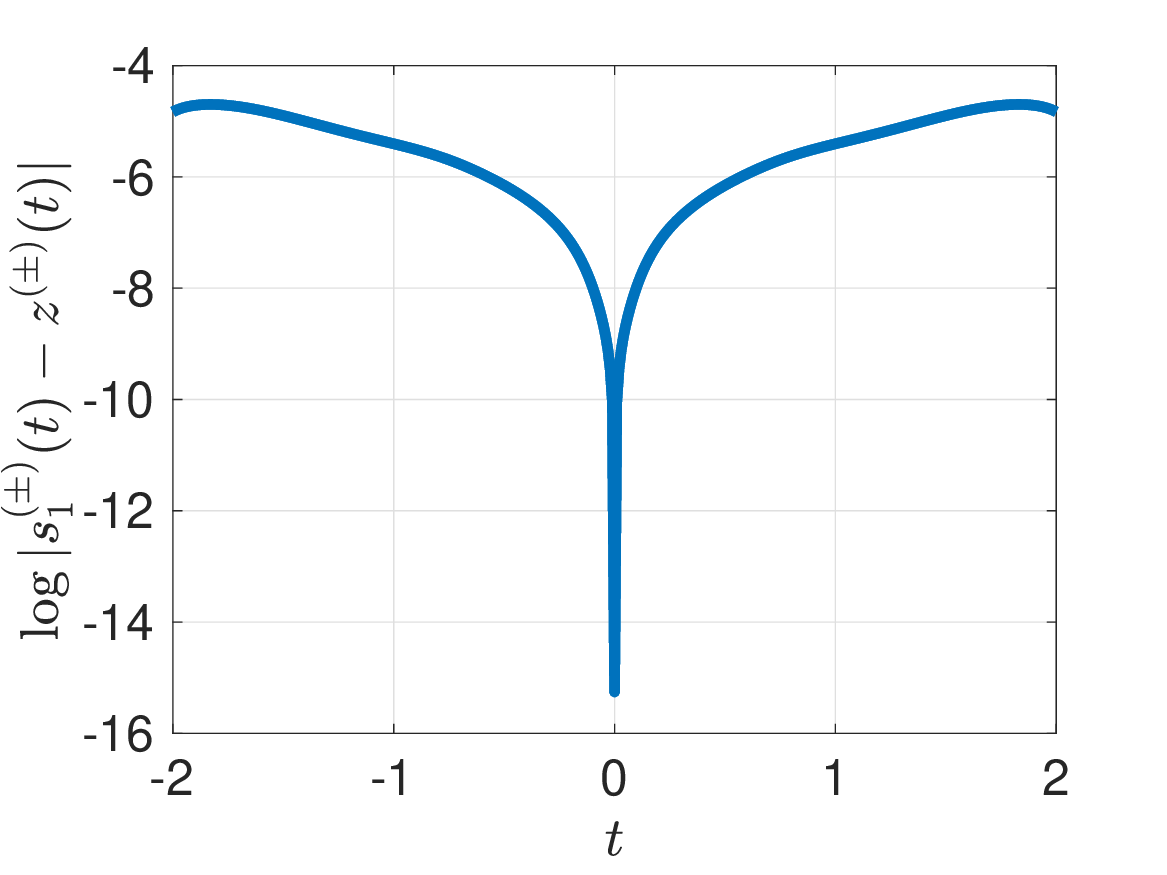}
  % \caption{$t=0.5$}
  \end{subfigure}
      \caption{ Trajectories of complex singularities $s_1^{(\pm)}(t)$ as functions of time (left) and the corresponding absolute approximation errors $|s_{1}^{(\pm)}(t)-z^{(\pm)}(t)|$ in the logarithmic scale (right) of the the extended solution  of  the NLS (\ref{nls}) subject to the initial condition $u(x,0)$ as in  (\ref{solnls}) computed by the equations (\ref{spl})-(\ref{spl1})  and  (\ref{nlspoles}).  }
  \label{sindynfignls}
\end{figure}

\begin{figure}[h!]
  \centering
      \begin{subfigure}[b]{0.45\linewidth}
    \includegraphics[width=1\linewidth]{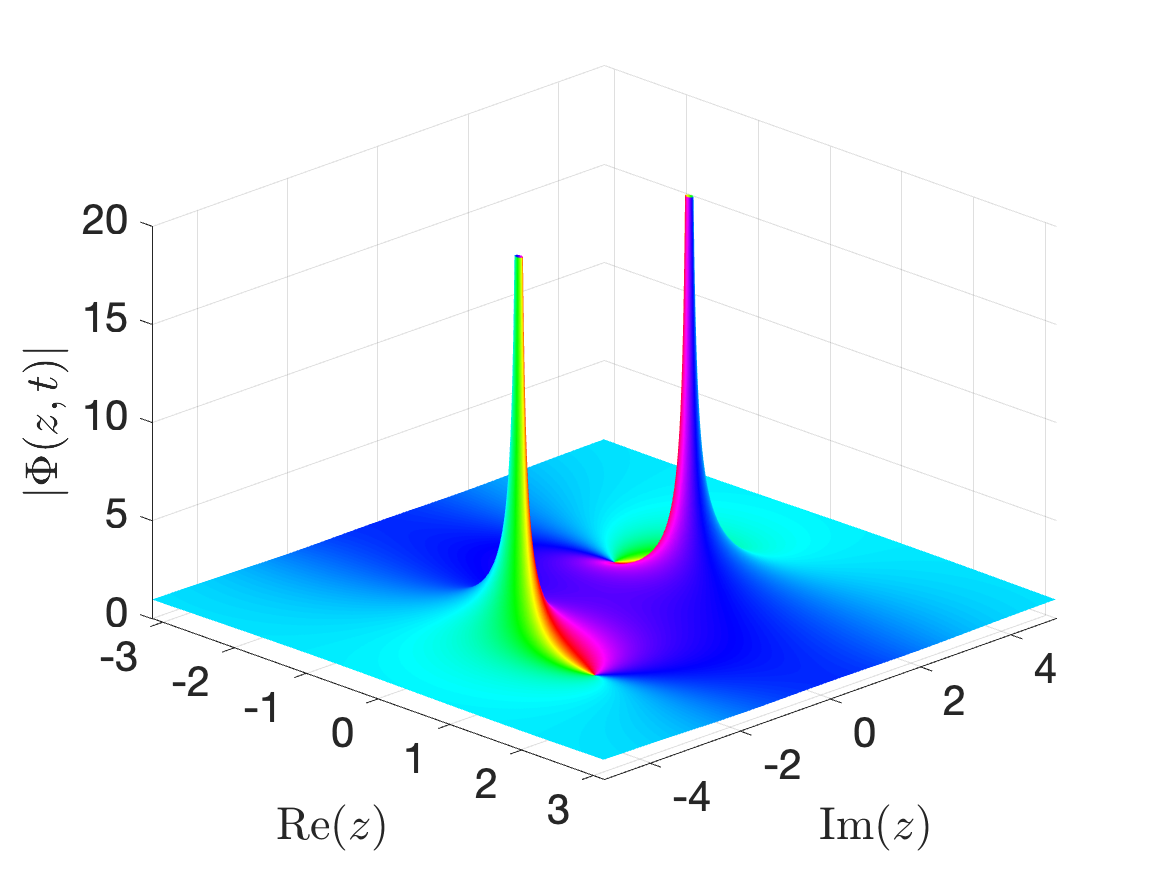}
     \caption{$t=-1$}
  \end{subfigure}
    \begin{subfigure}[b]{0.45\linewidth}
    \includegraphics[width=1\linewidth]{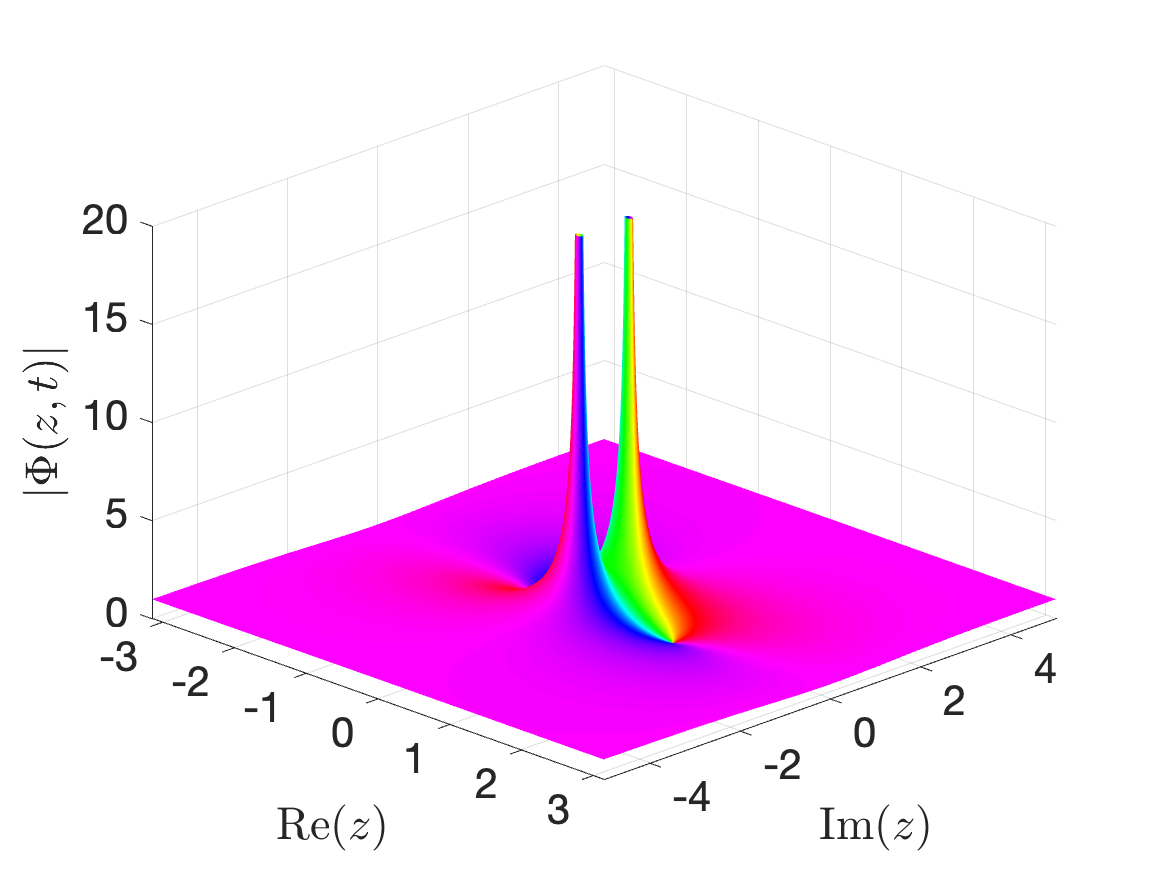}
     \caption{$t=0$}
  \end{subfigure}
  \begin{subfigure}[b]{0.45\linewidth}
    \includegraphics[width=1\linewidth]{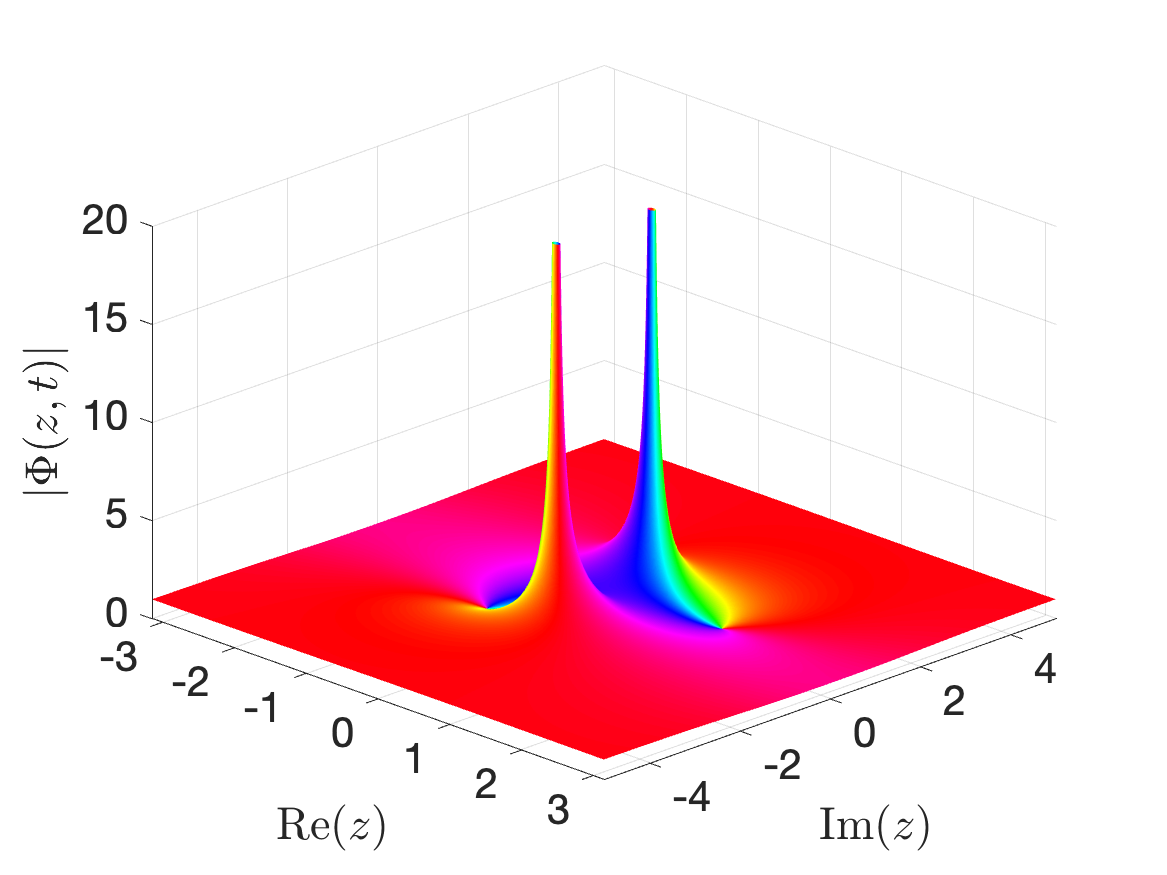}
   \caption{$t=0.5$}
  \end{subfigure}
    \begin{subfigure}[b]{0.45\linewidth}
    \includegraphics[width=1\linewidth]{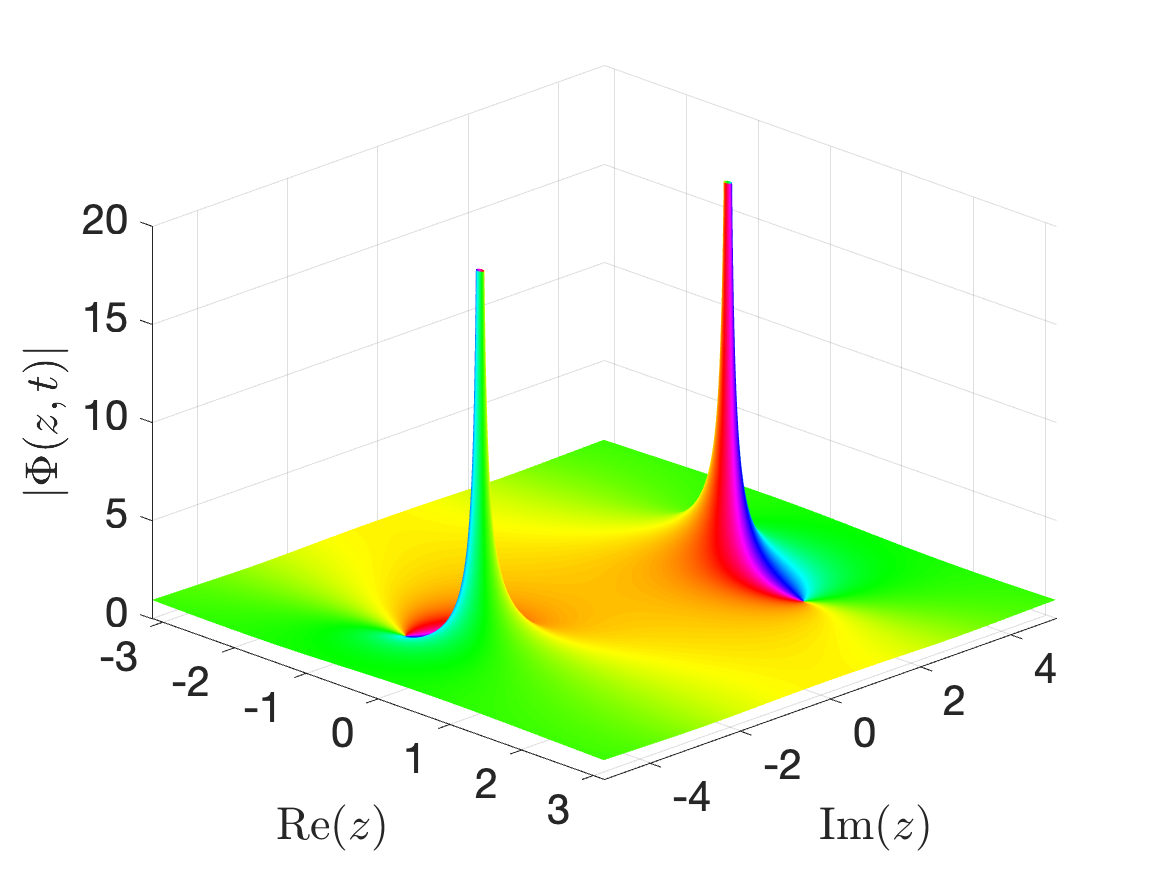}
\caption{$t=1.5$}
  \end{subfigure}
 \caption{The  analytic landscapes of the analytic continuation  $\Phi(z,t)$  as in (\ref{nnapr}) of the  solution  of  the NLS (\ref{nls}) with the initial condition $u(x,0)$ as in  (\ref{solnls}) for time $t=-1,0,0.5,1.5$. At time $t=0$ (b) two poles of the extended solution $u(z,t)$ reach their minimum distance to the real axis. At this point, the rogue wave achieves its maximum height.}
  \label{figpolesnlsp}
\end{figure}

The function $u(z,t)$ has singularities  at $ z^{(\pm)}(t)=\pm \mathrm{i} \, \operatorname{arcosh}\left(\frac{2}{\sqrt{3}} \cosh(\sqrt{3} t) \right) +2 \pi n$ for $n=0,\pm  1, \pm 2,\dots$.
In the domain $D$, $u(z,t)$ has  two simple poles of the form
\begin{equation}\label{nlspoles}
 z^{(\pm)}(t)=\pm \mathrm{i} \, \operatorname{arcosh}\left(\frac{2}{\sqrt{3}} \cosh(\sqrt{3} t) \right).
\end{equation}
The SVD-based procedure from \cite{GGT13} (or Algorithm 2 in \cite{DKD2025}) allows us to compute  the optimal values for parameters  $N^{(\pm)}(t)=M^{(\pm)}(t)=1$. It means that 
each neural network component $\Phi^{(\pm)}$ detects one pole $s_1^{(\pm)}(t)$, respectively, according to the formulas (\ref{spl})-(\ref{spl1}). Trajectories of those poles as well as the corresponding absolute approximation errors $|s_1^{(\pm)}(t)-z^{(\pm)}(t)|$  are presented in Figure   \ref{sindynfignls}. The neural network approximation $\Phi(z,t)$ of the extended solution is shown in Figure \ref{figpolesnlsp}, while the corresponding pole-scaled approximation error, defined by (\ref{errorrel}), is shown in Figure \ref{nlserror} for time $t=-1,0,0.5,1.5$. It should be noted that the errors presented in Figures~\ref{sindynfignls} (left) and~\ref{nlserror} result from the fact that the Akhmediev breather \eqref{solnls} provides only an approximate solution to \eqref{nls} in some small neighborhood of $t=0$ \cite{GS18}.

%Note that the errors presented in Figures \ref{sindynfignls} and \ref{nlserror}  results from the fact that   the Akhmediev breather (\ref{solnls}) is the approximate solution to (\ref{nls}) in some small neighborhood of $t=0$ \cite{GS18}.  
%The neural network approximation $\Phi(z,t)$ of the extended solution  is presented in Figure  \ref{figpolesnlsp} and the corresponding pole-scaled approximation  error as in (\ref{errorrel})   for time $t=-1,0,0.5,1.5$ in Figure \ref{nlserror}. 

As time evolves from $- \infty$ to $+\infty$, the poles $s_1^{(\pm)}(t)$ move along the imaginary axis toward the real axis. With the poles approaching the real axis (Figure \ref{figpolesnlsp} (a)), the rogue wave starts to grow (Figure \ref{fig_nls_2} (a)) and achieves its maximum height at time $t=0$ (Figure \ref{fig_nls_2} (b)), when the poles reach their minimum distance to the real axis (Figure \ref{figpolesnlsp} (b)). Next, the singularities turn around and travel back along the imaginary axis for positive $t$ (Figure  \ref{figpolesnlsp} (c) and (d)). The pair of poles never actually reaches the real axis.  The maximum height of the Akhmediev breather in the physical space occurs at the location $x=0$, which is the real part of the point in the pole trajectories in the complex plane where they turn around. Hence, similarly as for the first-order Peregrine soliton (non-periodic solution) \cite{PCC20}, the following conjecture also holds for the first-order  Akhmediev breather: the spatial locations of the points of maximum heights of a rogue wave in physical space will coincide, or closely correlate, with the real parts of the poles of the rogue wave solutions in the complex plane at points where the pole trajectories reverse directions.

\begin{figure}[h!]
  \centering
  \begin{subfigure}[b]{0.24\linewidth}
    \includegraphics[width=1.1\linewidth]{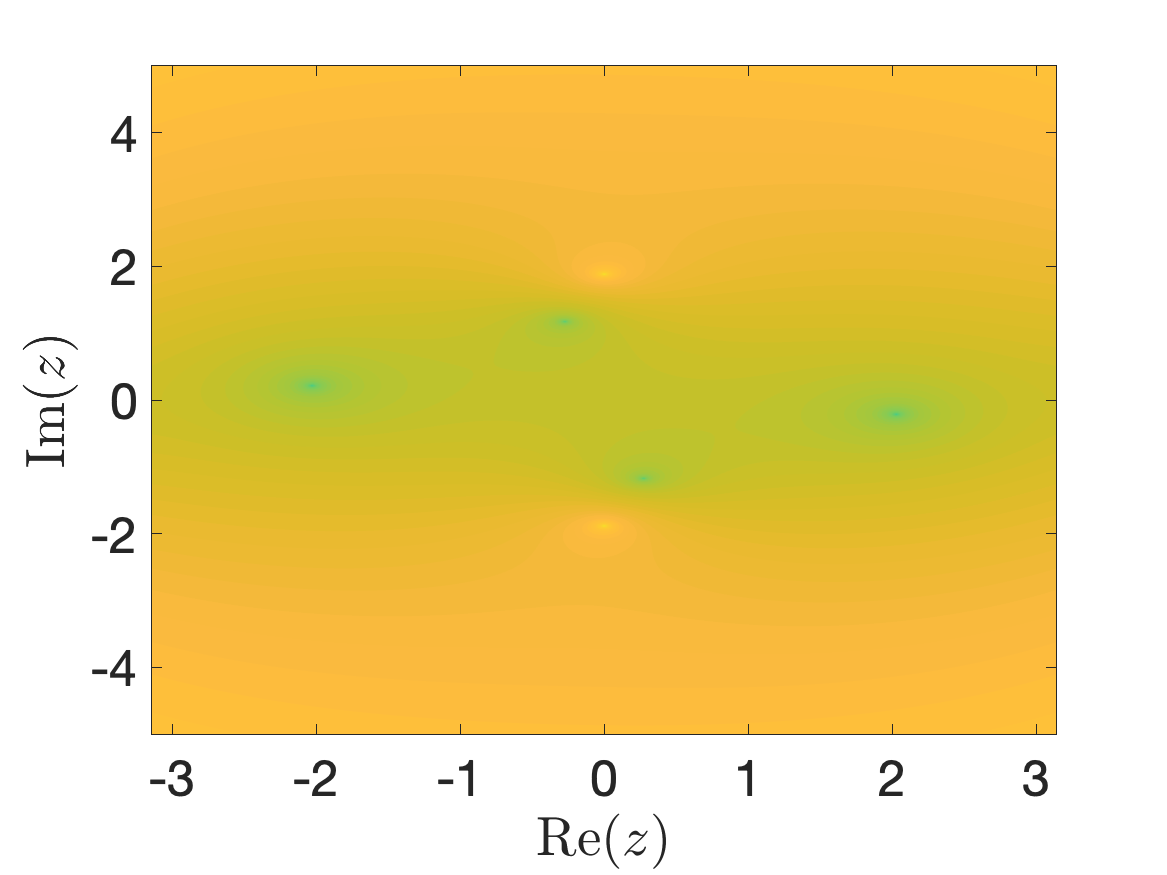}
    \caption{$t=-1$}
  \end{subfigure}
  \begin{subfigure}[b]{0.24\linewidth}
    \includegraphics[width=1.1\linewidth]{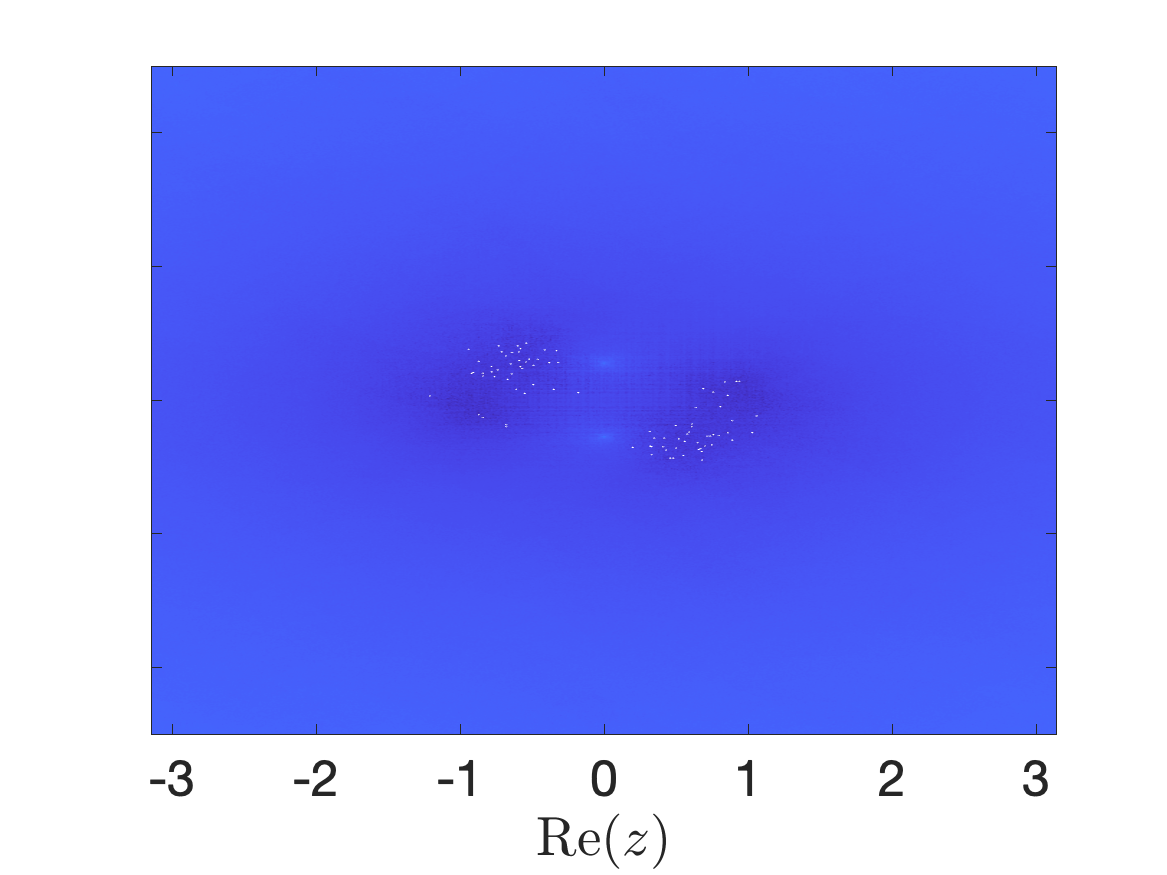}
    \caption{$t=0$}
  \end{subfigure}
  \begin{subfigure}[b]{0.24\linewidth}
    \includegraphics[width=1.1\linewidth]{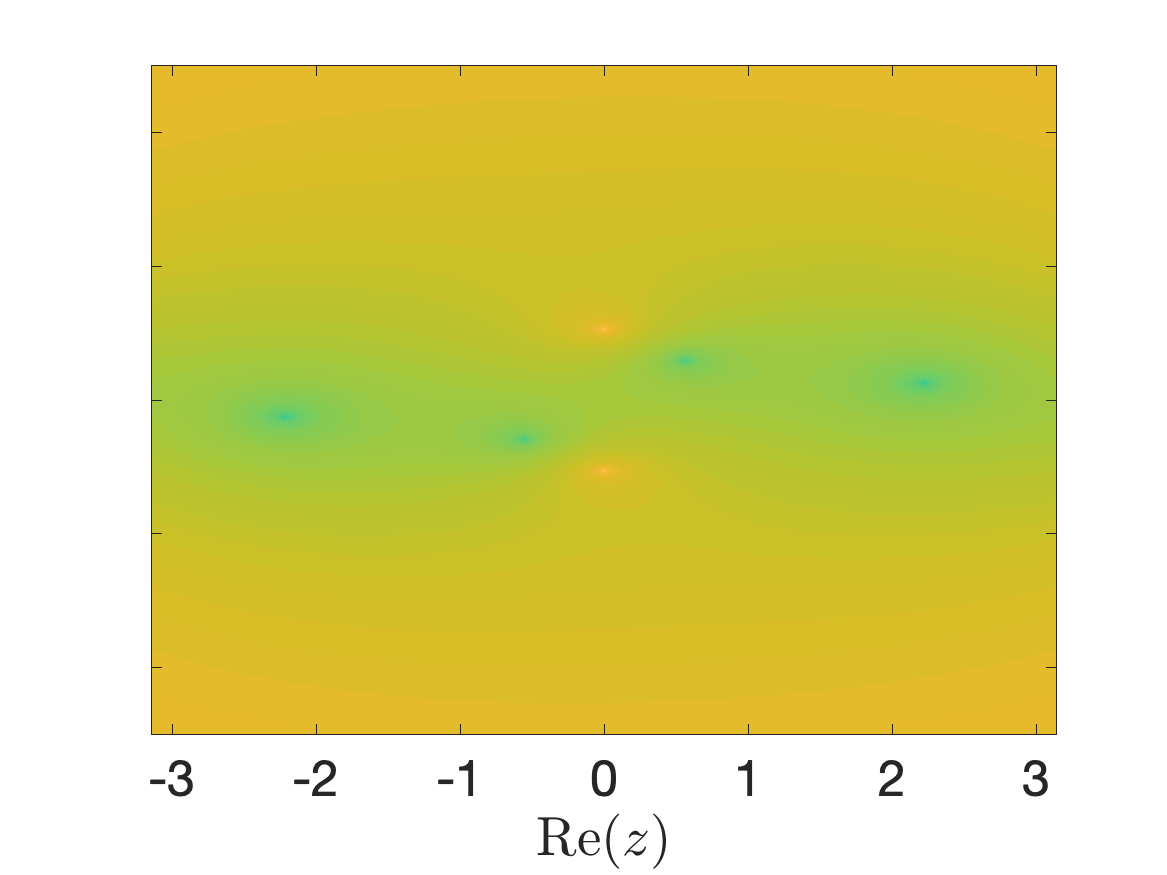}
    \caption{$t=0.5$}
  \end{subfigure}
  \begin{subfigure}[b]{0.24\linewidth}
    \includegraphics[width=1.1\linewidth]{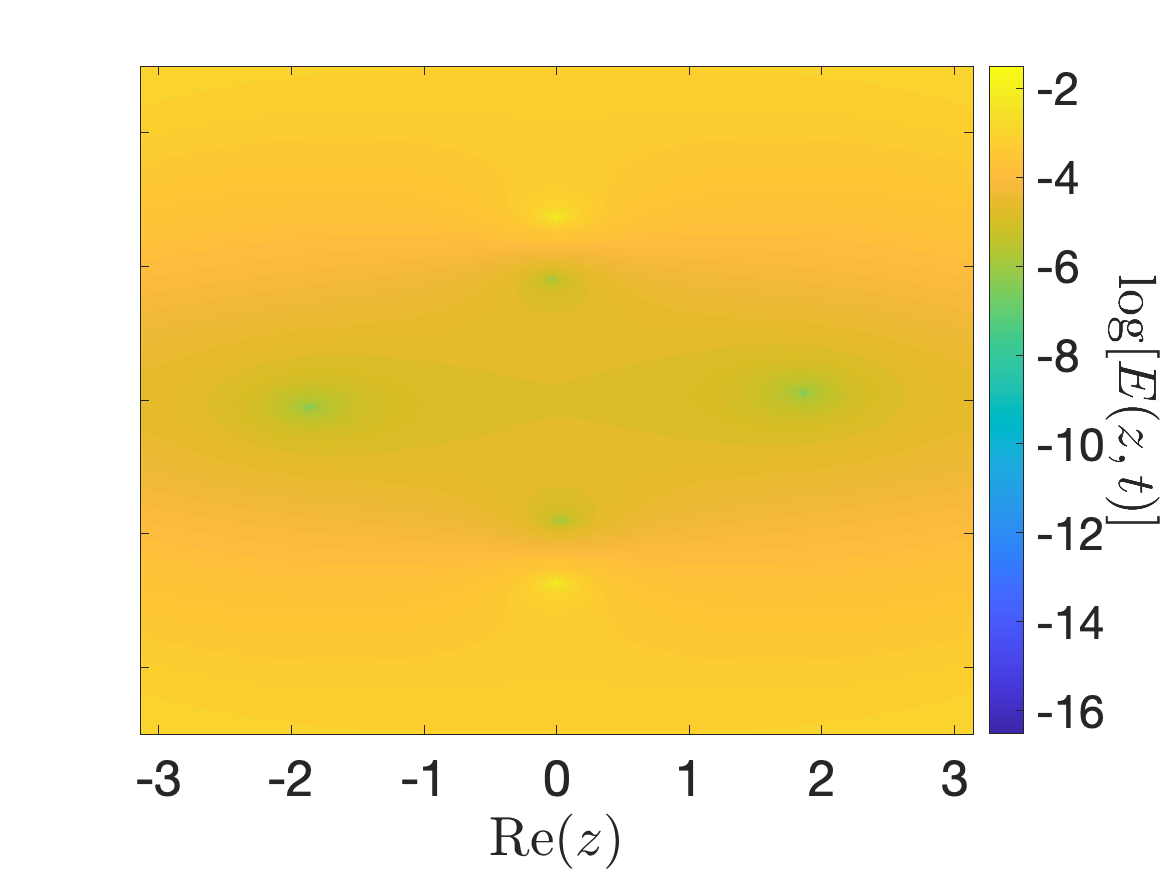}
    \caption{$t=1.5$}
  \end{subfigure}
  \caption{The pole-scaled approximation error $E(z,t)$ as in (\ref{errorrel}) in the logarithmic scale of the extended solution  of  the NLS (\ref{nls}) subject to the initial condition $u(x,0)$ as in  (\ref{solnls}) for time $t=-1, 0, 0.5, 1.5$. }
  \label{nlserror}
\end{figure}

Finally, we study  trajectories of the weights $w_{11}^{(\pm)}(t)$ and biases $b_{11}^{(\pm)}(t)$ of the hidden layers of the components $\Phi^{(\pm)}$ as functions of time in the complex plane. 
Note that, with only a single neuron in the hidden layer, there is no freedom to choose the parameters $C_{10}^{(\pm)}(t)$ randomly. Instead, for each choice of the pole $z_0$ of the activation functions, the parameters $w_{11}^{(\pm)}(t)$ and $b_{11}^{(\pm)}(t)$ are uniquely determined.
%Note that with a single neuron in the hidden layer, a random choice of the parameters $C_{j0}^{(\pm)}$ is not possible, so that $w_{11}^{(\pm)}(t)$ and $b_{11}^{(\pm)}(t)$ are uniquely determined for each choice of the pole $z_0$ of the activation functions.
If we choose $z_0$ to be real, the resulting weights satisfy $w_{11}^{(-)}(t)=w_{11}^{(+)}(t)<0$, while the biases satisfy $z_0+b_{11}^{(\pm)}(t)<0$, with $z_0+b_{11}^{(\pm)}(t)$ remaining constant in time (Figure \ref{fig_nls_w_b}).
%If we choose $z_0$ to be a real number, we obtain numerical values for the weights, such that $w_{11}^{(-)}(t) = w_{11}^{(+)}(t) < 0$, and for the biases, $z_0 + b_{11}^{(\pm)}(t) < 0$ and constant. 
Then, according to formula \eqref{spl}, as $w_{11}^{(-)}(t)$ decreases (increases) ($w_{11}^{(+)}(t)$ decreases (increases)), $\mathrm{Im}\big(s_1^{(-)}(t)\big)$ likewise decreases (increases) ($\mathrm{Im}\big(s_1^{(+)}(t)\big)$ increases (decreases)) (Figure \ref{fig_nls_w_b}). If, instead, the pole $z_0$ is chosen to be complex, the weights and biases become complex as well, but their overall behavior remains  unchanged (Figure \ref{fig_nls_w_b}).

\begin{figure}[h!]
    \centering  
     \begin{subfigure}[b]{0.45\linewidth}
    \includegraphics[width=0.8\linewidth]{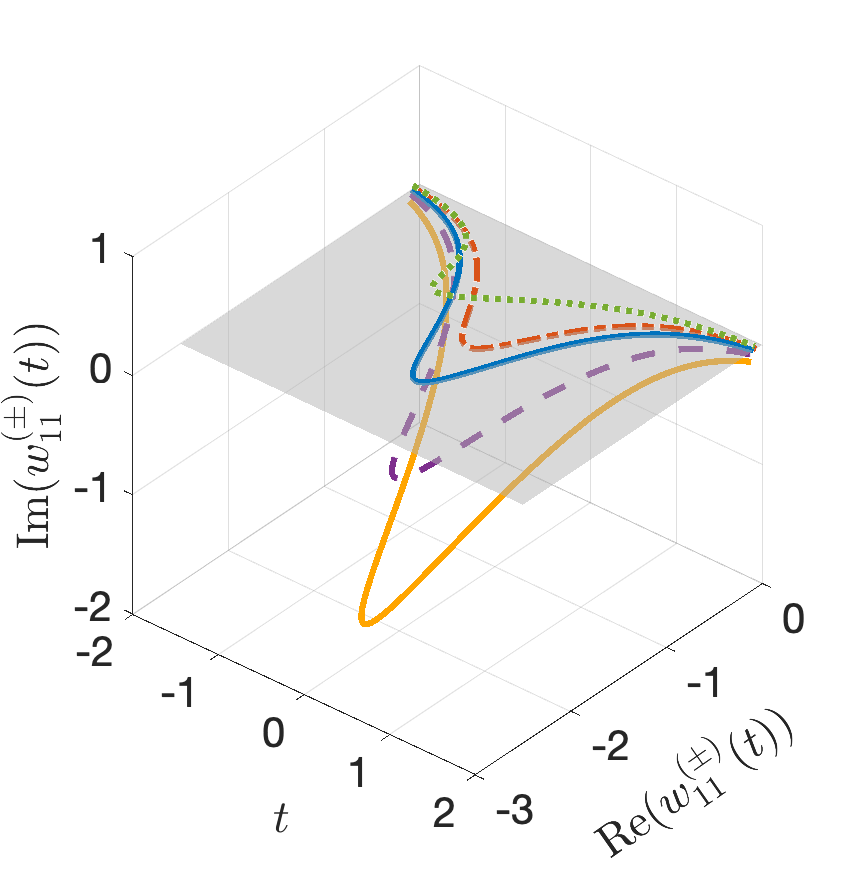}
   %\caption{$t=0.5$}
  \end{subfigure}
     \begin{subfigure}[b]{0.45\linewidth}
    \includegraphics[width=1.1\linewidth]{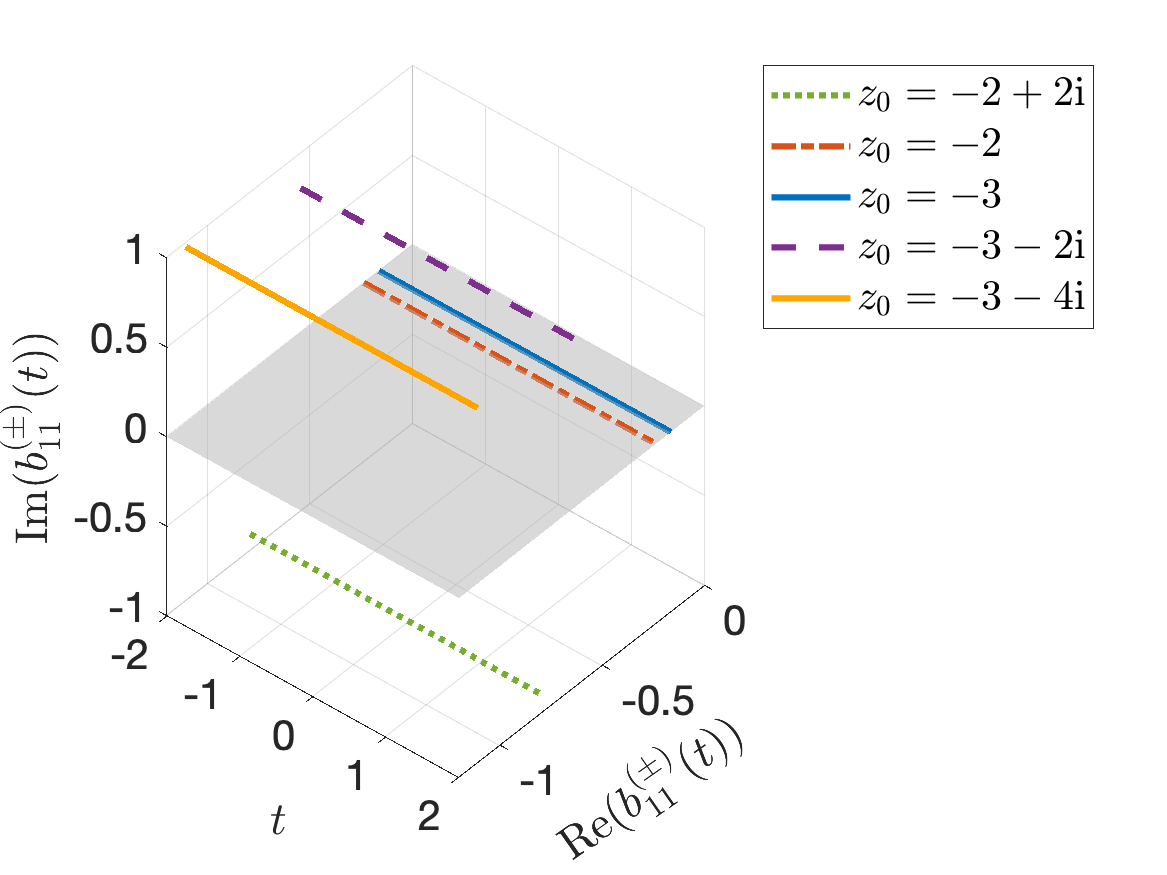}
  % \caption{$t=0.5$}
  \end{subfigure}
      \caption{ 
      Locations of the weights $w_{11}^{(\pm)}(t)$ (left) and biases $b_{11}^{(\pm)}(t)$ (right) used in \eqref{spl}--\eqref{spl1} for computing the singularities $s_1^{(\pm)}(t)$ of the extended solution of the NLS \eqref{nls} subject to the initial condition $u(x,0)$ as in \eqref{solnls} over the time interval $[-2,2]$ for different choices of the pole $z_0$ of the activation functions.  }
  \label{fig_nls_w_b}
\end{figure}

\end{example}

Within the proposed framework,
the emergence of rogue waves can be interpreted through the dynamics of complex singularities of the analytic continuation of the solution. As the nearest singularities approach the real axis, the region of analyticity around it narrows, producing a rapid, localized amplification of the wave amplitude.
 Unlike in the case of blow up or shock formation, the singularities do not reach the real axis. Instead, they recede from it after the peak amplitude is attained, and the solution returns to its background state.

\section{Conclusions and plans for  future work}

We present a neural network-based method for numerical analytic continuation and singularity detection, and apply it to study  phenomena arising in nonlinear PDEs. We consider three cases: (1) the nonlinear heat equation, which solutions exhibit finite-time blow up; (2) the nonlinear Burgers equation, which develops a shock; and (3) the Schrödinger equation, in which rogue waves emerge.
Solving the given PDE on the segment $[-\pi,\pi]$  in the form of a truncated Fourier series \eqref{appsr}, and using the corresponding Fourier coefficients $\hat{u}_k(t)$, we construct an analytic continuation $\Phi(z,t)$, as in \eqref{nnapr}, of the solution to the complex plane. Studying the dynamics of the complex singularities (poles, branch points, and branch cuts) of $\Phi(z,t)$  computed via formulas \eqref{spl}-\eqref{spl1}, helps to explain the occurence of these phenomena in the corresponding solutions.

In this paper, we focus on PDEs with univariate spatially periodic solutions. As in \cite{CG15}, our method can also be applied in the two-dimensional case to study the formation of singularities in Prandtl and Navier–Stokes wall shears (at different Reynolds numbers), in order to explain the various viscous–inviscid interactions characterizing the separation process in the corresponding solutions.
Another possibility is to study analogous problems in the non-periodic case. The method can be readily adapted by replacing the Fourier expansion \eqref{appsr} with a Taylor or Chebyshev expansion.
Yet another direction is to apply our method to the study of complex singularities of dynamical systems -- where time itself is treated as a parameter extended to the complex plane -- in order to investigate some of their properties (see, for example, \cite{W2011}).

Note again that the method for analytic continuation presented in this paper is based on our recent neural network-based approach for learning meromorphic functions of a single complex variable \cite{DKD2025}. Another promising direction is to extend the ideas developed in \cite{DKD2025} to the multivariate setting by employing the generalization of the SVD-based approach from \cite{GGT13} to several variables proposed in \cite{AK2021}, together with sparse grids \cite{BG24} and  sparse sampling \cite{DH25} techniques  to cope with the curse of dimensionality. We believe that these research problems deserve further investigation and leave them for future work.

\section*{Statements and Declarations}
\textbf{Competing Interests:} The authors declare that they have no competing interests.

%\section*{Acknowledgement}   

%F.D. is supported by \dots

\small
%\bibliography{myBib}{}
%\bibliographystyle{plain}

%\newpage

%\appendix
%\section{Appendix  }

\end{document}